\documentclass[English]{article}
\usepackage[T1]{fontenc}
\usepackage[latin9]{inputenc}
\usepackage{amsmath}
\usepackage{amssymb}
\usepackage{amsthm}
\usepackage{graphicx}
\usepackage{tikz}
\usepackage[all]{xy}
\usepackage{bm}
\usepackage{comment}
\usetikzlibrary{calc,patterns,angles,quotes}
\usepackage[toc,page]{appendix}
\usepackage{hyperref}

\newtheorem{prop}{Proposition}

\newcommand{\R}{\mathbb{R}}

\newcommand{\C}{\mathbb{C}}

\newcommand{\arccosh}{\operatorname{arccosh}}
\newcommand{\sign}{\operatorname{sign}}

\title{{Boltzmann's Billiard Systems: Computation of the Billiard Mapping and Some Numerical Results}}
\author{Michael Plum, Airi Takeuchi, Lei Zhao}

\begin{document}
	%\graphicspath{{Fig_Numerics/}}
	\date{}
	\maketitle	
	\begin{abstract}
		{L. Boltzmann proposed in \cite{Boltzmann} a billiard model with a planar central force problem reflected against a line not passing through the center.  He asserted that such a system is ergodic, which thus illustrates his ergodic hypothesis. However, it {has been} recently shown that when the underlying central force problem is the Kepler problem, then the system is actually integrable \cite{Gallavotti-Jauslin}. This raises the question {of} whether Boltzmann's assertion holds true for \emph{some} central force problems that he considered. In this paper, we present some geometrical and numerical analysis on the dynamics of several of these systems. As indicated by the numerics, many of these systems show chaotic dynamics and a system seems to be ergodic. }
		\end{abstract}
	
	\section{Introduction}

	In \cite{Boltzmann}, Boltzmann examined a billiard model in a potential field, which is defined with the central force problem in $\R^2$ with the potential 
	$$V_{\alpha,\beta}:= -\dfrac{\alpha}{2r} + \dfrac{\beta}{2 r^{2}},$$ 
	in which $r$ is the distance of the moving particle to the origin $O$ and $\alpha,\beta \in \R$ are parameters. The motion of the moving particle is assumed to reflect elastically against a line in $\mathbb{R}^{2}$ with distance $\gamma>0$ to the origin $O$. Physically, the potential $V_{\alpha, \beta}$ describes a Kepler-Coulomb problem with an additional centrifugal correction. When $\alpha > 0, \beta =0$, the potential $ V_{\alpha,\beta}$ {is that of }the attractive Kepler problem. 
		%In such a case, the billiard systems in $y \leq \gamma$ and $y \geq \gamma$ are equivalent by mirroring each reflected arc against the wall at $y = \gamma$.
	
	Boltzmann considered this as a simple model which illustrates his ergodic hypothesis. When the energy of the system is properly fixed {which ensures that the orbits consecutively hit the line of reflection}, he asserted that
	\begin{itemize}
		\item the billiard mapping of the system preserves a measure, and
		\item The dynamics are ergodic with respect to this measure.
	\end{itemize}
	
	Boltzmann explicitly computed the billiard mapping and its Jacobian to establish the first assertion. The computation was actually incomplete. We discuss this issue in Section \ref{sec: Boltzmann_billiard_mapping}. This assertion is nevertheless true, as follows from these computations. Nowadays we understand that this is a more general feature related to symplectic reduction{, which is discussed} in Section \ref{sec: symplectic}.
	
	The second assertion of Boltzmann has been proven false when $\beta = 0$, which is actually integrable. In \cite{Gallavotti}, Gallavotti suspected that in this scenario the system is actually integrable based on numerical evidence. In \cite{Gallavotti-Jauslin}, Gallavotti and Jauslin explicitly constructed a conserved quantity in addition to the total energy of the system, which proves its integrability. The integrable behavior of the system is analyzed by Felder \cite{Felder}, who shows the Poncelet property of the system. {Two alternative proofs of the integrability of the system are provided in \cite{Zhao2021} (which was extended to more general systems in \cite{TakeuchiZhao2}) and \cite{Takeuchi-Zhao1}.} Moreover, the analysis in \cite{Felder} shows that KAM stability holds for systems with $(\alpha, \beta)$ such that $|\beta/\alpha|$ is sufficiently small. Therefore, for Boltzmann's ergodic assertion to hold, the parameter $\beta$ must have a larger norm compared to $\alpha$.
	
        It is an open question to determine whether Boltzmann's billiard system is ergodic for some $(\alpha, \beta)$. In this paper, we conduct numerical simulations with different parameters. The numerical results demonstrate a diverse range of dynamical behavior, and suggest that the system might indeed be ergodic for some parameters.

	\section{{Canonical Coordinates} for Central Force Problem}
	\label{sec: cnst_coodinate}
	In this section, we construct some canonical coordinates for general central force problem in the plane. {For Kepler problem, these type of coordinates can be already found in the work of Lagrange \cite{LagrangeCoordinates}. }
	%In the following, we describe the construction of the canonical coordinates, as is explained in \cite{Valtonen} for the Kepler problem, but here for more general central force problems. 
	
	%They are also considered as a generalization of the canonical coordinates associated to the Kepler problem in the plane explained in 
	
	Consider a central problem in the plane with a general radial force function $U = U(r)$. The potential is $V(r)=-U(r)$.
	The kinetic energy is given {in polar coordinates} by 
	\[
	K =\frac{1}{2}(\dot{r}^2+ r^2 \dot{\phi}^2). 
	\]
       The conjugate momenta are 
	\[
	p_{r}= \frac{\partial K}{\partial \dot{r}}=\dot{r}, \quad p_{\phi}= \frac{\partial K}{\partial \dot{\phi}}=C,
	\]
	where $C := r^2 \dot{\phi}$ is the angular momentum. 
	
	In the coordinates $(r, \phi, p_{r}, p_{\phi})$, the Hamiltonian is given by 
	\[
	H =\frac{1}{2}\left(p_{r}^2+ \dfrac{p_{\phi}^2}{r^2}\right) - U(r).
	\]
	
        The coordinates $(Q_1, Q_2, P_1,P_2)$ will first be constructed via the generating function $S= S(r, \phi, P_1,P_2).$
	
	The time-independent Hamilton-Jacobi equation $$H\Bigl(r, \phi, \frac{\partial S}{\partial \phi},\frac{\partial S}{\partial r} \Bigr) =  E,$$ takes the form
	\begin{equation}
		\label{eq: H-J}
		\frac{1}{2}\left( \left(\frac{\partial S}{\partial r}\right)^2  +\left( \frac{1}{r}\cdot \frac{\partial S}{\partial \phi} \right)^2  \right) - U(r)   = E.
	\end{equation}
        Assume that the function $S$ is {separated} into
	\[
	S(r, \phi)= S_r(r) + S_{\phi}(\phi).
	\]
	Substituting into the equation \eqref{eq: H-J}, we get
	\begin{equation}
		\label{eq: S_1}
		\frac{1}{2}\left( \left(\frac{\partial S_r}{\partial r}\right)^2  +\left( \frac{1}{r}\cdot \frac{\partial S_{\phi}}{\partial \phi} \right)^2  \right) - U(r)  =E.
	\end{equation}
	
	\begin{comment}
		
		Since the LHS is independent of $t$ and RHS is independent of $r$ and $\phi$, both hand sides must be a constant that is independent of $r$, $\phi$, and $t$.
		Hence we have 
		\begin{equation}
			\label{eq: S_1}
			\frac{\partial S_t}{\partial t}= -\kappa_1
		\end{equation}
		and 
		\begin{equation}
			\label{eq: S_2}
			\frac{1}{2}\left( \left(\frac{\partial S_r}{\partial r}\right)^2  +\left( \frac{1}{r}\cdot \frac{\partial S_{\phi}}{\partial \phi} \right)^2  \right) - V(r) = \kappa_1,
		\end{equation}
		where $\kappa_1$ is constant.

	\end{comment}

	In the above equation, only  $\frac{\partial S_{\phi}}{\partial \phi}  $  in the LHS is dependent on $\phi$, thus can be set as a constant. We write, 
	\begin{equation}
		\label{eq: S_3}
		\frac{\partial S_{\phi}}{\partial \phi}  = \kappa
	\end{equation} 
	On the other hand, since $S$ is the generating function, we have
	\[
	\frac{\partial S}{\partial \phi} = p_\phi=C.
	\]
	Thus, we have $\kappa= C$.
	
	Substituting this into the equation \eqref{eq: S_1} and assuming {$\dot{r}(=\partial S/ \partial r)\geq0$},  we obtain 
	\begin{equation}
		\label{eq: S_4}
		\frac{\partial S_r}{\partial r}=  \sqrt{2( E+U(r))-\frac{C^{2}}{r^{2}}}.
	\end{equation}
	
		%From the equations \eqref{eq: S_3},\eqref{eq: S_4}, the 
		
		Thus we may choose 
			\[
	S = C \phi +\int_{r_{min}}^r \sqrt{2(E+U(r))-\frac{C^{2}}{r^{2}}}.
	\]
	{We set the two constants of integration $E$ and $C$ as the new momenta}
	\[
	P_1 = E, \quad P_2 = C.
	\]
	
	%\rd{We get}
	%\[
	%Q_1 = \frac{\partial S}{\partial P_1} = \frac{\partial S}{\partial E} , \quad Q_2 = \frac{\partial S}{\partial P_2} = \frac{\partial S}{\partial C}.
	%\]

	{We get} 
	%$Q_1$ and $Q_2$ can be computed as
	\begin{equation}
		\label{eq: new_coordinate1}
		Q_1 {= \frac{\partial S}{\partial E}}=  \int_{r_{min}}^r \frac{dr}{\sqrt{2({E} + U(r)) -  \frac{C^2}{r^2}} } 
	\end{equation}%\marginpar{\rd{Check your thesis}}
	and 
	\begin{equation}
		\label{eq: new_coordinate2}
		Q_2 {=\frac{\partial S}{\partial C}}= \phi - C\int_{r_{min}}^r \frac{dr }{r^2 \sqrt{2({E} + U(r)) - \frac{C^2}{r^2}}}.
	\end{equation}
	
	 The conservations of the energy and the angular momentum can be presented as
	\[
	\dot{r}^2+r^2 \dot{\phi}^2 = 2 \cdot {E} + U(r),
	\]
	\[
	r^2 \dot{\phi} = C.
	\]
	With the assumption $\dot{r}>0$, we obtain
	%(see also Section \ref{subsec: solution_curves} for the derivations)
	\[
	dt = \frac{dr}{\sqrt{2(E + U(r)) -  \frac{C^2}{r^2}} },
	\]
	and 
	\[
	d\phi = \frac{C\, dr }{r^2 \sqrt{2(E + U(r)) -  \frac{C^2}{r^2}}}.
	\]
	Thus from Equations \eqref{eq: new_coordinate1} and \eqref{eq: new_coordinate2} we have
	\[
	Q_1 =  \int_{t_{peri}}^t dt = t - t_0 = \tilde{t}
	\]
	\[
	Q_2 = \phi - \int_{g}^{\phi} d\phi = g,
	\]
	where $t_{peri}$ represents the time of the pericenter passage from a fixed direction and $g$ represents the argument of the pericenter (the angle of the pericenter from the first coordinate direction).
	%By taking the time derivative of these components, we obtain, 
	%\[
	%\dot{Q_1}=1, \quad \dot{Q_2}= 0.
	%\]
	The case $\dot{r} < 0$ can be treated similarly {and the meaning of the variables retains. We should nevertheless remember that in the case of multi-pericenters, the angle $g$ is assigned to a fixed one and is subject to a choice. }
	%The physical meaning of these new coordinates is as follows: $Q_1$ represents the time $\tilde{t}$ of the passage from the pericenter, and $Q_2$ represents the argument of pericenter $g$ (the angle of the pericenter from a fixed direction). %For the direct computation of the integrals appeared in \eqref{eq: new_coordinate1} and \eqref{eq: new_coordinate2}, see  \cite[Section 4.7]{Valtonen}.

	In this way, we obtain the canonical coordinates 
	$${(P_{1}, P_{2}, Q_{1}, Q_{2}):=(E, C, \tilde{t}, g)}.$$

	\section{Symplectic Property of the Billiard Mapping}
	\label{sec: symplectic}
	%\rd{Rewrite this part.}
	
	{In this section, we discuss the billiard mapping from the viewpoint of symplectic geometry.} A main assertion by Boltzmann in \cite{Boltzmann} is that the billiard mapping preserves a measure. {This is easily deduced from the fact that} the billiard mapping preserves {an}
	 explicit symplectic 2-form, {and therefore preserves the associated Liouville measure.} {In the last part of this paper}, Boltzmann remarked that this preservation holds for more general force function $U= U(r)$, and with any curve $\mathcal{C}: r= \psi (\theta)$. In our discussion, we assume that $U$ {is $C^{1,1}$} and $\mathcal{C}$ is of class $C^1$, {so that the reflection law is well-defined.} 
	
	%In the following, we explain the method of symplectic reduction.
	%which will be use to establish the symplecticty for the new coordinates used for the billiard mapping in the latter argument. 
	%We start with Hamiltonian symplectic geometry. Recall the followings:  
	%\marginpar{\rd{These sentences are not informative.}}
	A symplectic manifold is a pair $(M, \omega)$ with {$M$ a smooth manifold} with a closed, non-degenerate 2-form $\omega$. A vector field $X$ on $M$ is called a Hamiltonian vector field with the Hamiltonian $H$ which we assume to be of class {$C^{1,1}$},  if there holds
	$$\omega(X, \cdot )= -d H.$$
	
	 In a natural mechanical system such as our central force problems, the symplectic manifold
	is the cotangent bundle of the configuration space, equipped with
	a canonical symplectic form, and the Hamiltonian function is the total
	energy. 
	
	%Now we take quotients of symplectic manifolds under group actions that is called symplectic reduction. 
	
	Let $G$ be a Lie group acting on a symplectic manifold $(M, \omega)$. We say that the action is Hamiltonian if for every $\xi \in T_e G${, the} associated vector field $X_\xi$ given by 
	\[
	(X_\xi )_x =  \left. \frac{d}{dt} \right| _{t =0} \exp(t \xi ) \cdot x 
	\]
	is { a Hamiltonian vector field, the Hamiltonian of which is denoted by $H_\xi$. Let $c$ be a regular value of $H_\xi$. The level set $H^{-1}_\xi(c)$ is thus a codimension-1 submanifold of $M$ on which $G_\xi := \{\exp(t\xi)  \}$ acts freely in the kernel direction of the restriction of $\omega$ to $H^{-1}_\xi(c)$. The quotient space $H^{-1}_\xi(c)/ G_\xi$ is thus again symplectic.} %This can be checked from the facts that $\omega $ is invariant under $G$ and $H_{\xi}$ is a constant of motion (\emph{i.e.} $\{H, H_\xi  \}=0$ ).
	
	{For our} central force problem, the Hamiltonian {is}
	\[ 
	H(p,q)=\dfrac{\|p\|^{2}}{2}-U(r),\quad(p,q)\in\mathbb{R}^{2}\times(\mathbb{R}^{2}\setminus O),r=\|q\|.
	\]
	{The} canonical symplectic form is 
	\[
	\omega=dp_{1}\wedge dq_{1}+dp_{2}\wedge dq_{2}.
	\]
	%This symplectic form is invariant under $SO(2)$ action. 
	{The Hamiltonian vector field of the conserved angular momentum $C$ generates an $SO(2)$-symmetry of the system by simultaneously rotating $p, q$.} We may thus apply the {above} symplectic reduction {procedure} and {obtain} the reduced Hamiltonian 
	\[
	H_{r}(p,q;C)=\dfrac{p_{r}^{2}}{2}-U(r)+\dfrac{C^{2}}{r^{2}},
	\]
	with the reduced symplectic form {$ d p_r \wedge d r$}. The reduced energy
	level $\left\{ H_{r}=h\right\} $ projects into the Hill's region
	$\left\{ -U(r)+\dfrac{C^{2}}{r^{2}}\le h\right\} $ in the reduced
	configuration space {$\mathbb{R}_{+}:=\{s \in \R; s>0\}$}. We assume that this projection
	is not the full $\mathbb{R_{+}}$ and we consider a connected component
	of this projection that is not merely a point: then it is a closed interval
	{$[a_{h},b_{h}]$ with $a_{h}>0$ and $b_{h} \le \infty$.}
	%(Analogously we may consider the case that the interval has a boundary point $0$ and the other boundary point is finite.)
	%\marginpar{\rd{I am not sure with the commented sentence. Have a check please.}}
	 We assume in addition
	that the boundary points of this component depend continuously on
	$h$ and $C$. We localize our system near this component for energy
	close to $h$ and for angular momentum close to $C$. In this way, we get 
	the localized system defined on the localized phase space, that preserves the symplectic structure of the original system.
	
	{After this localization} we use the coordinates $(H, C, \tilde{t}, g)$ {constructed in Section \ref{sec: cnst_coodinate}}. Note that for general $U(r)$ and for a fixed orbit, pericenters and apocenters are not unique. Therefore this canonical variables {$(E, C, \tilde{t}, g)$} is defined on a covering space on the localized phase space, {which is naturally also symplectic. }
	%We remember that the localized phase space is a symplectic therefore its covering space is also symplectic. 
	{In these coordinates, the symplectic form is}
	\[
	{d E \wedge d\tilde{t}+dC\wedge dg.}
	\]

	{Here is another way to deduce these coordinates}. Since $\partial_{\tilde{t}},\partial_{g}$ are respectively
	the Hamiltonian vector fields of {$E(=H)$} and $C$, the symplectic form
	has to take the form 
	\[
	{d E} \wedge d\tilde{t}+dC\wedge dg+f({E},C,\tilde{t},g)dH\wedge dC;
	\]
	now since $\left\{ E,C\right\} =0$ we conclude that $f(H,C,\tilde{t},g)=0$.
	%Consequently, the reduced symplectic form by fixing $H$ is simply
	%$dC\wedge dg$. 
	{Fixing $E$ and passing through the quotient by time, we get the reduced symplectic form $dC\wedge dg$ on the space of orbits.}
	
	We now add the wall of reflection. We study the reflection at a point $x_{0}$
	in $\mathcal{C}$. Without loss of generality, we may {put} this point at the origin {with }
	the tangent line {$T_{x_{0}} \mathcal{C}$ being} the first coordinate axis. {The reflection is given by the involution} $(p_{1},p_{2},q_{1},q_{2})\mapsto(p_{1},-p_{2},q_{1},-q_{2})$ {which is symplectic. }
	
	Consequently, {the symplectic form is preserved in the billiard system}. Since {$E$} is invariant along
	the orbit {as well as at reflections at $\mathcal{C}$}, we conclude that the reduced form $dC\wedge dg$ is preserved under {the reflections}. 
	
	Note that {this} argument holds only when ${E}$ and $C$ are functionally independent.
	{Otherwise, the} orbit is circular {and the angle $g$ is not well-defined. It is nevertheless not hard to see that this corresponds to a set of measure zero in the phase space. }
	
	{The first assertion of Boltzmann follows}. 
	
	%On the other hand, for almost circular motions we may take a set of
	%coordinates following a trick of Poincar\'e. We shall complete the details
	%when necessary at a later stage. 

	\section{Computation of the Billiard Mapping}
	\label{sec: Boltzmann_billiard_mapping}
	\subsection{Solutions for the Central Force Problem: Kepler Problem with Centrifugal Force}
	\label{subsec: solution_curves}
	%\subsection{Solutions in the case $\beta =0$: Kepler Problem}
	%\subsection{Solutions in more general case $\beta \neq 0$}
	As the first step of computing billiard maps, we solve the central force problem in the plane with the force function $U= \dfrac{\alpha}{2r} - \dfrac{\beta}{2 r^{2}}$.
	
	We {first} recall the analysis of Boltzmann in {\cite{Boltzmann}}. He $(r,\phi)$ and {{wrote down} the equations {on} the preservation of the energy $E$ and the angular momentum $C$ {in polar coordinates $(r,\phi)$}:
		%\footnote{{Normally I take the convention to use $C$ for the angular momentum. }} 
		}

	\begin{align}
		\label{eq: first_der}
		\dot{r}^{2}+r^{2}\dot{\phi}^{2}&=2\cdot E+\dfrac{\alpha}{r}-\dfrac{\beta}{r^{2}},\\
		\label{eq: anglmomenta}
		r^{2}\dot{\phi}&=C,
	\end{align}
	
	in which a dot denotes the derivative with respect to
	time.
	
	Note that here $E$ is always conserved in the billiard system {since the kinetic energy does not change at reflections}, while
	$C$ changes from orbit arcs to orbit arcs when a reflection at the
	wall takes place. {We shall only consider bounded orbits, so we set $E<0$.}
	
	Boltzmann then writes ``from here it follows that''
	
	\[
	\dot{r}=\sqrt{2\cdot E+\frac{\alpha}{r}-\frac{C^{2}+\beta}{r^{2}}},
	\]
	
	thus
	
	\[
	dt=\frac{dr}{\sqrt{2\cdot E+\frac{\alpha}{r}-\frac{C^{2}+\beta}{r^{2}}}}.
	\]

	This deduction is problematic, {as in general $x^{2}=a$ does not
		imply $x=\sqrt{a}$}.
		% as well as the computations {that follow}. 
		%\marginpar{\rd{I deleted a small sentence. Please have a check.}}
		Along
	an arc containing either an pericenter or an apocenter, the quantity
	$\dot{r}$ changes its signs. 
	
	From (\ref{eq: anglmomenta}), it follows that
	
	\[
	dt=\frac{r^2 \, d\phi} {C}.
	\]
	
	By equating these equations for $dt$, we have

	\[
	d\,\phi=\dfrac{C\,dr}{r\sqrt{2\cdot Er^{2}+\alpha r-C^{2}-\beta}}.
	\]
	
	Also, this formula 
	is problematic, as it uses the previous formula. Indeed, the LHS
	has the same sign as $C$, which is positive resp. negative when the
	corresponding orbit is oriented counterclockwise resp. clockwise.
	On the other hand, when an arc contains a peri- or apo-center, the
	monotonicity of $r$ changes while the monotonicity of $\phi$ does
	not change.
	
	We now restrict our system to the case that these formulas are valid. Namely, we consider an arc between a pericenter and the consecutive apocenter. In this case, one can rewrite the above equation as

	\[
	d\,\phi=\dfrac{dr/ r^2}{\sqrt{\frac{2 E}{C^2} + \frac{\alpha}{C^2 r} - \frac{C^2 + \beta}{C^2r^2}}} 
	=\sqrt{\frac{C^2}{C^2 + \beta }} \cdot    \frac{dr/ r^2}{\sqrt{-(\frac{1}{r_{min}} - \frac{1}{r})(\frac{1}{r_{max}} - \frac{1}{r})}},
	\]
	assuming that $C^2+ \beta >0$.
	Here, $r_{min}$ and $r_{max}$ are {respectively distances of the pericenter and apocenter to the center of the system. We have
		$$r_{min} = \frac{-\alpha + \sqrt{\alpha^2 + 8 E (C^2 + \beta )}}{4E},$$
		$$r_{max} = \frac{-\alpha - \sqrt{\alpha^2 + 8 E (C^2 + \beta )}}{4E} .$$ }
	
	{We now set} $\rho = 1/  r$ so that {$d \rho = -dr / r^2  $}
	% \footnote{{There was a minus sign missing}}
	, and we get 
	
	\[
	d\,\phi
	= \sqrt{\frac{C^2}{C^2 + \beta }} \cdot  \frac{-d \rho }{\sqrt{-(\rho_{min}- \rho)(\rho_{max} - \rho)}},
	\]
	
	naturally, $\rho_{min} = 1/ r_{min}$ and $\rho_{max} = 1 / r_{max}$. Note that $\rho_{min} \geq \rho_{max} $. 
	
	Finally, we change the integration variable {to} $\chi = \rho - \frac{1}{2}(\rho_{max} - \rho_{min})$ and set $\chi_0 = \frac{1}{2} (\rho_{min} - \rho_{max} )$.
	%\footnote{{Normally $\widetilde{\chi}$ denotes a variant of $\chi$, so the notation here may have to be modified. I think $\chi_{0}$ is a good choice} } 
	{We have}
	
	\[
	d\,\phi
	= \sqrt{\frac{C^2}{C^2 + \beta }} \cdot  \frac{- d\chi}{\sqrt{\chi_0^2- \chi^2}}.
	\]
	
	If $C \geq 0$ so that the particle moves in the counterclockwise direction, then the integration from the pericenter to another point {on the orbit arc between the pericenter and the successive apocenter} becomes
	
	\begin{align*}
		\phi - \varepsilon
		&= \int_{r_{min}}^{r} \dfrac{C\,dr}{r\sqrt{2\cdot Er^{2}+\alpha r-C^{2}-\beta}} \\
		&= \sqrt{\frac{C^2}{C^2 + \beta }} \cdot  \int_{\chi_0}^{\chi} \frac{- d\chi}{\sqrt{\chi_0^2- \chi^2}} \\
		&=  \sqrt{\frac{C^2}{C^2 + \beta }} \cdot \arccos \frac{\chi}{\chi_0} \\
		&= \sqrt{\frac{C^2}{C^2 + \beta }} \cdot \arccos \dfrac{2\frac{C^{2}+\beta}{r}-\alpha}{\sqrt{{\normalcolor \textcolor{black}{\alpha}}^{2}+8\cdot E(C^{2}+\beta)}},
	\end{align*}
	where $\varepsilon$ denotes {the argument of a pericenter, {\emph{i.e.}} }the angle that the pericenter makes {from} the x-axis. In the third equality we used $\arccos 1 = 0 $.
	Recall that $0 \leq \arccos x \leq \pi $ for $-1 \leq x  \leq 1 $.
	
	For $C \leq 0$, the particle moves in the clockwise direction. Considering that the sign of the left hand side will be changed for this case, we get 
	\[ 
	\phi - \varepsilon
	= - \sqrt{\frac{C^2}{C^2 + \beta }} \cdot \arccos \dfrac{2\frac{C^{2}+\beta}{r}-\alpha}{\sqrt{{\normalcolor \textcolor{black}{\alpha}}^{2}+8\cdot E(C^{2}+\beta)}}.
	\]
	
	{These two cases} can be rewritten in a uniform way as
	
	\begin{equation}
		\label{eq:angle_shift_1}
		\phi - \varepsilon
		=  {\frac{C}{\sqrt{C^2 + \beta} }} \cdot \arccos \dfrac{2\frac{C^{2}+\beta}{r}-\alpha}{\sqrt{{\normalcolor {\alpha}}^{2}+8\cdot E(C^{2}+\beta)}}.
	\end{equation}
	{This equation {appeared} in Boltzmann's paper \cite{Boltzmann}. {However mind the typo therein.}

		%It can be easily checked that when $\phi = \varepsilon$ the particle is placed at the perihelion and when $ \phi = \frac{C^2}{C^2 + \beta } \pi + \varepsilon$ it is placed at the aphelion. \footnote{\rd{I think this can be deleted. The pericenter part of the information has been stated before. The apocenter part can be recalled when needed.}}

		%It is verified (by Maple) that $\phi=\varepsilon\Rightarrow\dot{r}=0$,
		%and in our mechanical system that the particle is placed at a peri-center. 
		{To complete his study, }we still have to consider the case 
		
		\[
		\dot{r}= -\sqrt{2\cdot E+\frac{\alpha}{r}-\frac{C^{2}+\beta}{r^{2}}},
		\]
		when the distance from the center decreases with time. {Boltzmann did not consider this case, making his analysis incomplete.}
		
		In this case, the sign of {the LHS} of (\ref{eq:angle_shift_1}) {needs to be} switched, and we have
		\begin{equation}
			\label{eq:angle_shift_2}
			\phi - \varepsilon
			=  -{\frac{C}{\sqrt{C^2 + \beta} }} \cdot \arccos \dfrac{2\frac{C^{2}+\beta}{r}-\alpha}{\sqrt{{\normalcolor \textcolor{black}{\alpha}}^{2}+8\cdot E(C^{2}+\beta)}}.
		\end{equation}
		
		% \footnote{\rd{It can be a good idea to include the figures that you draw on the arcs containing the pericenter/apocenter etc here.}}
		
		%Besides this sign change, there is another problem for general cases $\beta \neq 0$. In the case of pure-Keplerian motions (\rd{$\alpha<0, \beta  = 0$}), the angle $\phi - \varepsilon$ lies between $[-\pi, \pi]$, hence this parametrization is actually okay, as well as for those systems such that the perihelion proceeds directly, so that the central angle between two consecutive perihelion is at least $2 \pi$. When the perihelion proceeds retrogradely, the central angle between two consecutive perihelion is below $2 \pi$, and such a parametrization may not cover all the orbit arc.
		
		%Overall there are many problems with the deduction of Boltzmann. 
		
		We may then solve the problem further from (\ref{eq:angle_shift_1}) and (\ref{eq:angle_shift_2}).
		%It can be preferred to work this problem out in another way. 
		To be consistent with modern convention in celestial mechanics, we denote the angle which the particle makes {from the} x-axis by $\theta$ and denote the angle of (one of) the pericenter makes {from the} x-axis by $g$.
		%%%\footnote{\rd{A standard notation for argument of pericenter is $g$.}}
		
		From 
		\[
		\theta - g
		=  \pm{\frac{C}{\sqrt{C^2 + \beta} }} \cdot \arccos \dfrac{2\frac{C^{2}+\beta}{r}-\alpha}{\sqrt{{\normalcolor \textcolor{black}{\alpha}}^{2}+8\cdot E(C^{2}+\beta)}},
		\]
		
		we get 
		\[
		\pm{\frac{\sqrt{C^2 + \beta} }{C}}  (\theta - g)
		=  \arccos \dfrac{2\frac{C^{2}+\beta}{r}-\alpha}{\sqrt{{\normalcolor \textcolor{black}{\alpha}}^{2}+8\cdot E(C^{2}+\beta)}}.
		\]
		By taking cosine in both sides {we get}
		\begin{equation}
			\label{eq: solution_trans}
			\cos \left( \sqrt{\frac{{C^2 + \beta} }{C^2}}  (\theta - g) \right)
			=  \dfrac{2\frac{C^{2}+\beta}{r}-\alpha}{\sqrt{{\normalcolor \textcolor{black}{\alpha}}^{2}+8\cdot E(C^{2}+\beta)}}.
		\end{equation}
		Solving this equation for $r$ in the case of $\alpha >0$, we get 
		\begin{equation}
			\label{eq: solution1}
			r = \frac{p}{e \cos (\omega (\theta - g)) +1},
		\end{equation}
		here,  $p = \frac{2 (C^2 + \beta )}{\alpha}$, $\omega = \sqrt{\frac{{C^2 + \beta} }{C^2}}$, 
		and $e = \sqrt{1 + \frac{8E(C^2 + \beta )}{\alpha^2}} $.  %\footnote{\rd{This solution does not allow negative $\alpha$.}} 
		Figure \ref{fig:orb1}, \ref{fig:orb2} illustrate the orbits for $p=1, e= 0.8, \omega = 1.1$ and $p=1, e= 0.2, \omega = 10.1$, respectively.
		
		%Note that $\theta = \phi \pm \omega^{-1} \cdot 2 \pi \cdot  n$ with the number $n$ of perihelion that orbit has until the next reflection. Remember that $-\omega^{-1} \cdot \pi \leq \phi \leq \omega^{-1} \cdot \pi$. 

		{Note that for the repulsive case $\alpha < 0$, we necessarily have $E>0$, and in this case we get 
			\begin{equation}
				\label{eq: solution_repulsive}
				r = \frac{p}{e \cos (\omega (\theta - g)) -1},
			\end{equation}
			where $p = \frac{2 (C^2 + \beta )}{-\alpha}$, $\omega = \sqrt{\frac{{C^2 + \beta} }{C^2}}$ and $e = \sqrt{1 + \frac{8E(C^2 + \beta )}{\alpha^2}} $, from (\ref{eq: solution_trans}). 
			Note that in this case, we have the corresponding billiard system only in the region {$\{y \leq \gamma\} $}.}
		
		%\footnote{\rd{I am not sure I fully understand, but it seems to me that between two reflections at the wall, the orbit cannot have more than one pericenter/apocenter. }}
		%\footnote{\bl{Orbits like in Figure \ref{fig:orb2} can have multiple pericenters / apocenters between two reflections when the wall is sufficiently close to the origin. }}
		%\footnote{\rd{I think it is a good idea to uniformize the notations for the angles.} }
		%\footnote{\rd{In a strict sense, this solution of the central force problem is interesting, but incomplete when the motion is retrograde, due to the problem of angles that I remarked in previous version (though now I prefer to delete). It also exclude the case of unbounded orbits. I think it is a good idea to combine your solution with what I wrote that is now as marked in comments. }} 
		%\footnote{\bl{I think this parametrization $\theta ~(-\pi \leq \theta \leq \pi)$ covers full orbits even when they system proceeds retrogradely unlike the original parameter $\phi ~ (-\pi/ \omega \leq \theta \leq \pi / \omega)$ in Boltzmann's paper.}}

		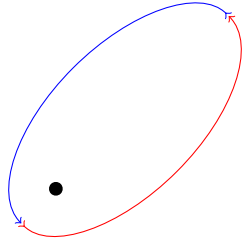
\begin{figure}
			\centering
			\begin{tikzpicture}
				\begin{scope}[scale=1.3, rotate=45]
					\draw[blue, >->]  (1.5,0)arc(0:180:1.5cm and 0.75cm);
					\draw[red,>->] (-1.5,0) arc(180:360:1.5cm and 0.75cm);
					\fill (-1,0) circle (2pt) coordinate (A);
				\end{scope}
			\end{tikzpicture}
			\caption{$r$-increasing direction (in red) and decreasing direction (in blue) for an ellipse}
		\end{figure}
		
		\begin{figure}
			\begin{tabular}{cc}
				\begin{minipage}[t]{0.49\hsize}
					\centering
					\includegraphics[keepaspectratio, scale=0.32]{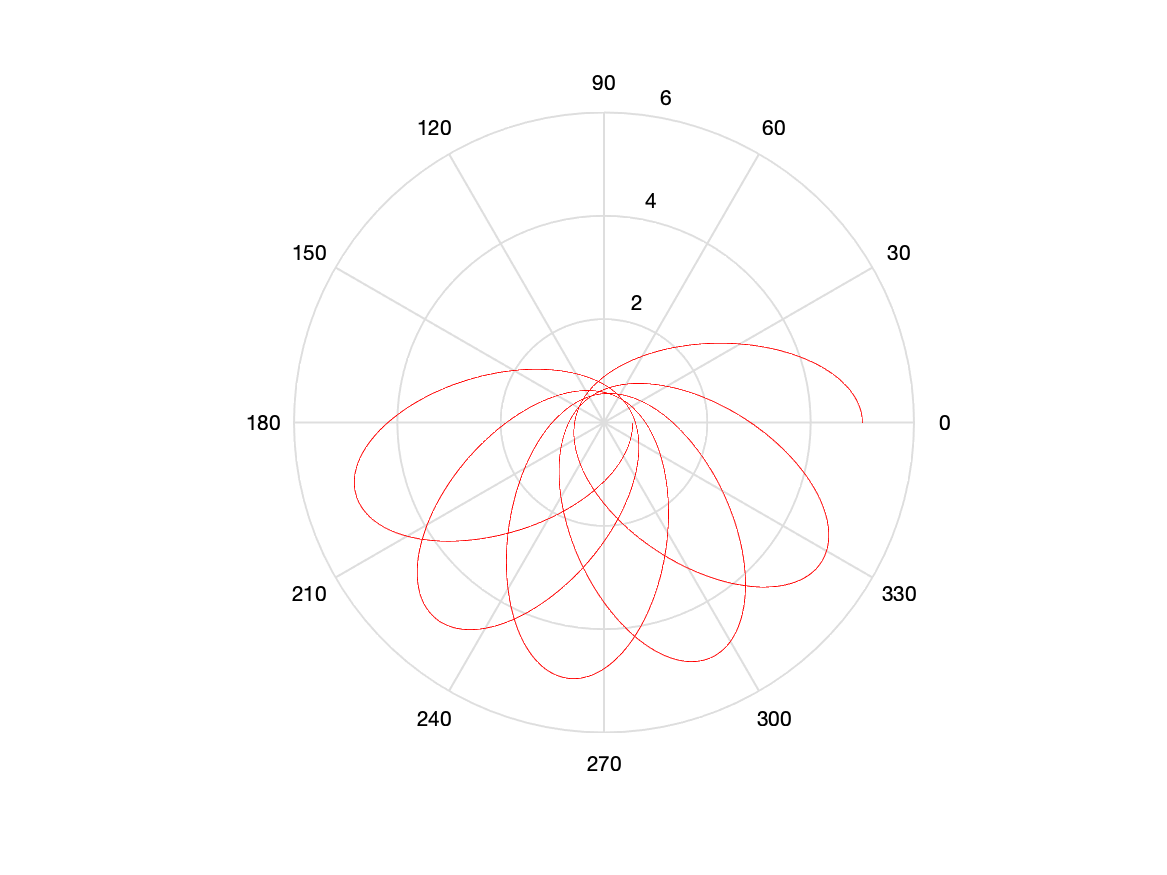}
					\caption{Orbit for $p=1, e= 0.8, \omega = 1.1$}
					\label{fig:orb1}
				\end{minipage} &
				\begin{minipage}[t]{0.49\hsize}
					\centering
					\includegraphics[keepaspectratio, scale=0.32]{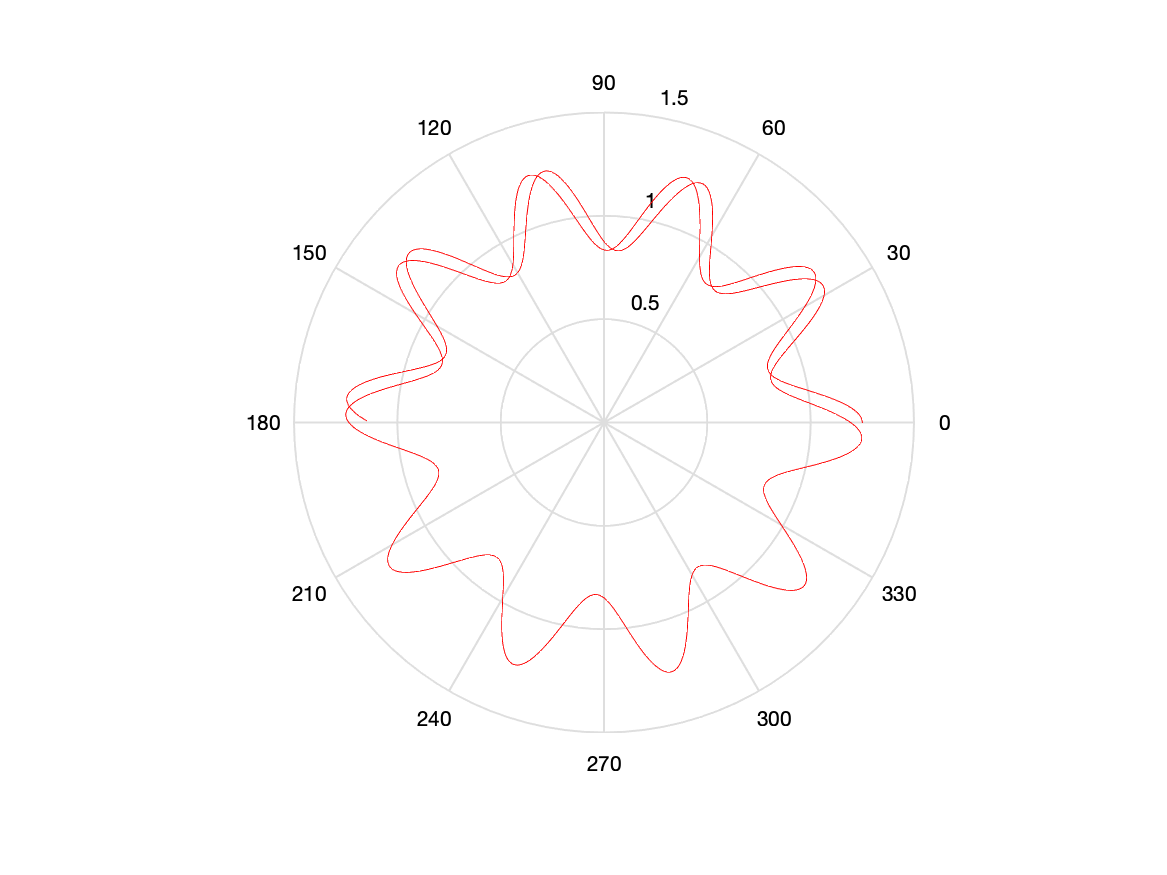}
					\caption{Orbit for $p=1, e= 0.2, \omega = 10.1$}
					\label{fig:orb2}
				\end{minipage}
			\end{tabular}
		\end{figure}
		
		%%%%%%%%%more general cases%%%%%%%%%%%%%%%%%%%%

		\begin{comment}
			From the equation (\ref{eq: first_der}), we have 
			\[
			\dot{r}^2 + \frac{C^2}{r^2} = - \frac{\alpha}{2 r^2} + \frac{\beta}{r^3}.
			\]
			By differentiating both sides with respect to $t$, we have 
			\[
			\ddot{r} - \frac{C^2 + \beta }{r^3} =  -\frac{\alpha}{2 r ^2}
			\]
			(remembering that $C$ is a constant of motion). Note that up to the change of constants this reduced system takes the same form as that of the Kepler problem.
			Taking the Clairaut variable $\rho = 1/r $ and having $C^{2}\dfrac{d^{2}\rho}{d\theta^{2}}=-r^{2}\ddot{r}$,
			the above equation is transformed into
			\begin{equation}
				\dfrac{d^{2}\rho}{d\theta^{2}}+\dfrac{C^{2}+\beta}{C^{2}}\rho=\dfrac{\alpha}{2C^{2}}.\label{eq:Clairaut}
			\end{equation}
			This is a second order linear differential equation whose general
			solution can be discussed case by case. We assume $\alpha\neq0$.
			We consider the case $C^{2}+\beta>0$. In this case the general solution
			is 
			\[
			\dfrac{1}{r}=\rho=\dfrac{\alpha}{2C^{2}}k\cos\omega(\theta-g)+\dfrac{\alpha}{2(C^{2}+\beta)},
			\]
			where $\omega=\sqrt{{(C^{2}+\beta)/C^{2}}}$ . 
			For $\alpha > 0$, this representation coincides with the previous form (\ref{eq: solution1}) of solutions. For repulsive cases, $\alpha < 0$, we write the solution equation as 
			\[
			r = \frac{p}{e \cos (\omega (\theta - g)) -1}. 
			\]
			
		\end{comment}
		
		Indeed, {since $r>0$}, we find solutions only in the case $e > 1$. For $e >1$, we have unbounded hyperbolic orbits look like those shown in Figure \ref{fig:orb_repel}. Notice that $r \to \infty $ when $e \cos (\omega (\theta - g)) \to 1$.
		
		\begin{figure}
			\centering
			\includegraphics[keepaspectratio, width=5cm]{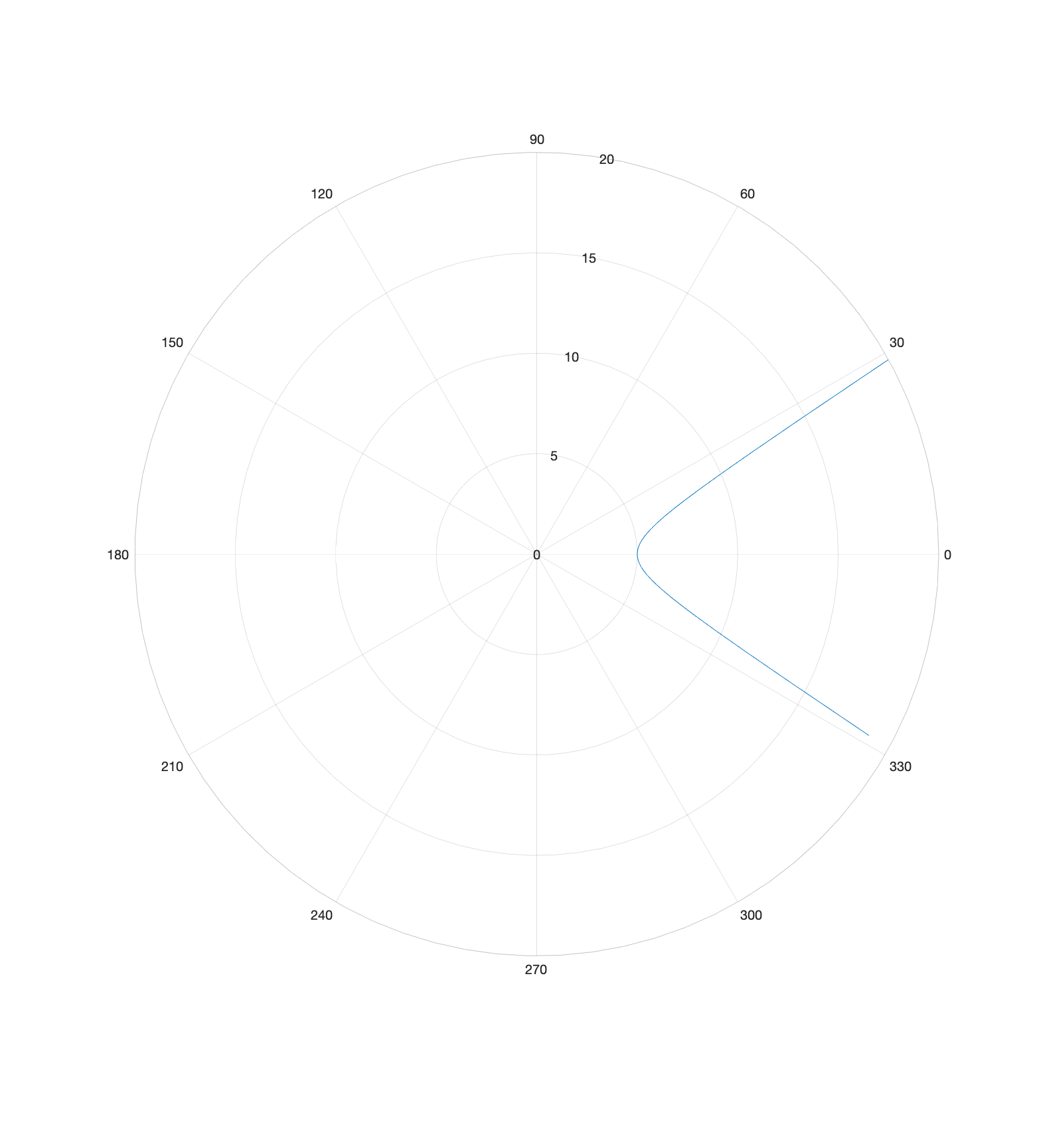}
			\caption{Orbit for $\alpha < 0 $}
			\label{fig:orb_repel}
		\end{figure}
		
		We now consider the other cases when {$C^2 + \beta \le 0$}.
		%From the equation (\ref{eq: first_der}), we have 
		%\[
		%\dot{r}^2 + \frac{C^2}{r^2} = - \frac{\alpha}{2 r^2} + \frac{\beta}{r^3}.
		%\]
		By differentiating the equation (\ref{eq: first_der}) with respect to $t$, we have 
		\[
		\ddot{r} - \frac{C^2 + \beta }{r^3} =  -\frac{\alpha}{2 r ^2}.
		\]
		
		Taking the Clairaut variable $\rho = 1/r $ and using $C^{2}\dfrac{d^{2}\rho}{d\theta^{2}}=-r^{2}\ddot{r}$,
		the above equation is transformed into
		\begin{equation}
			\dfrac{d^{2}\rho}{d\theta^{2}}+\dfrac{C^{2}+\beta}{C^{2}}\rho=\dfrac{\alpha}{2C^{2}}.\label{eq:Clairaut}
		\end{equation}

		When $C^{2}+\beta=0$, the equation reduces to
		
		\[
		\dfrac{d^{2}\rho}{d\theta^{2}}=\dfrac{\alpha}{2C^{2}}.
		\]
		
		So, the orbit takes the form 
		
		\[
		\dfrac{1}{r}=\rho=\dfrac{\alpha}{2C^{2}}\theta^{2}+k_{1}\theta+k_{2}
		\]
		
		which determines a spiral.
		
		When $C^{2}+\beta<0$, {the general solution is} in the
		form
		
		\[
		\dfrac{1}{r}=\rho=k\cos\omega(\theta-g)+\dfrac{\alpha}{2C^{2}},
		\]
		where $\omega = \sqrt{\frac{{C^2 + \beta} }{C^2}}$.
		Note that in this case $\omega$ is purely-imaginary and thus, the $\cos$
		appearing in the above formula {is actually} a $\cosh$. We may again
		put it into the form
		
		\[
		r=\dfrac{p}{1+e\cos\omega(\theta-g)},
		\]
		where $p = \frac{2 C^2}{\alpha}$ and $e = \frac{2k C^2}{\alpha}$.
		Note that $e$ may be either positive, negative, or zero. 
		
		We first discuss the case $\alpha<0$. In this case we have $p<0$ . When
		$e\ge0$, there are no solutions. When $-1\leq e<0$, the value of $\theta$ is restricted to $\theta < \theta_1 \leq \theta_2< \theta  $ with two limiting values $\theta_1, \theta_2$ such that $r \to \infty$ when $\theta \to \theta_1 -0$ or $\theta \to \theta_2 +0$, and $r \to 0$ when $\theta \to \pm \infty$, thus the orbit is an unbounded spiral.  When $e<-1$, we have $r < \dfrac{p}{1+e}$ and $r \to 0$ when $\theta \to \pm \infty$, thus the orbit is a bounded spiral  biasymptotic to the origin. 
		
		Secondly, we discuss the case $\alpha>0$. In this case $p>0$ . When $e>0$,
		we see that $r\to0$ when $\theta\to\pm\infty$, so the orbit is a
		spiral which is biasymptotic to the origin. When $e=0$, the orbit is a circle.
		When $-1<e<0$, $\theta$ is confronted between two limiting values,
		and the orbit is unbounded and has two asymptotic directions. There
		are no solutions when $e\le-1$.

		\subsection{Computation of General Boltzmann's Billiard Mapping}
		\label{sec:billiard}
		%\subsection{Billiard Mapping in Kepler case: $\beta =0$}
		%\subsection{Billiard Mapping in more general cases: $\beta \neq 0$}
		The billiard system is defined by adding a wall of reflection to the central force problem. 
		In this section, we assume $\alpha >0$ and $C^2 + \beta >0$.
		We define an arc as part of an orbit with
		starting and ending points on the reflection wall, and no other points
		hit the wall in between. The billiard mapping that sends a reflection point and a reflection velocity to the next extends to a mapping that maps an arc to another arc which then extends to a mapping of orbits. We shall analyze this mapping. 
		
		In polar coordinates, the wall of reflection $\left\{ y=\gamma>0\right\} $
		is represented by the equation $$r\sin\theta=\gamma.$$ {We
		compute the billiard mapping for which a given orbit}
		is reflected after reaching a point $(r_*, \theta_*)=( \gamma/ \sin\theta_{*},\theta_*)$ on the
		%\footnote{\rd{I find $_{1-2}$ and the like rather bad as notations:  the slash there can be thought of as a minus sign, which does not make sense for subscripts. Please change them to $_{1,2}$ or $_{1,2}$, or maybe just $_{*}$ }} at the
		wall.
		%, \emph{i.e.}, $(r_*,\theta_*)$ is the reflection point between the first orbit and the proceeding second orbit. 
		%\marginpar{\rd{I think the commented out sentence is not needed. Please check.}}
		
		Once fixing the energy, a {non-circular, non-singular} orbit is characterized
		by the coordinates $(g,C)$, {constructed in Section \ref{sec: cnst_coodinate}}. 
		Define the billiard mapping $S$ as $S(g_1, C_1) = (g_2,C_2)$, where $(g_1,C_1)$ and $(g_2,C_2)$ correspond to the orbits before and after the reflection.
		The derivatives $$\dfrac{dr}{d\theta}(r_*,\theta_*)=\dfrac{p\cdot e\cdot\omega\cdot\sin\omega(\theta_*-g_{1,2})}{(1+e\cos\omega(\theta_*-g_{1,2}))^{2}}$$
		will be denoted as $r'_{1,2}$ %in the two orbits 
		respectively. The
		derivatives $\dfrac{d\theta}{dt}(r_*,\theta_*)$ will be denoted
		as $\dot{\theta}_{1,2}$% in the two orbits
		 respectively. We also write  the corresponding $p,e,\omega$ in the two orbits as  $p_{1,2},$ $ e_{1,2},$ and $ \omega_{1,2}$.
		
		We get the following equations from the law of reflection.
		\[
		\begin{cases}
			r_*= \dfrac{\gamma}{\sin \theta_*} =\dfrac{p_{1}}{1+e_{1}\cos\omega_{1}(\theta_*-g_1)}=\dfrac{p_{2}}{1+e_{2}\cos\omega_{2}(\theta_*-g_2)}\\
			(r'_{1}\sin\theta_*+r_*\cos\theta_*)\dot{\theta}_{1}=-(r'_{2}\sin\theta_*+r_{*}\cos\theta_*)\dot{\theta}_{2}\\
			(r'_{1}\cos\theta_*-r_*\sin\theta_*)\dot{\theta}_{1}=(r'_{2}\cos\theta_*   {-}  r_*\sin\theta_*)\dot{\theta}_{2}\\
			\dot{\theta}_{1}=\dfrac{C_{1}}{r_*^{2}},\qquad\dot{\theta}_{2}=\dfrac{C_{2}}{r_*^{2}}
		\end{cases}
		\]
		
		From these {we deduce}
		
		\begin{align*}
		&C_{2}=\dfrac{C_{1}(-2r_*\cos^{2}\theta_*-2r'_{1}\sin\theta_*\cos\theta_*+r_*)}{r_*},\\
		&r'_{2}=\dfrac{-r_*r'_{1}\tan^{2}\theta_*-2r_*^{2}\tan\theta_*+r'_{1}r_*}{r_*\tan^{2}\theta_*-2r'_{1}\tan\theta_*-r_*},
		\end{align*}
		
		and  \textcompwordmark
		\[
		p_{2}=\dfrac{2(C_{2}^{2}+\beta)}{\alpha},\quad\omega_{2}=\sqrt{\frac{(C_{2}^{2}+\beta)}{C_{2}^{2}}}, \quad e_2 = \sqrt{1 + \frac{8E(C_2^2 + \beta )}{\alpha^2}}.
		\]
		Consequently, we obtain
		
		\[
		e_{2}\cos\omega_{2}(\theta_*-g_2)=\dfrac{p_{2}-r_*}{r_*},\qquad e_{2}\sin\omega_{2}(\theta_*-g_2)=\dfrac{p_{2}r'_{2}}{\omega_{2}r_*^{2}}
		\]
		
		From these {we} solve $g_2$ as
		\begin{equation}
			\label{eq: pericenter}
			g_2= \theta_* - \frac{ \operatorname{sign}\left(\dfrac{p_{2}r'_{2}}{e_2 \omega_{2}r_*^{2}}\right) \arccos \left(\dfrac{p_{2}-r_*}{e_2 r_*}\right)}{\omega_{2}}.
		\end{equation}
		{Remember that, when $\omega \neq 1$, there are multiple pericenters and apocenters. We choose the closest pericenter from the current reflection point, that is defined in (\ref{eq: pericenter}) as the next argument of pericenter.}
		%The
		%second orbit
		%\footnote{\rd{I deleted the word elliptic. They are not ellipses...}} 
		%is thus determined by $(C_2, g_2)$.
		In order to complete this inductive step, we compute the next reflection point $(r_{**}, \theta_{**})$ {with $0 < \theta_{**} < \pi$} from 
		\begin{equation}
			\label{eq:next_reflex}
			e_{2}\cos\omega_{2}(\theta_{**}-g_2)=\dfrac{p_{2}-r_{**}}{r_{**}} , \quad r_{**} \sin \theta_{**} = \gamma.
		\end{equation}
		 {{In general, Equation (\ref{eq:next_reflex}) have multiple solutions, with $(r_*, \theta_*)$ being one of them. We therefore} add the following condition to determine {the} next reflection point $(r_{**}, \theta_{**})$ %corresponding to the current one $(r_*, \theta_*)$
		 :
			\begin{equation}
				r \cdot \sin \theta  = \frac{p \sin \theta}{e \cos(\omega(\theta -g ))-1} \geq \gamma
			\end{equation}
			for all $\theta $ such that $\theta_* \leq \theta \leq \theta_{**} $ if $C > 0$ (for all $\theta$ such that $\theta_{**} \leq \theta \leq \theta_{*} $ if $C < 0$).}
		
		%After this, the billiard mapping is defined via the second elliptic
		%orbit, by sending $(r_*,\dot{\theta_{1}})$ to the other point of intersection
		%of the orbit with the velocity of incidence $(r_{**}, \dot{\theta_{2}})$. This is equivalent to the mapping %with arcs $(C_1, g_{1})$, by identifying
		%each arc with its end point. 
		%\footnote{It can be good to write these discussions at a more general level in a separate section on billiards.}

		\subsection{Solutions {of} the Central Force Problem: {The Case of Cotes' Spirals}}
		We {here} consider the case $\alpha=0, {\beta \neq 0}, C\neq0$, in this special case, solution curves of the central force problem with a force function $\beta/r^2$ are \emph{Cotes' spiral} \cite[Chapter IV]{Earnshaw}. 
		
		%From the equation (\ref{eq: first_der}), we have 
		%\[
		%\dot{r}^2 + \frac{C^2}{r^2} = - \frac{\alpha}{2 r^2} + \frac{\beta}{r^3}.
		%\]
		{Differentiating} the equation (\ref{eq: first_der}) with respect to $t$, we have 
		\[
		\ddot{r} - \frac{C^2 + \beta }{r^3} =  -\frac{\alpha}{2 r ^2}.
		\]
		%(remembering that $C$ is a constant of motion). 
	
		Taking the Clairaut variable $\rho = 1/r $ and having $C^{2}\dfrac{d^{2}\rho}{d\theta^{2}}=-r^{2}\ddot{r}$,
		the above equation is transformed into
		\begin{equation}
			\dfrac{d^{2}\rho}{d\theta^{2}}+\dfrac{C^{2}+\beta}{C^{2}}\rho=\dfrac{\alpha}{2C^{2}}.\label{eq:Clairaut}
		\end{equation}
		
		\begin{comment}
			
			This is a second order linear differential equation whose general
			solution can be discussed case by case. We assume $\alpha\neq0$.
			We consider the case $C^{2}+\beta>0$. In this case the general solution
			is 
			\[
			\dfrac{1}{r}=\rho=\dfrac{\alpha}{2C^{2}}k\cos\omega(\theta-g)+\dfrac{\alpha}{2(C^{2}+\beta)},
			\]
			where $\omega=\sqrt{{(C^{2}+\beta)/C^{2}}}$ . 
			For $\alpha > 0$, this representation coincides with the previous form (\ref{eq: solution1}) of solutions. For repulsive cases, $\alpha < 0$, we write the solution equation as 
			\[
			r = \frac{p}{e \cos (\omega (\theta - g)) -1}. 
			\]
		\end{comment}
		
		By substituting $\alpha=0$, the equation (\ref{eq:Clairaut}) can be written into

		\[
		\dfrac{d^{2}\rho}{d\theta^{2}}+\left(\dfrac{C^{2}+\beta}{C^{2}}\right)\rho=0.
		\]
		
		We discuss different subcases. 
		
		When $C^{2}+\beta>0$,  the general solution of the equation is written
		as
		
		\[
		\dfrac{1}{r}=\rho=k\cos{\omega(\theta-\psi)},
		\]
		where $k \in \R$ and $\psi \in [0, 2 \pi)$.
		When $C^{2}+\beta=0$, the general solution reduces to the form
		
		\[
		\dfrac{1}{r}=\rho=k_{1}\theta+k_{2}.
		\]
		
		When $C^{2}+\beta<0$, the general solution is
		
		\[
		\dfrac{1}{r}=\rho=k_{1}\exp\left({i\omega(\theta-\psi}\right)+k_{2}\exp\left({-i\omega(\theta-\psi)}\right),
		\]
		with a purely imaginary $\omega$.
		
		To make further analysis observe that 
		\[
		h:= \left( \dfrac{d\rho}{d\theta} \right) ^{2}+\left(\dfrac{C^{2}+\beta}{C^{2}}\right)\rho^{2}
		\]
		is a first integral of the equation. Drawing its level sets in the
		phase space with coordinates $\bigl(\dfrac{d\rho}{d\theta},\rho\bigr)$, we see
		that the level sets are hyperbolae and {bifurcate} at {the} zero-level $\left\{ h=0\right\} $
		through a degeneration into a pair of lines, and {then continue as hyperbolae with the major
		axis switched.}
		
		When $h<0$, the hyperbola has {the $\rho-$axis as major axis and}
		admits a parametrization with hyperbolic functions. The corresponding
		solution in polar form is
		
		\[
		\dfrac{1}{r}=\rho=k\cos\omega(\theta-\psi).
		\]
		
		Similarly, when $h>0$, we get 
		
		\[
		\dfrac{1}{r}=\rho=k\cdot i\cdot \sin\omega(\theta-\psi).
		\]
		
		{And}, when $h=0$ we have\textcompwordmark{}
		
		\[
		\dfrac{1}{r}=\rho=k\exp\left(\pm i\omega(\theta-\psi)\right).
		\]
		
		We thus get the five classes of Cotes' spirals as orbits of the problem
		with $\alpha=0$.
		
		\subsection{Computation of the Billiard Mapping: Cotes' Spiral Case}
		\label{sec:cotes}
		{We here} compute the billiard mapping for the special case $\alpha = 0$. We again {only} consider bounded orbits, thus we assume $E < 0$.
		
		The doubled total energy is written as
		\[
		{2 E= } \,\dot{r}^2 + \frac{C^2 + \beta}{r^2}, 
		\]
		{which} leads to 
		\[
		h =  \left( \dfrac{d\rho}{d\theta} \right) ^{2}+(\dfrac{C^{2}+\beta}{C^{2}})\rho^{2} = \frac{2E}{C^2}{<0}.
		\]
		%From these equations and $E < 0 $, we always have 
		{Thus $C^2 + \beta <0 $} and 
		%$$h:= \left( \dfrac{d\rho}{d\theta} \right) ^{2}+(\dfrac{C^{2}+\beta}{C^{2}})\rho^{2} <0.$$
		 %Therefore, 
		 
		the  orbits are given in the form
		\[
		\dfrac{1}{r}=\rho=k\cos\omega(\theta-\psi).
		\] 
		{with} $\omega = \sqrt{\frac{C^2 + \beta}{C^2}}, k = \sqrt{\frac{2E}{\omega^2 C^2}}$ and $\psi$ the {argument of the apocenter}. 
		
		{We} consider the billiard mapping $(\psi_1,C_1) \mapsto (\psi_2,C_2)$. Let $(r_*, \theta_*)$ be the {reflection} point.  
		%between the first orbit and the second one. 
		The derivatives $$\dfrac{dr}{d\theta}(r_*,\theta_*)=\dfrac{\omega \sin \omega (\theta_* -\psi_{1,2} ) }{k \cos^2  \omega (\theta_* -\psi_{1,2})}$$
		{are} denoted as $r'_{1,2}$ respectively. The
		derivatives $\dfrac{d\theta}{dt}(r_*,\theta_*)$ {are} denoted
		as $\dot{\theta}_{1,2}$ respectively. We also write the corresponding $p,e,\omega$ in the two orbits as  $p_{1,2},$ $ e_{1,2},$ and $ \omega_{1,2}$. {The next reflection point is computed using} the following equations 
		\[
		\begin{cases}
			r_*= \dfrac{\gamma}{\sin \theta_*} =\dfrac{1}{k_1 \cos \omega_1 (\theta_* - \psi_1)}=\dfrac{1}{k_2 \cos \omega_2 (\theta_* - \psi_2)}\\
			(r'_{1}\sin\theta_*+r_*\cos\theta_*)\dot{\theta}_{1}=-(r'_{2}\sin\theta_*+r_{*}\cos\theta_*)\dot{\theta}_{2}\\
			(r'_{1}\cos\theta_*-r_*\sin\theta_*)\dot{\theta}_{1}=(r'_{2}\cos\theta_*   {-}  r_*\sin\theta_*)\dot{\theta}_{2}\\
			\dot{\theta}_{1}=\dfrac{C_{1}}{r_*^{2}},\qquad\dot{\theta}_{2}=\dfrac{C_{2}}{r_*^{2}}.
		\end{cases}
		\]
		From these one deduces that
		
		\begin{align*}
		&C_{2}=\dfrac{C_{1}(-2r_*\cos^{2}\theta_*-2r'_{1}\sin\theta_*\cos\theta_*+r_*)}{r_*},\\
		&r'_{2}=\dfrac{-r_*r'_{1}\tan^{2}\theta_*-2r_*^{2}\tan\theta_*+r'_{1}r_*}{r_*\tan^{2}\theta_*-2r'_{1}\tan\theta_*-r_*},
		\end{align*}
		and then
		$$\omega_2 = \sqrt{\frac{C_2^2 + \beta}{C_2^2}},$$  $$k_2 = \sqrt{\frac{2E}{\omega_2^2 C_2^2}}.$$
		
		Consequently, we obtain 
		\[
		\frac{1}{r_*} = k_2 \cos \omega_2 (\theta_* - \psi_2), \quad - \frac{1}{r_*^2}  r'_2 = - k_2 \omega_2 \sin \omega (\theta_* - \psi_2).
		\]
		
		{Thus $\psi_2$ can be solved as}
		\[
		\psi_2 = \theta_* - \frac{\sign \left(-\frac{r'_2}{r_* k_2 (- i \omega)}\right) \arccosh \left(\frac{1}{k_2 r_*}\right)}{-i \omega_2}.
		\]
		
		%The
		%second orbit
		%\footnote{\rd{I deleted the word elliptic. They are not ellipses...}} 
		%is thus determined by $(C_2, \psi_2)$. In order to complete this inductive step, we finally compute
		 {The next reflection point $(r_{**}, \theta_{**})$, {$0 < \theta_{**} < \pi$} is then computed from }
		\begin{equation}
			\label{eq: next_reflex_cotes}
			\frac{1}{r_{**}} = k_2 \cos \omega_2 (\theta_{**} - \psi_2) , \quad r_{**} \sin \theta_{**} = \gamma.
		\end{equation}
		%Remember that $0 \leq \theta_{**} \leq \pi$.

		\section{{Numerical Results of Boltzmann's Billiard Trajectories}}
		\label{sec: numerical_trajectory}
		We here present some numerical {simulations} of Boltzmann's billiard mapping {based on  Section \ref{sec:billiard}}. In our simulations, we set $\alpha = 4, E= -0.5, \gamma= 0.5$, and vary the parameter $\beta \geq0$.
		Figures in this section
		% represents two pictures for each $\beta$, the one shows the graphs of  $e \cos (\theta -g) - \dfrac{p \sin \theta}{\gamma} $ from the orbit equation (\ref{eq:orbeq}) for $-\pi \leq \theta \leq \pi$, and the other one 
		illustrate numerically computed trajectories of the billiard mapping \emph{i.e.} the evolving values of 
		$$(g,C) \in [0, 2 \pi ) \times [C_{min}, C_{max}]$$ 
		at each reflection. {All numerical computations here have been operated by MATLAB and the interval
		arithmetic \cite{Rump} has been used to count the solutions satisfying Equation \eqref{eq:next_reflex}  or \eqref{eq: next_reflex_cotes}.}
		
		For $\beta = 0$, our simulation shows periodic behavior of a trajectory as we illustrated in Figure \ref{fig:beta000}, which is compatible with the integrability of the system %for this parameter setting.
		 For small $\beta$, for example $\beta = 0.5$ the system remains quasi-periodicity and seems not to be transitive, see figure \ref{fig:beta050}. For a bigger value of $\beta$, the {(quasi-)}periodicity may {break}, and chaotic behavior {appears}, as we illustrated for the case $\beta = 2.6$ in Figure \ref{fig:beta260}. In this case, {it seems possible} that a single orbit densely covers the whole energy hypersurface.
		 % in state space (for fixed $H = E$). 
		 Therefore {for big enough $\beta$, it is possible to have ergodic systems.} However, chaotic behavior does not always show up for large $\beta$. {Figure \ref{fig:beta240} shows both quasi-periodic (Subfig. a) and chaotic (but not transitive) behavior (Subfig. b) for $\beta = 2.4$, with different initial values. Subfig. c indicates the existence of 2-period orbit for this parameter setting. }
		
		\begin{figure}
			%	\begin{tabular}{cc}
				%		\begin{minipage}[t]{0.49\hsize}
					\centering
					%			\includegraphics[keepaspectratio, width=6cm]{orb_beta000.png}
					%			\hspace{1.6cm}
					%			a. The orbits in the billiard system
					%		\end{minipage} &
				%		\begin{minipage}[t]{0.49\hsize}
					\centering
					\includegraphics[keepaspectratio, width=6cm]{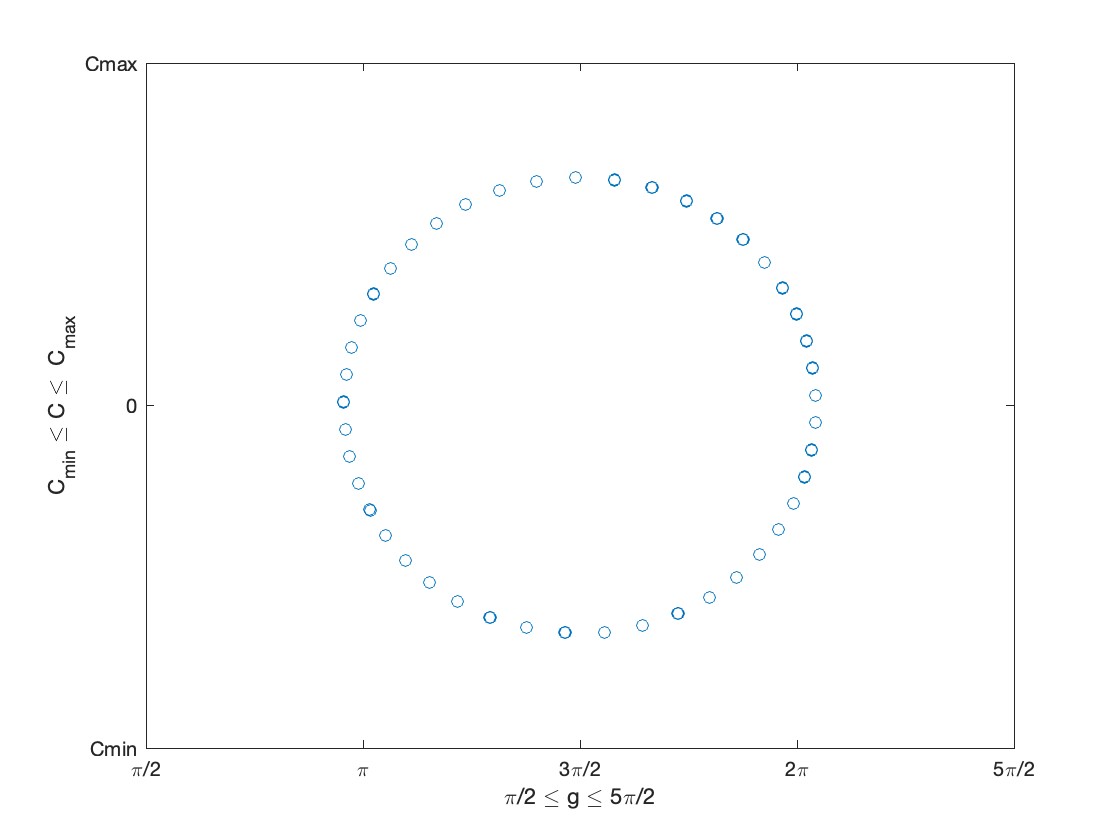}
					%		\end{minipage}
				%	\end{tabular}
			\caption{Periodic behavior of the mapping trajectory for $\beta = 0$}
			\label{fig:beta000}
		\end{figure}

		\begin{figure}
			%\begin{tabular}{cc}
			\centering
			\includegraphics[keepaspectratio, width=6cm]{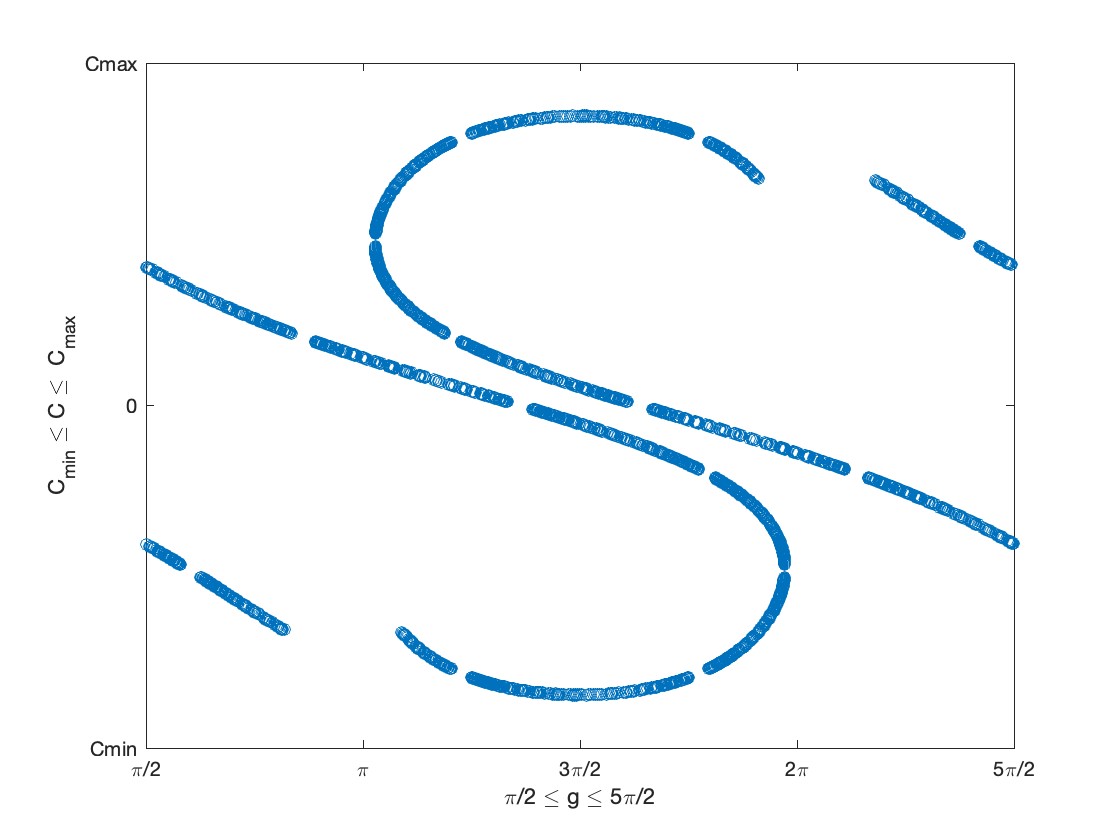}
			%\end{tabular}
			\caption{(Quasi-)periodic behavior of the mapping trajectory for $\beta= 0.5$}
			\label{fig:beta050}
		\end{figure}

		\begin{figure}
				\begin{tabular}{cc}
				\begin{minipage}[t]{0.49\hsize}
			\centering
			\includegraphics[keepaspectratio, width=6cm]{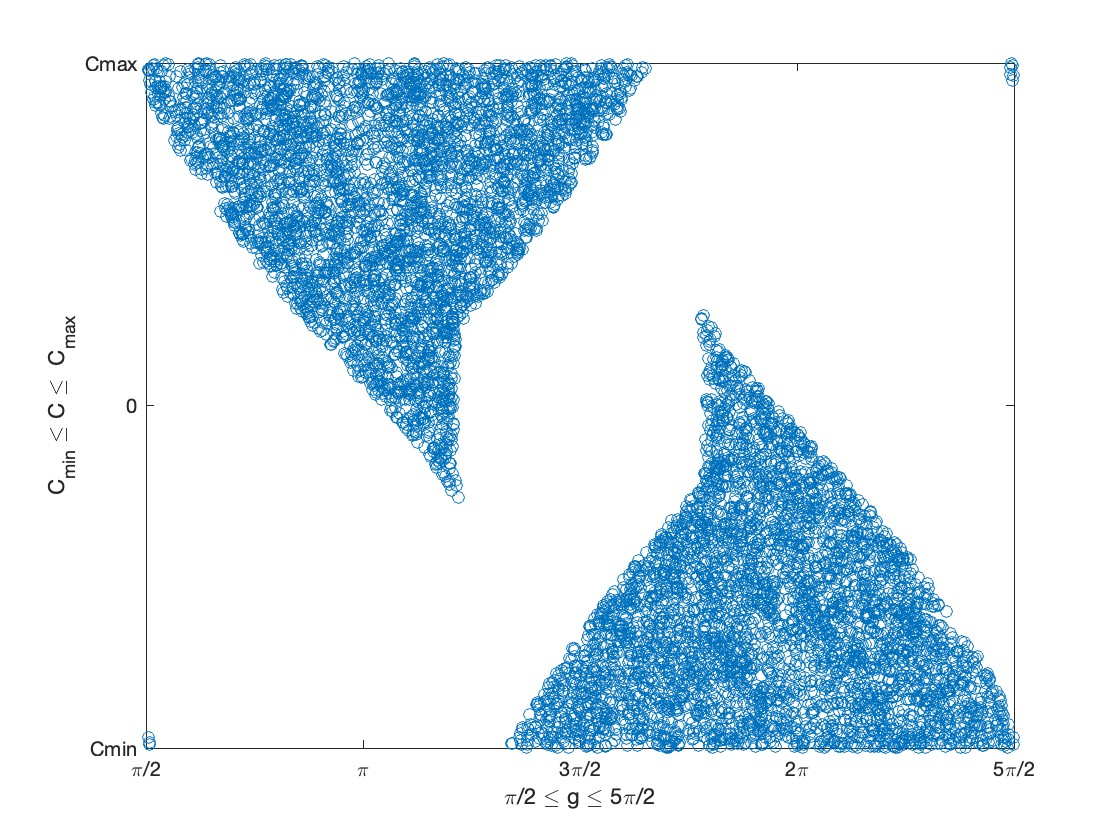}
			a. Trajectory with initial value $(g_0, C_0)=(0.2, 0.8)$
			%\end{tabular}
		\end{minipage} &
	\begin{minipage}[t]{0.49\hsize}
	\centering
	\includegraphics[keepaspectratio, width=6cm]{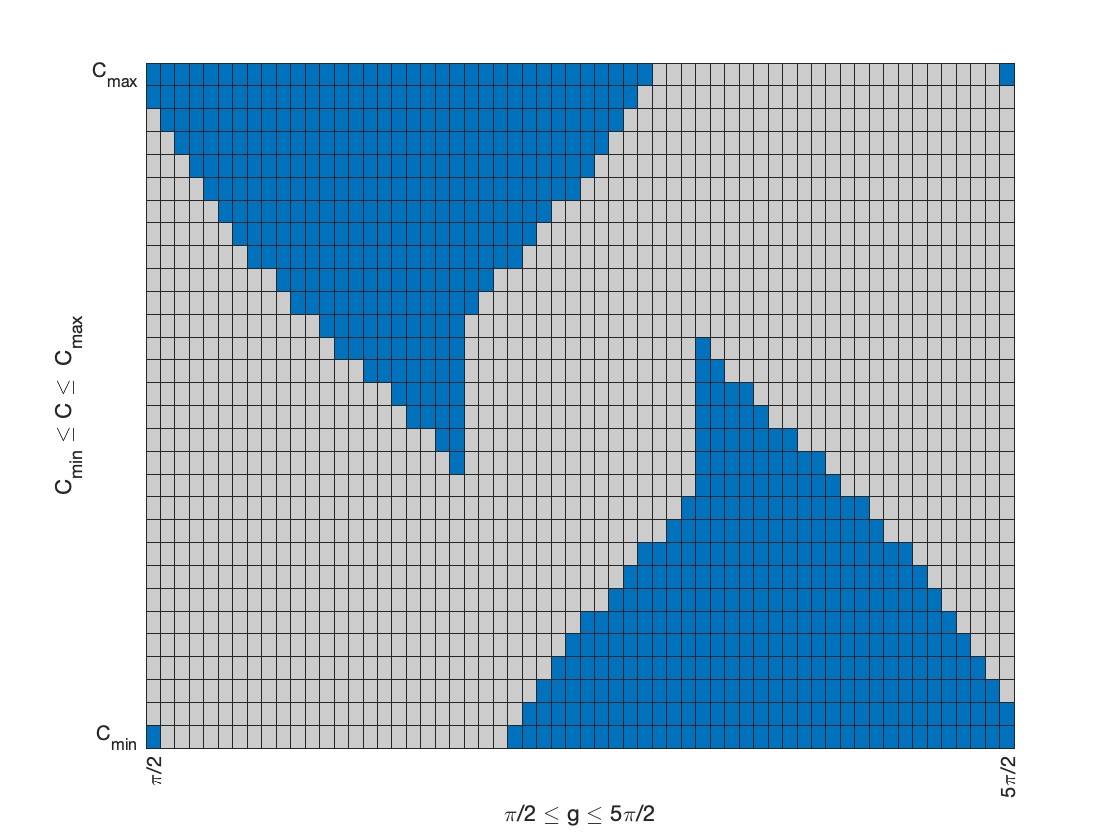}
	b. Discretized allowed region (in blue)
\end{minipage}
\end{tabular}
			\caption{Transitive behavior of the mapping trajectory for $\beta= 2.6$}
			\label{fig:beta260}
		\end{figure}

		\begin{figure}
			\begin{tabular}{cc}
				\begin{minipage}[t]{0.49\hsize}
					\centering
					\includegraphics[keepaspectratio, width=6cm]{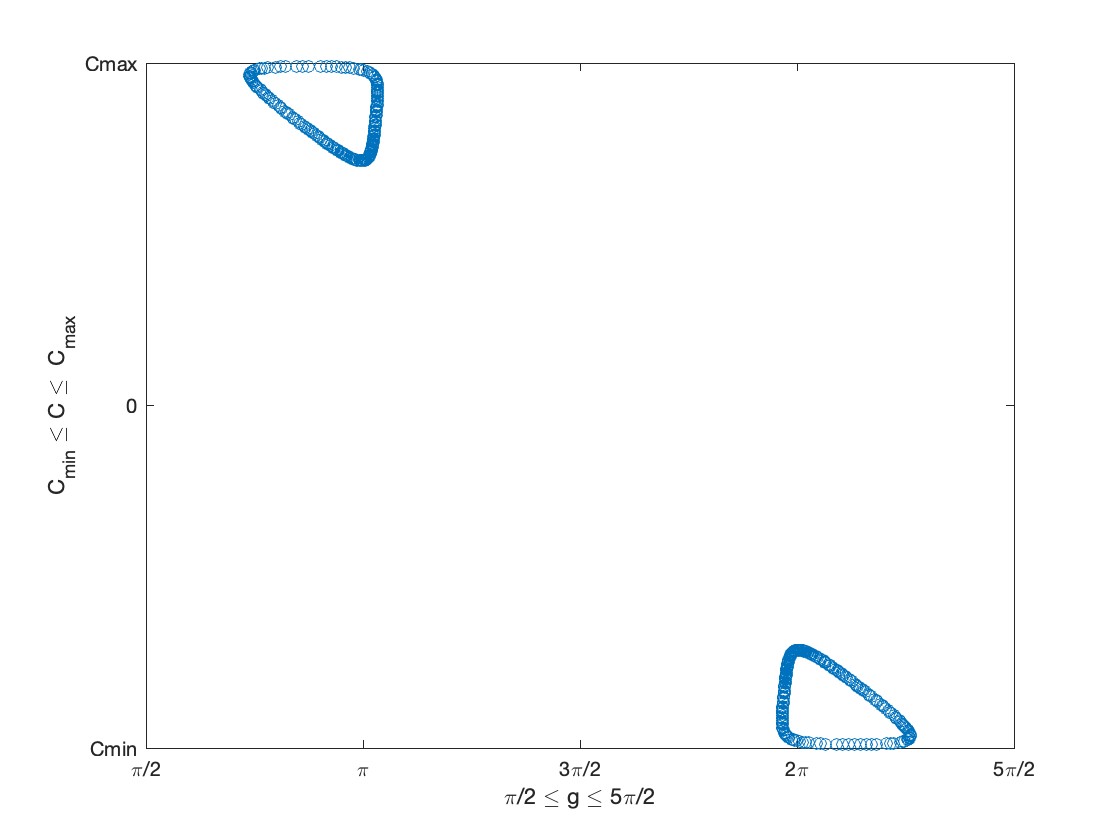}
					\hspace{1.6cm}
					a. Trajectory with initial value $(g_0, C_0)=(0.1, 1.1)$
					\includegraphics[keepaspectratio, width=6cm]{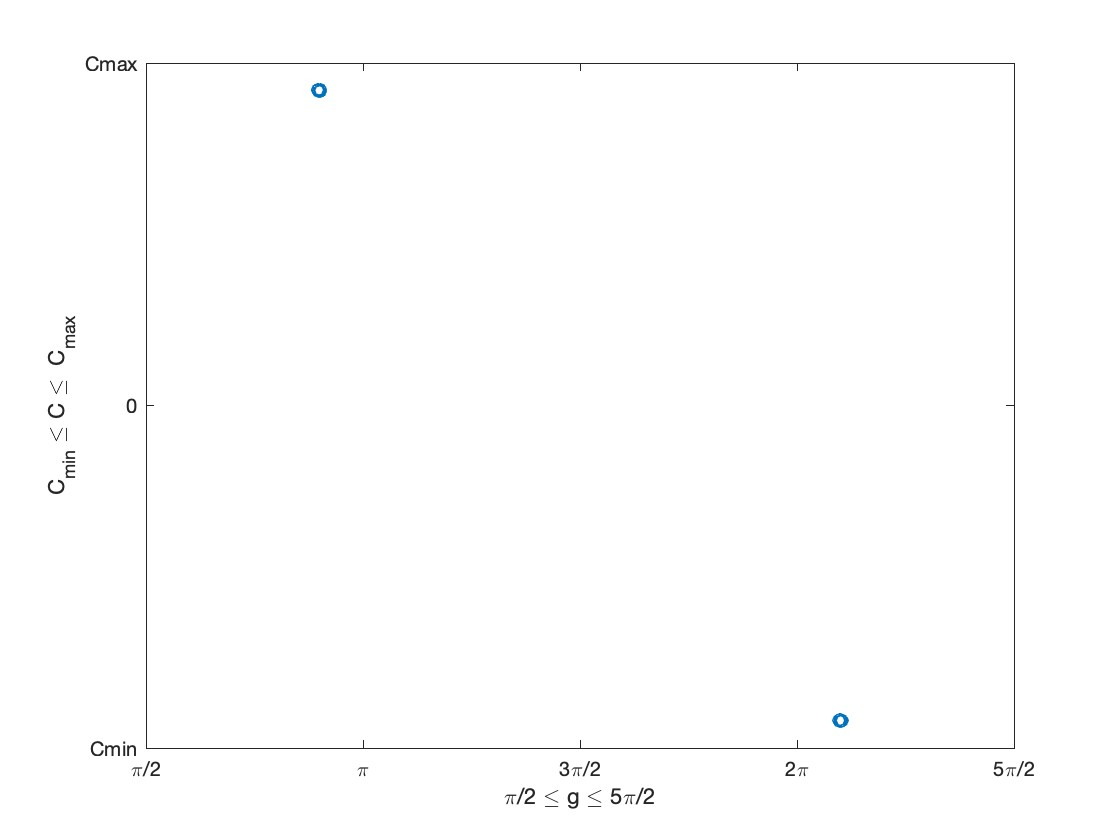}
					\hspace{1.6cm}
					c. Periodic trajectory with initial value $(g_0, C_0)=(3.45, -1.16)$
				\end{minipage} &
				\begin{minipage}[t]{0.49\hsize}
					\centering
					\includegraphics[keepaspectratio, width=6cm]{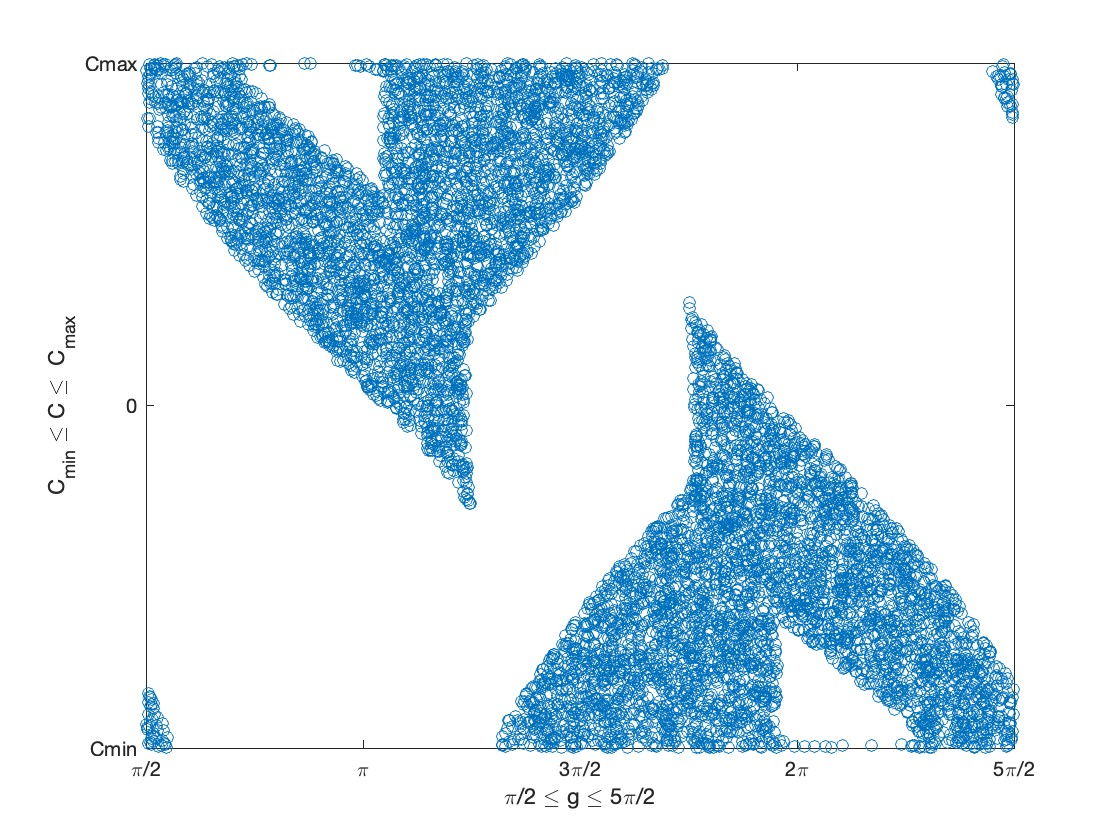}
					\hspace{1.6cm}
					b. Trajectory with initial value $(g_0, C_0)=(0.2, 1.0)$
				\end{minipage}
			\end{tabular}
			\caption{(Quasi-)periodic and chaotic behavior for $\beta = 2.4$}
			\label{fig:beta240}
		\end{figure}

		\begin{comment}
			
			We also simulate the numerical trajectories for the billiard mapping of Cotes' spiral, that was discussed in section \ref{sec:cotes}. For this case we set our parameter as $\alpha = 0, \beta = 2.0$. From Figure \ref{fig:spiral} we can see that single orbits may not cover the whole space and the system may have non-trivial invariant subset. Therefore, for Cotes' spiral case ($\alpha = 0$), the billiard system is likely non-ergodic.

			\begin{figure}
				%\begin{tabular}{cc}
				\centering
				\includegraphics[keepaspectratio, width=6cm]{spiral_beta200_800.png}
				%\end{tabular}
				\caption{The values of $(g, C)$ in the billiard system of Cotes' spiral}
				\label{fig:spiral}
			\end{figure}
		\end{comment}
		
		\section{{The Koopman Operator and Its Eigenvalue Problem}}
		\label{sec: Koopman}
		
		For any measure-preserving map $S$ on a probability measure space $(X, \mu, \Sigma)$, the Koopman operator can be defined as the transfer operator on $L^2(X):=L^2(X, {\mu;} \, \C)$ by 
		\begin{equation}
			K f := f \circ S, \qquad f \in L^2(X).
		\end{equation}
		%From the preservation of measure (i.e. $\mu(A) = \mu(S^{-1}(A)), \forall A \in \Sigma$), one \rd{sees} 
		{Since $S$ is measure preserving,} the Koopman operator $K: L^2(X) \to L^2(X)$
		% preserves the \rd{$L^{2}-$}norm. Therefore, this operator 
		is unitary and has its spectrum on the unit circle. 
		The spectrum of the Koopman operator carries essential {dynamical} information {of} the {map} $S$. %(\emph{e.g}., ergodicity, weakly mixing, invariant sets). 
		%The following proposition connects the ergodic property of the original mapping $S$ and the eigenvalue problem of the corresponding Koopman operator.
		{In particular, we have}
		
		\begin{prop}
			Let $S$ be a measure-preserving map on a probability measure space $(X, \mu, \Sigma)$ and let $K: L^2(X) \to L^2(X)$ be the corresponding Koopman operator.  Then $1$ is an eigenvalue of $K$. Moreover, the map $S$ is ergodic if and only if eigenvalue 1 is simple.
		\end{prop}
		See \cite[Proposition 7.15]{Eisner} for the proof.

		%Our aim here is to prove ergodicity of the billiard mapping by showing the simplicity of eigenvalue 1 of the corresponding Koopman operator. 
		%Before working on rigorous computation of eigenvalue problem, we prefer to solve approximated eigenvalue problem by discretizing the operator to get rough pictures of eigenvalues and eigen functions. 
		In the following, we {numerically investigate} the eigenvalue problem of {the} Koopman operator with {the} Galerkin method. 
		{All numerical computations here have been operated by MATLAB and the interval
			arithmetic \cite{Rump} has been used to count the solutions satisfying Equation \eqref{eq:next_reflex}.}
		
		\subsection{Approximation of Koopman Eigenvalue Problem with Galerkin Method }
		\label{sec: Koopman_approx}
		We here explain the approximation procedure of the Koopman eigenvalue problem using Galerkin method \cite{Ern} with piecewise constant basis functions. 
		\paragraph{Galerkin Method and Midpoint Quadrature with Uniform Weights}
		Consider the original eigenvalue problem of the Koopman operator on $L^2(X)$
		\begin{equation*}
			Ku = \lambda u, \qquad u \in L^2(X),
		\end{equation*}
		{which can be transformed }into the equivalent equation
		\begin{equation*}
			\langle K u , v \rangle  = \lambda \langle u , v \rangle , \qquad \forall v \in L^2(X),
		\end{equation*}
		where $\langle \cdot , \cdot \rangle $ is the inner product of the Hilbert space $L^ 2(X)$. 
		{We fix} finitely many basis functions $\{ f_1, \cdots, f_N  \}$ in $L^2(X)$. We look for approximate eigenfunctions in the form $u = \sum_{n = 1}^{N }  \alpha_i f_i$ and we restrict the above equation to the space which is spanned by the base functions. Then {we have} %it can be rewritten as 
		\begin{equation*}
			\sum _{n = 1}^{N} \alpha_n \langle f_n \circ S, f_m \rangle = \lambda \sum_{n = 1}^{N} \alpha _n \langle f_n, f_m \rangle, \qquad \forall m \in \{1, \cdots, N\}.
		\end{equation*}
		{In matrix form this is}
		\begin{equation}
			\label{eq: mat_eigen}
			\begin{pmatrix}
				\langle f_1 \circ S, f_1 \rangle & \cdots &  \langle f_N \circ S, f_1 \rangle \\
				\\
				\vdots & \ddots &\vdots  \\
				\\
				\langle f_1 \circ S, f_N \rangle & \cdots &  \langle f_N \circ S, f_N \rangle 
			\end{pmatrix}
			\begin{pmatrix}
				\alpha _1 \\
				\\
				\vdots \\
				\\
				\alpha _N
			\end{pmatrix}
			=
			\lambda
			\begin{pmatrix}
				\langle f_1 , f_1 \rangle & \cdots &  \langle f_N , f_1 \rangle \\
				\\
				\vdots & \ddots &\vdots  \\
				\\
				\langle f_1 , f_N \rangle & \cdots &  \langle f_N , f_N \rangle 
			\end{pmatrix}
			\begin{pmatrix}
				\alpha _1 \\
				\\
				\vdots \\
				\\
				\alpha _N
			\end{pmatrix}
		\end{equation}.
		
		For the computation of each entry of the matrices above, we divide the domain $ X$  of the mapping $S$ into finitely many disjoint {regions} $\Omega_1, \cdots \Omega _N $ so that $X = \sqcup_{n = 1}^{N} \Omega_n$. Suppose that our basis functions $f_n$ are the characteristic functions of each region $\Omega_n$ \emph{i.e.} $f_n(x) = 1 $ if $x \in \Omega_n$ and $f_n(x) = 0$ otherwise. Then the matrix in the left hand side of (\ref{eq: mat_eigen}) can be written as
		\begin{align}
			\begin{split}
				\label{eq: int_approx1}
				\langle f_n \circ S, f_m \rangle &= \int_{X} f_n(S(x)) \cdot f_m(x) dx\\
				& = \int_{\Omega_m} f_n(S(x)) dx \\
				&\approx \sum_{\ell = 1}^{L} w_{\ell}^{(m)} f_n(S(x_\ell)) 
			\end{split}
		\end{align}
		In the last line, we approximated the integral with the weighted summation of $f_n(S(x_\ell))$ over $L$ nodes in $\Omega_m$ which is chosen by the midpoint rule.

		%Note that $f_n(S(x_k)) = 1 $ if $S(x_k) \in \Omega_n$ and $f_n(S(x_k)) = 0$ otherwise. 
		If we set the same weight $w^{(m)} =  w_{\ell}^{(m)}$ at all nodes $\{x_\ell\}_{\ell =1}^L$ in $\Omega_m$, then we can simplify the above {formula} as
		\begin{align}
			\begin{split}
				\label{eq: int_approx2}
				\langle f_n \circ S, f_m \rangle 
				&\approx \sum_{\ell = 1}^{L} w_{\ell}^{(m)} f_n(S(x_\ell)) \\
				& =    w^{(m)} \cdot   \# \{ \ell \mid S(x_\ell) \in \Omega_n  \}\\
				& =  | \Omega_m |  \cdot  \frac{ \# \{ \ell \mid S(x_\ell) \in \Omega_n  \}}{ L},
			\end{split}
		\end{align}
		where $ | \Omega_m |$ is the measure of $\Omega_m$. In the last equation, we used 
		\[
		| \Omega_m | =  \int_{\Omega_m} dx = \sum_{\ell  ~\text{s.t.} x_\ell \in \Omega_m}  w_\ell^{(m)}=  \# \{  \ell \mid x_\ell \in \Omega_m \}  \cdot w^{(m)}  = L  \cdot w^{(m)}.
		\]
		
		The matrix in the left hand side of (\ref{eq: mat_eigen}) becomes
		\begin{align*}
			\langle f_n , f_m \rangle  & = \int_{\Omega} f_n(x) \cdot f_m(x) dx \\
			& = \int_{\Omega_m} f_n (x) dx\\
			& \approx \sum_{\ell =1}^{L} w_{\ell}^{(m)} f _n(x_\ell)\\
			& = w^{(m)} \cdot \#\{ \ell \mid x_\ell \in \Omega_n  \}\\
			& = 
			\begin{cases}
				| \Omega_n |  \qquad \text{if $n = m$} \\
				0 \qquad \text{otherwise.}
			\end{cases} 
		\end{align*}
		
		We call the matrix eigenvalue problem \eqref{eq: mat_eigen} approximated in the above way the {\emph{discretized Koopman eigenvalue problem.}}
		
		We now set $X= [0, 2 \pi ) \times [C_{min}, C_{max}]$, {$x= (g,C) \in X$} and let $S$ be {Boltzmann's} billiard mapping computed in Section \ref{sec: Boltzmann_billiard_mapping} and consider the approximated eigenvalue problem of the corresponding Koopman operator.
		
		In the following numerical computations, we divided the $(g,C)-$coordinate space $[0, 2 \pi ) \times [C_{min}, C_{max}]$ into $N=800$ partial sets. The number $L=25$ represents the number of the test nodes in each section used to approximate integrals, which appear in the equations \eqref{eq: int_approx1} and \eqref{eq: int_approx2}. We note that the billiard mapping $S$ is not defined on the whole space $[0, 2 \pi ) \times [C_{min}, C_{max}]$, therefore we need to restrict the divided space into the subset of all partitions where the corresponding orbits of the underlying mechanical system have at least two intersection points with the reflection wall $ y = \gamma$ %so that the reflection occurs at the wall. 
		In our computations, we set $\alpha = 4.0, E= -0.5, \gamma= 0.5$, and vary the parameter $\beta$.

		In Figure \ref{fig:eigen_beta000_gl}, Subfig. a shows the restricted region in a divided space $[0, 2 \pi ) \times [C_{min}, C_{max}]$ where the billiard mapping is well-defined for $\alpha = 4.0, \beta = 0.0$, and Subfig, b shows the all eigenvalues of the discretized Koopman eigenvalue problem, Subfig. c,d,e, and f show the level sets of all independent eigenfunctions corresponding to the three closest eigenvalues from 1. 
		Figure. \ref{fig:eigen_beta050_gl}, Figure. \ref{fig:eigen_beta240_gl} and Figure. \ref{fig:eigen_beta260_gl} show the same information on the discretized Koopman eigenvalue problem as Figure. \ref{fig:eigen_beta000_gl} but for the different parameter setting $\beta = 0.5, \beta = 2.4$, and $\beta = 2.6$, respectively.

		\begin{figure}
			\centering
			\begin{tabular}{cc}
				\begin{minipage}[t]{0.5\hsize}
					\includegraphics[keepaspectratio, height=5cm]{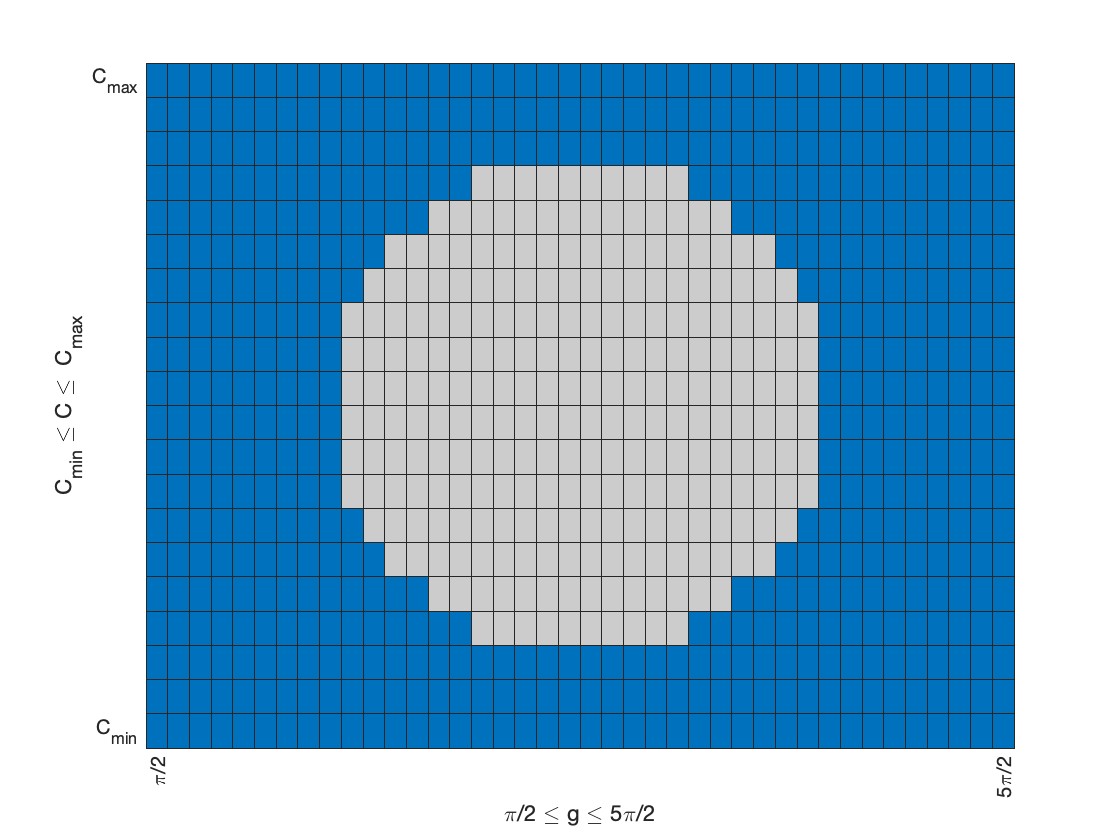}\\
					%\vspace{-4.5mm}
					\centering
					\text{  a. Allowed regions (in blue) }
				\end{minipage}
				\begin{minipage}[t]{0.5\hsize}
					\includegraphics[keepaspectratio, height=5cm]{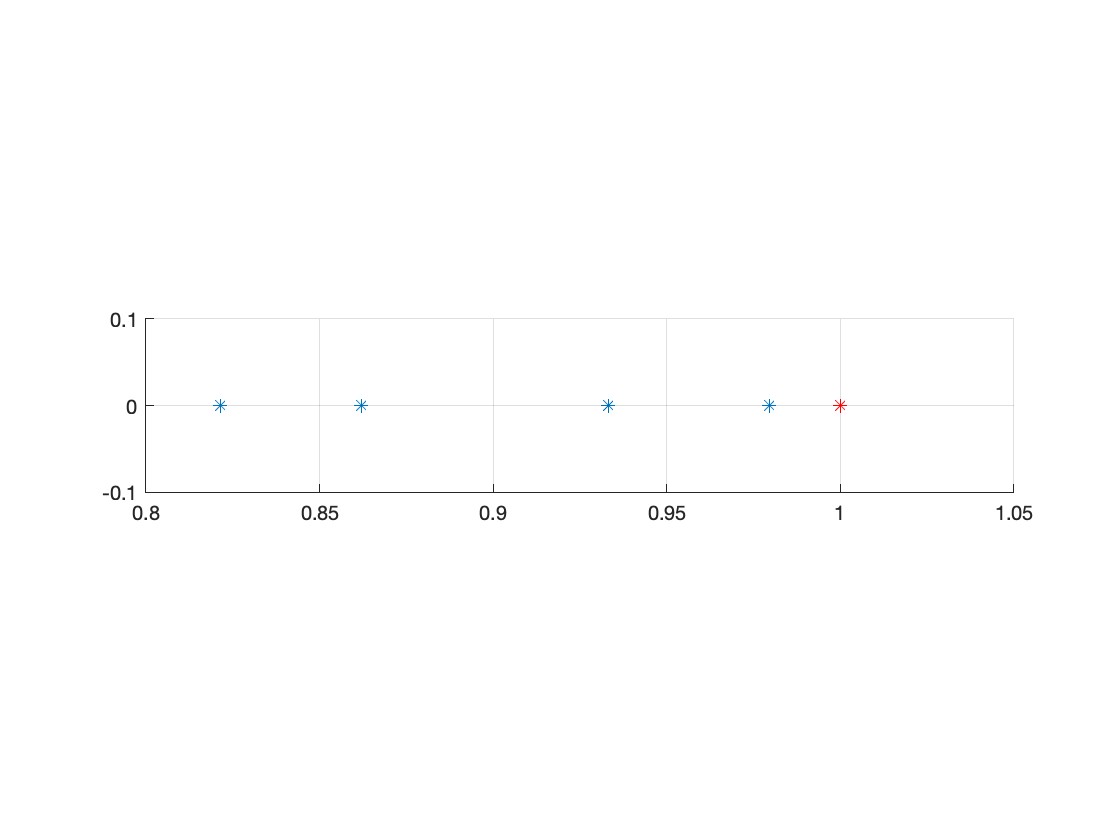}\\
					\vspace{-3.5mm}
					\centering
					\text{b. Approximated eigenvalues near $1$}
					%	\vspace{5mm}
				\end{minipage}
			\end{tabular}
			\begin{tabular}{cc}
				\begin{minipage}[t]{0.50\hsize}
					\centering
					\includegraphics[keepaspectratio, height=5cm]{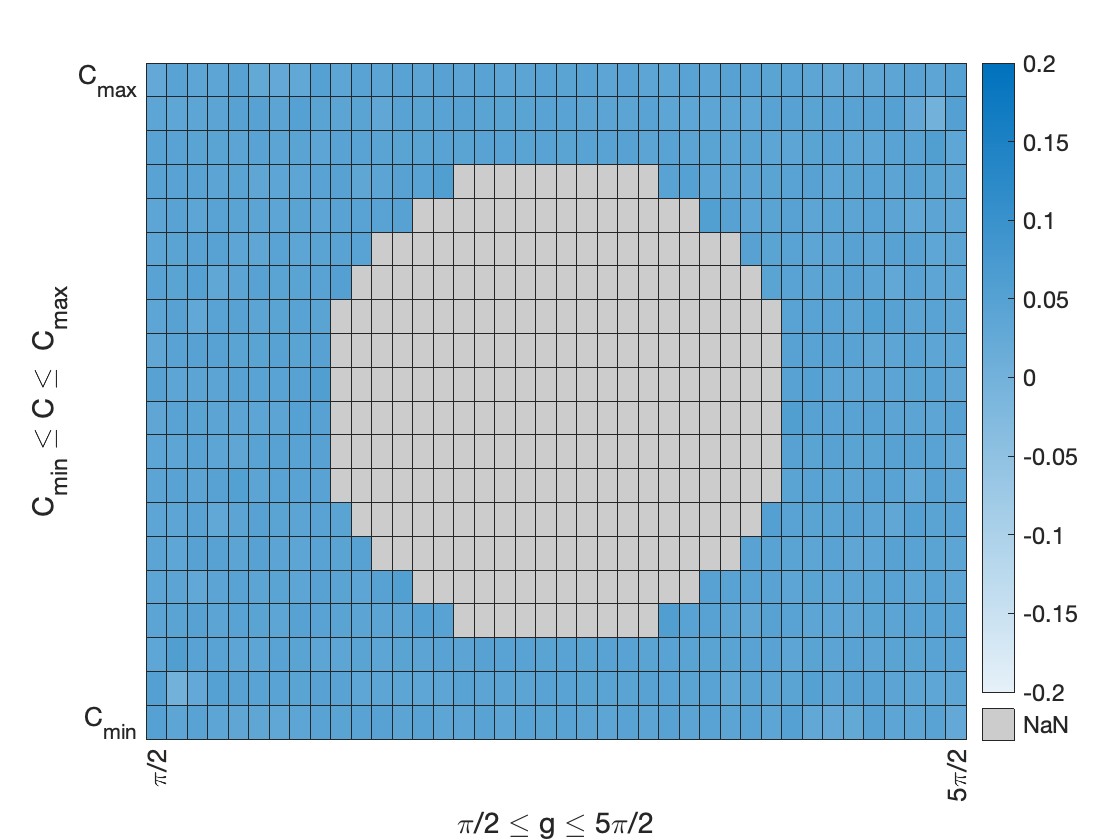}\\
					\centering
					\text{c. Eigenfunction for eigenvalue 1.00}
					\includegraphics[keepaspectratio, height=5cm]{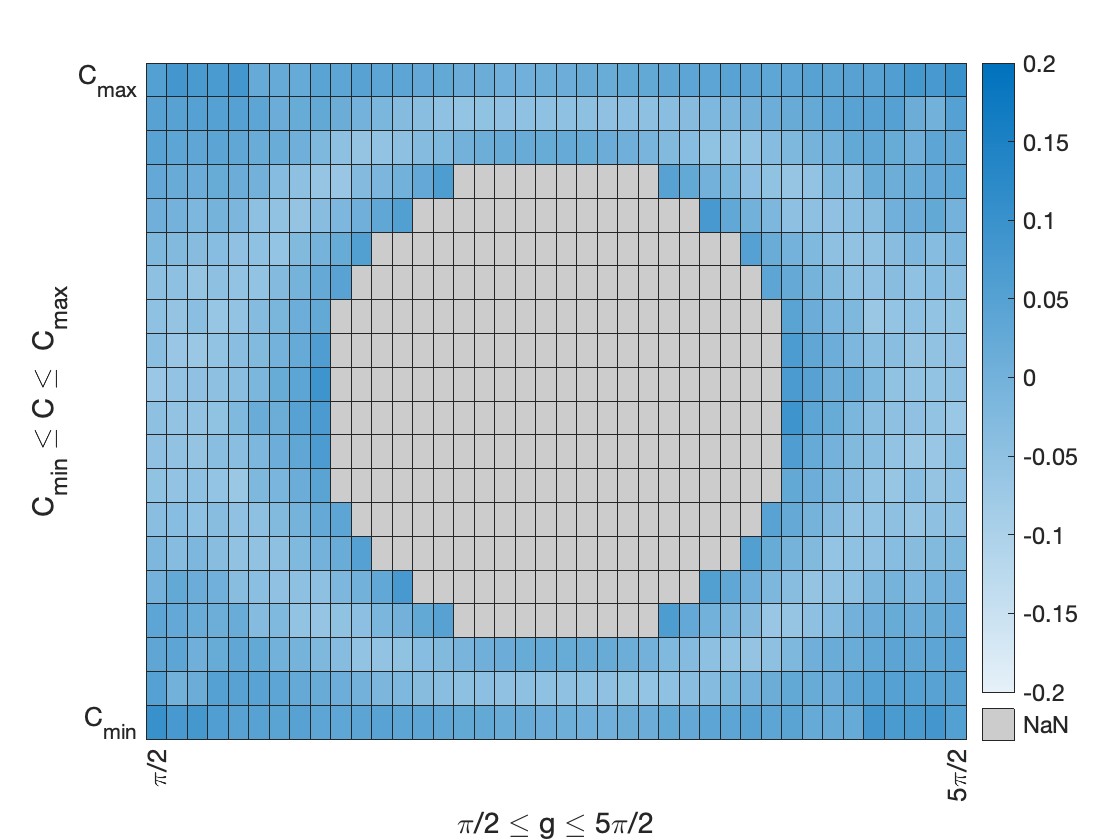}\\
					\centering
					\text{e. Eigenfunction for eigenvalue 0.93}
				\end{minipage} &
				\begin{minipage}[t]{0.50\hsize}
					\centering
					\includegraphics[keepaspectratio, height=5cm]{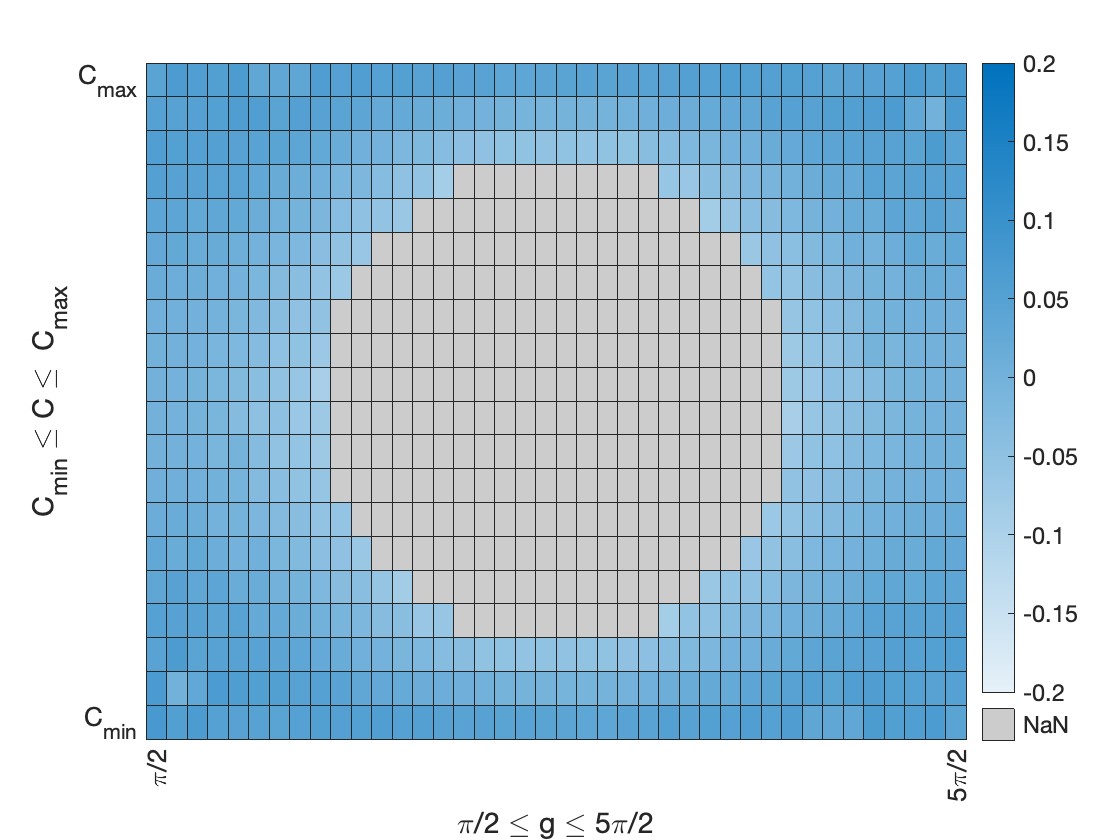}\\
					\centering
					\text{d. Eigenfunction for eigenvalue 0.98}
				\end{minipage}
			\end{tabular}
			\vspace{5mm}
			\centering
			%\text{c. level sets of independent eigenvectors of eigenvalue 1}
			\caption{Shifted allowed regions in $(g,C)-$coordinates space $[\pi/2, 5 \pi/2 ) \times [C_{min}, C_{max}]$.  Eigenvalues and eigenfunctions of approximated Koopman operator for $\alpha = 4.0, \beta = 0.0, N = 800, L = 25$, uniform weights.}
			\label{fig:eigen_beta000_gl}
		\end{figure}

		\begin{figure}
			\centering
			\begin{tabular}{cc}
				\begin{minipage}[t]{0.5\hsize}
					\includegraphics[keepaspectratio, height=5cm]{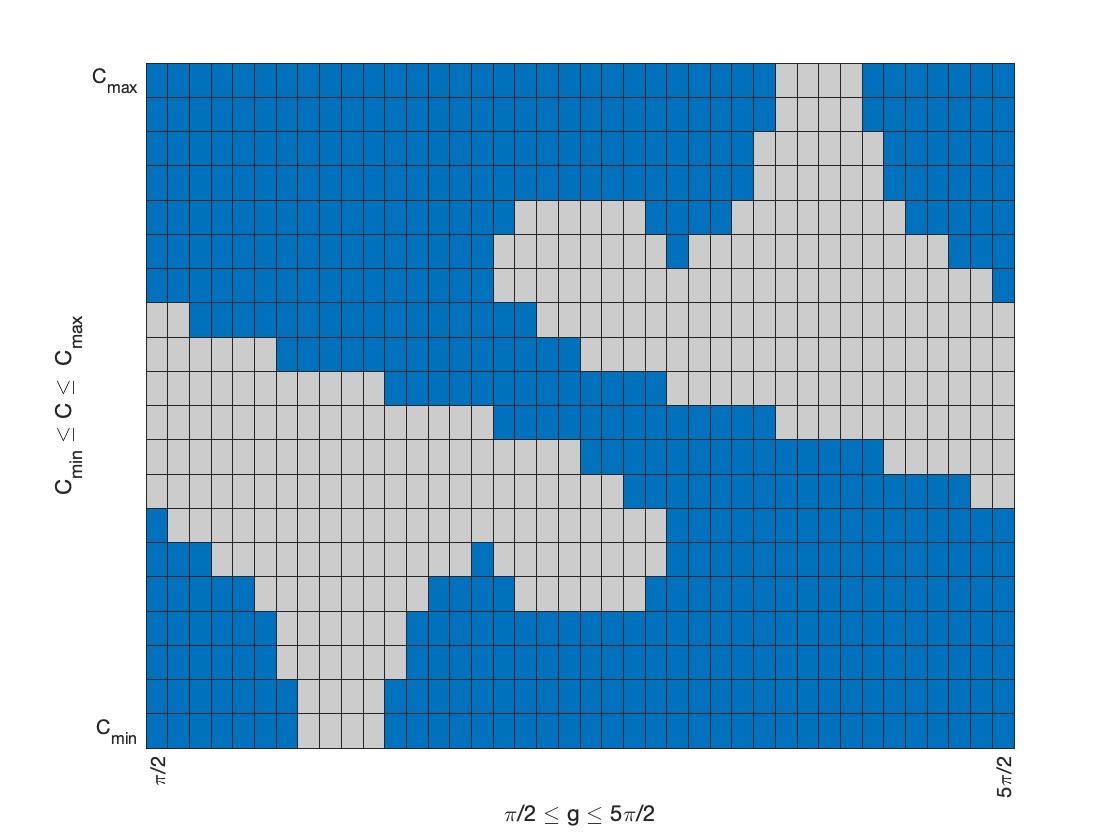}\\
					\centering
					\text{  a. Allowed regions (in blue)}
				\end{minipage}
				\begin{minipage}[t]{0.5\hsize}
					\includegraphics[keepaspectratio, height=5cm]{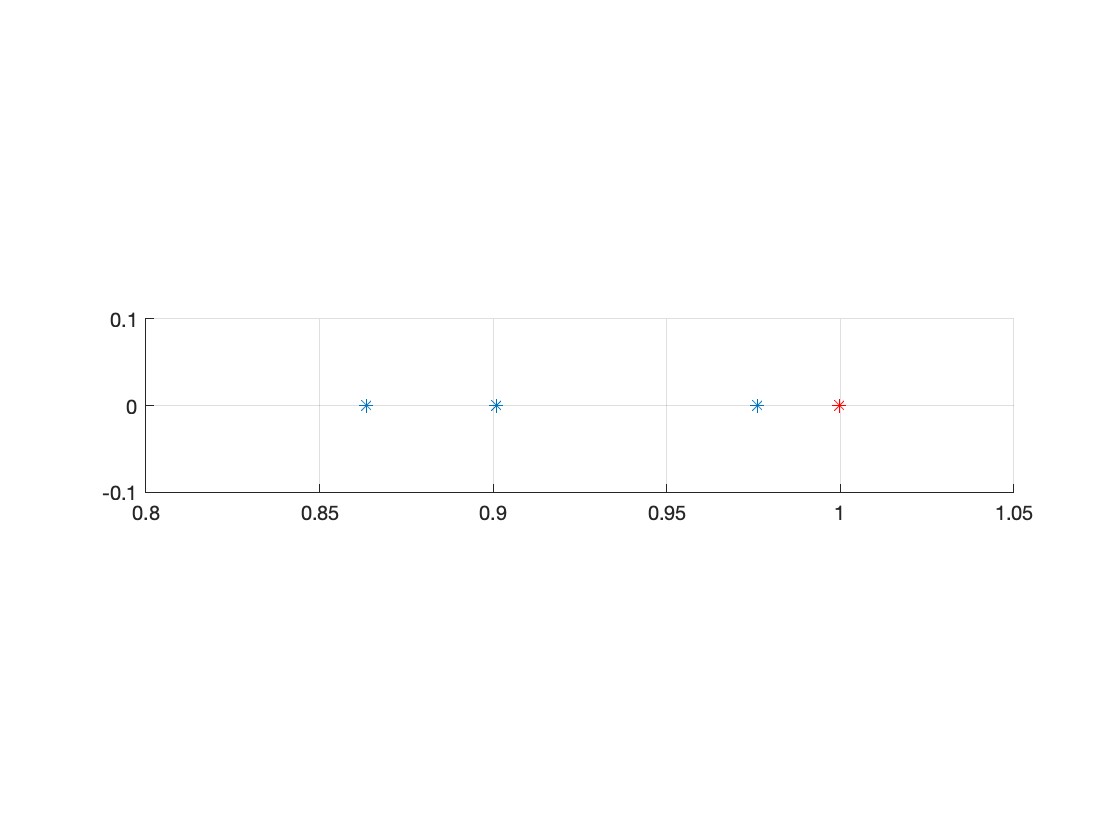}\\
					\vspace{-3.5mm}
					\centering
					\text{b. Approximated eigenvalues near $1$}
					%	\vspace{5mm}
				\end{minipage}
			\end{tabular}
			\begin{tabular}{cc}
				\begin{minipage}[t]{0.50\hsize}
					\centering
					\includegraphics[keepaspectratio, height=5cm]{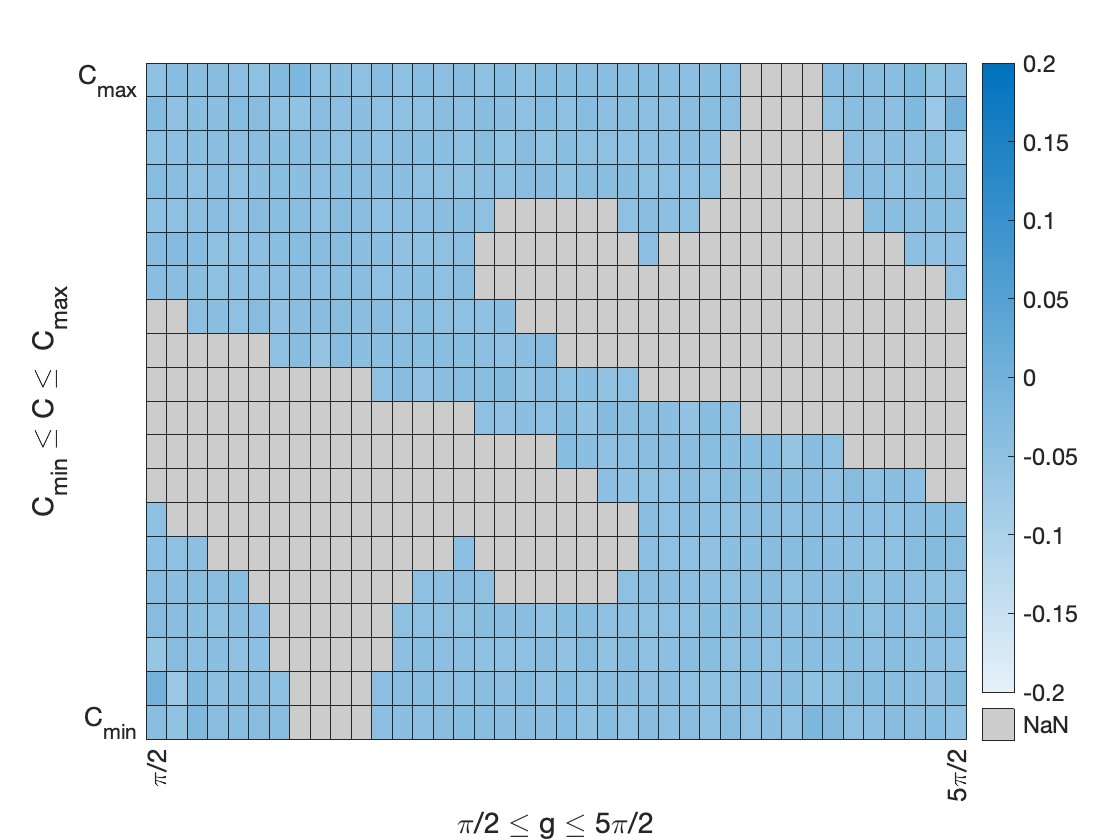}\\
					\centering
					\text{c. Eigenfunction for eigenvalue 1.00}
					\includegraphics[keepaspectratio, height=5cm]{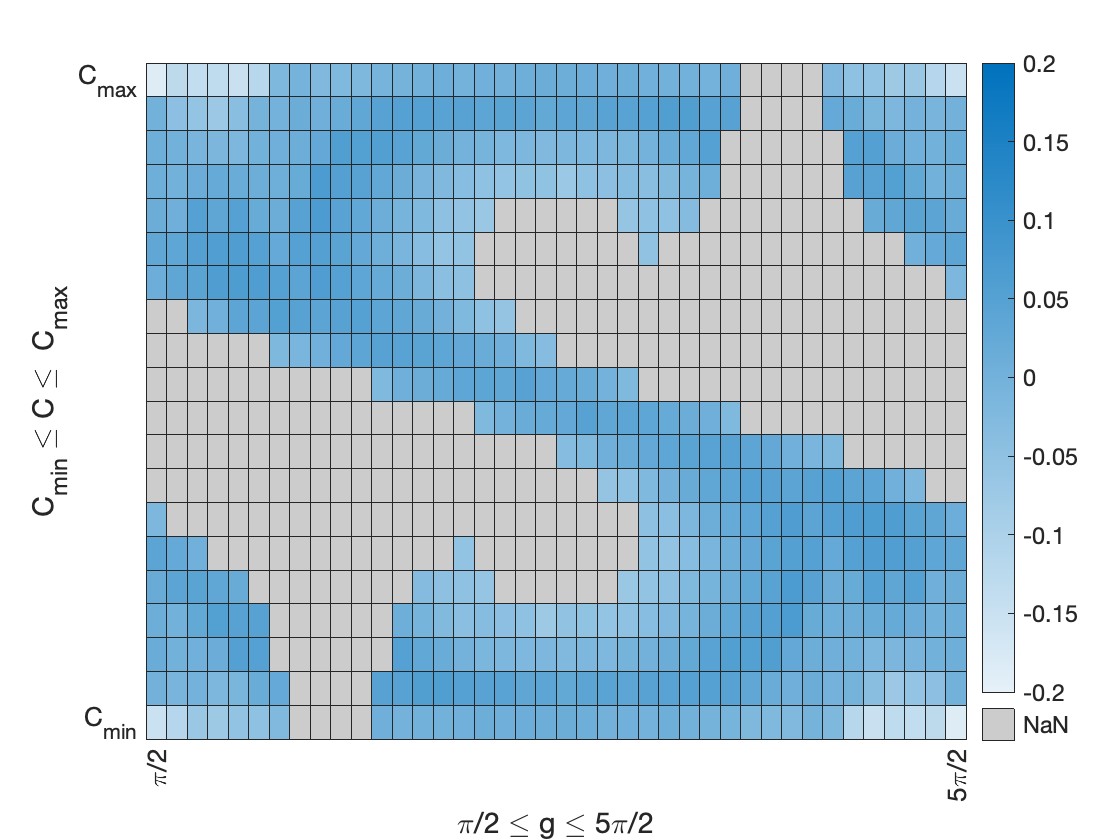}\\
					\centering
					\text{e. Eigenfunction for eigenvalue 0.90}
				\end{minipage} &
				\begin{minipage}[t]{0.50\hsize}
					\centering
					\includegraphics[keepaspectratio, height=5cm]{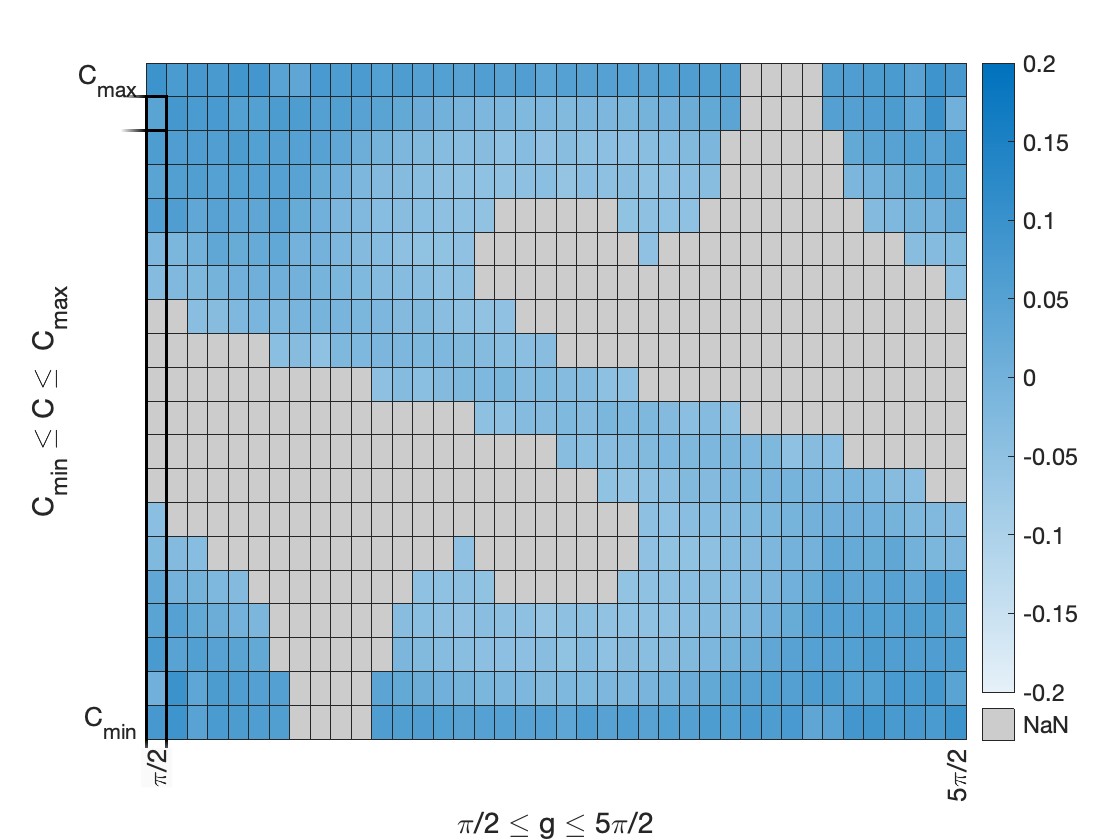}\\
					\centering
					\text{d. Eigenfunction for eigenvalue 0.98}
				\end{minipage}
			\end{tabular}
			\vspace{5mm}
			\centering
			%\text{c. level sets of independent eigenvectors of eigenvalue 1}
			\caption{Allowed regions in $(g,C)-$coordinates space $[0, 2 \pi ) \times [C_{min}, C_{max}]$.  Eigenvalues and eigenfunctions of approximated Koopman operator for $\alpha = 4.0, \beta = 0.5, N = 800, L = 25$, uniform weights.}
			\label{fig:eigen_beta050_gl}
		\end{figure}

		\begin{figure}
			\centering
			\begin{tabular}{cc}
				\begin{minipage}[t]{0.5\hsize}
					\includegraphics[keepaspectratio, height=5cm]{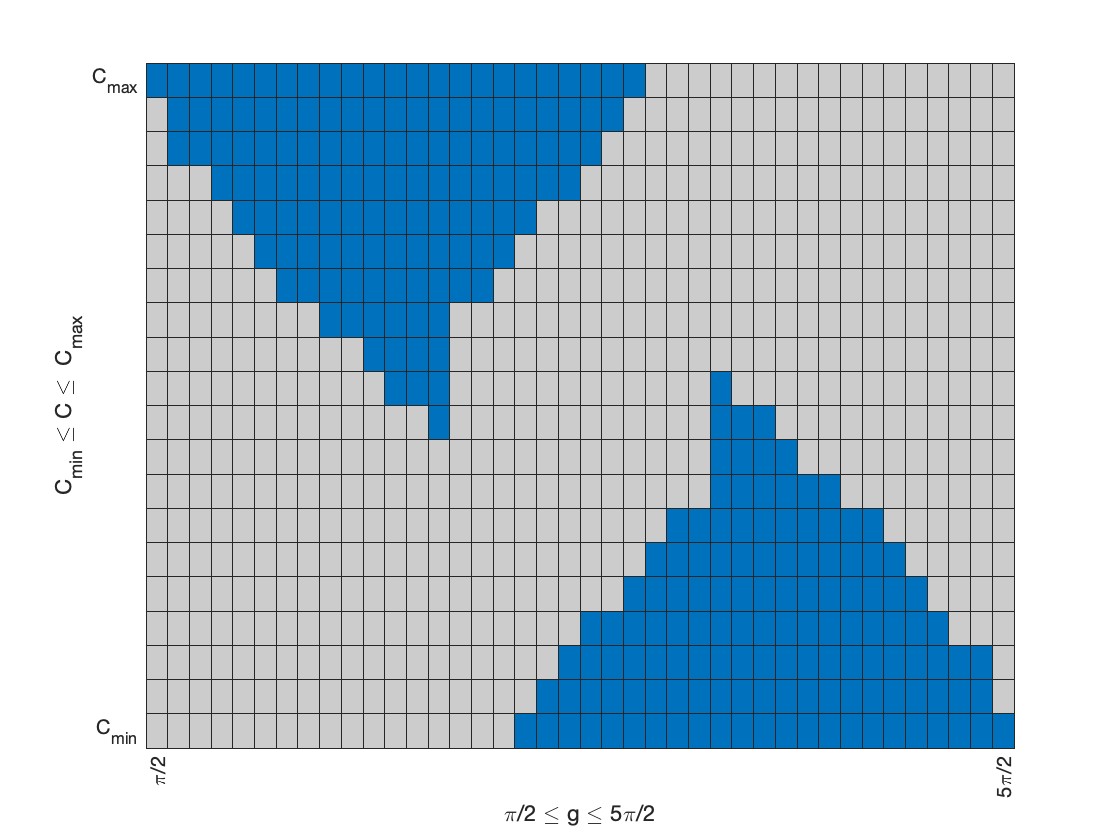}\\
					\centering
					\text{ a. Allowed regions (in blue)}
				\end{minipage}
				\begin{minipage}[t]{0.5\hsize}
					\includegraphics[keepaspectratio, height=5cm]{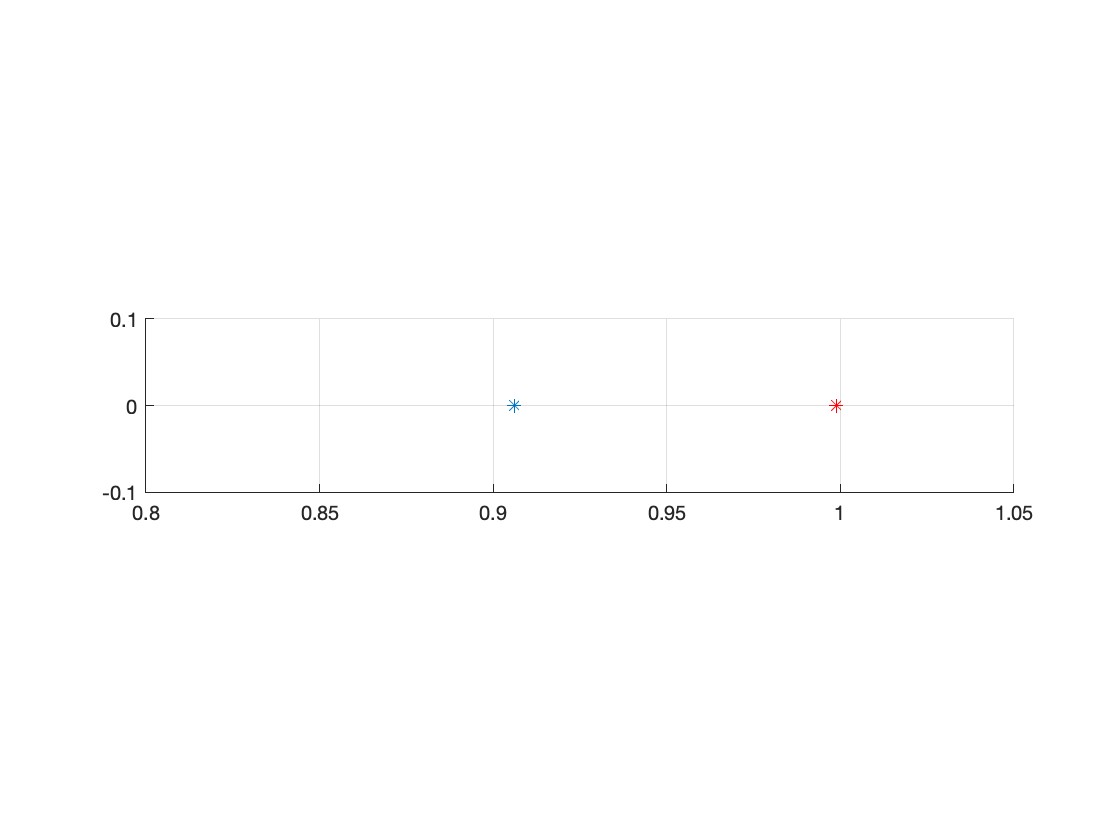}\\
					\vspace{-3.5mm}
					\centering
					\text{b. Approximated eigenvalues near $1$}
					%	\vspace{5mm}
				\end{minipage}
			\end{tabular}
			\begin{tabular}{cc}
				\begin{minipage}[t]{0.50\hsize}
					\centering
					\includegraphics[keepaspectratio, height=5cm]{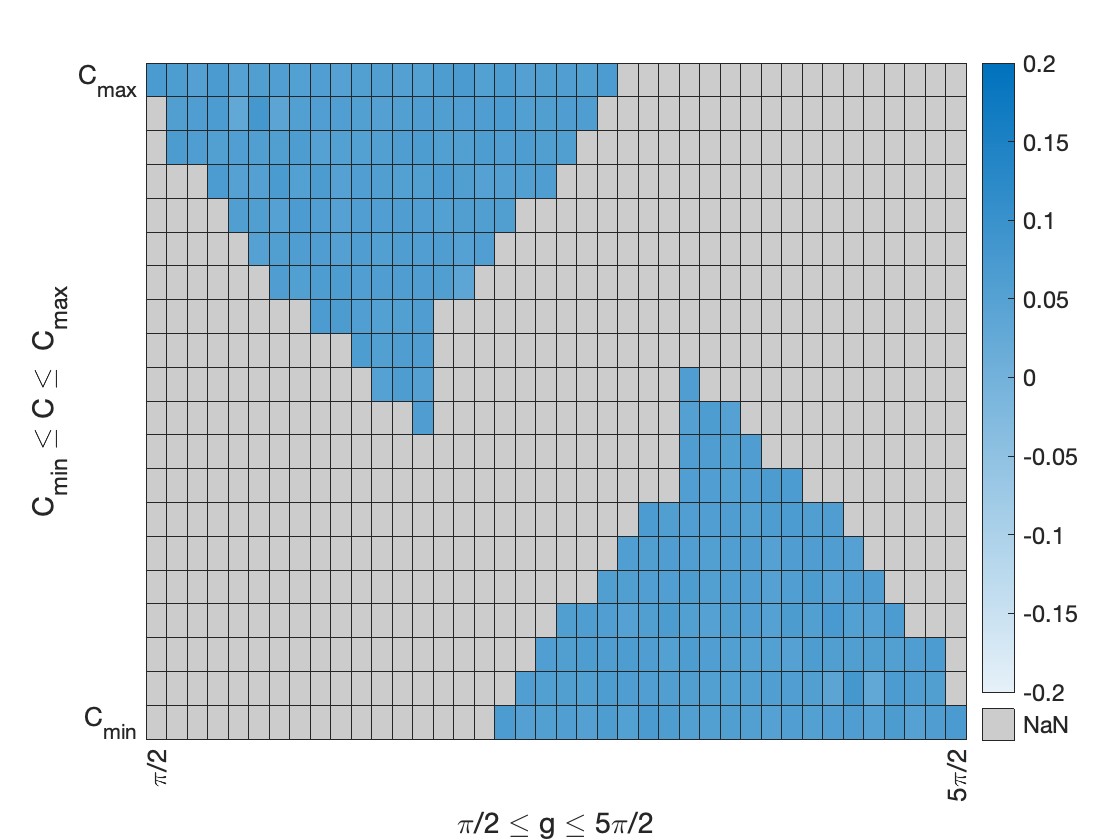}\\
					\centering
					\text{c. Eigenfunction for eigenvalue 1.00}
				\end{minipage} &
				\begin{minipage}[t]{0.50\hsize}
					\centering
					\includegraphics[keepaspectratio, height=5cm]{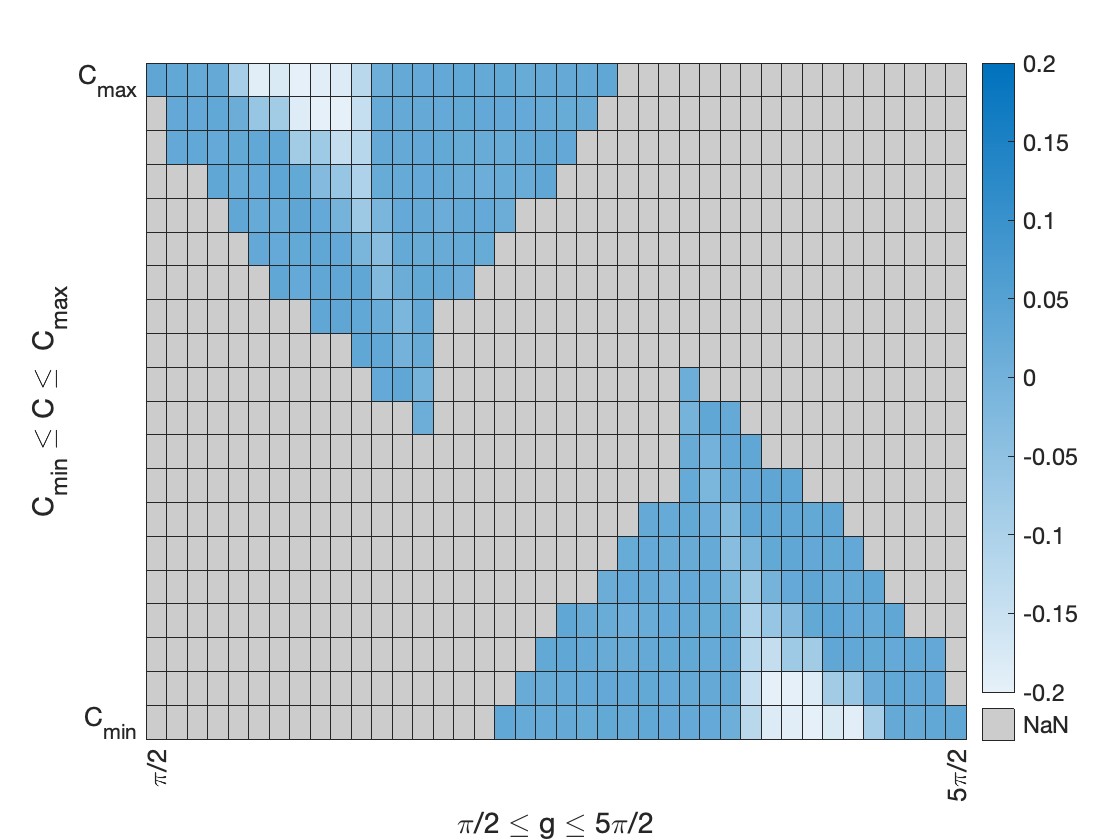}\\
					\centering
					\text{d. Eigenfunction for eigenvalue 0.91}
				\end{minipage}
			\end{tabular}
			\vspace{5mm}
			\centering
			%\text{c. level sets of independent eigenvectors of eigenvalue 1}
			\caption{Allowed regions in $(g,C)-$coordinates space $[0, 2 \pi ) \times [C_{min}, C_{max}]$.  Eigenvalues and eigenfunctions of approximated Koopman operator for $\alpha = 4.0, \beta = 2.4, N = 800, L = 25$, uniform weights.}
			\label{fig:eigen_beta240_gl}
		\end{figure}
		
		\begin{figure}
			\centering
			
			\begin{tabular}{cc}
				\begin{minipage}[t]{0.4\hsize}
					\hspace{-2.6cm}\includegraphics[keepaspectratio, width=6cm]{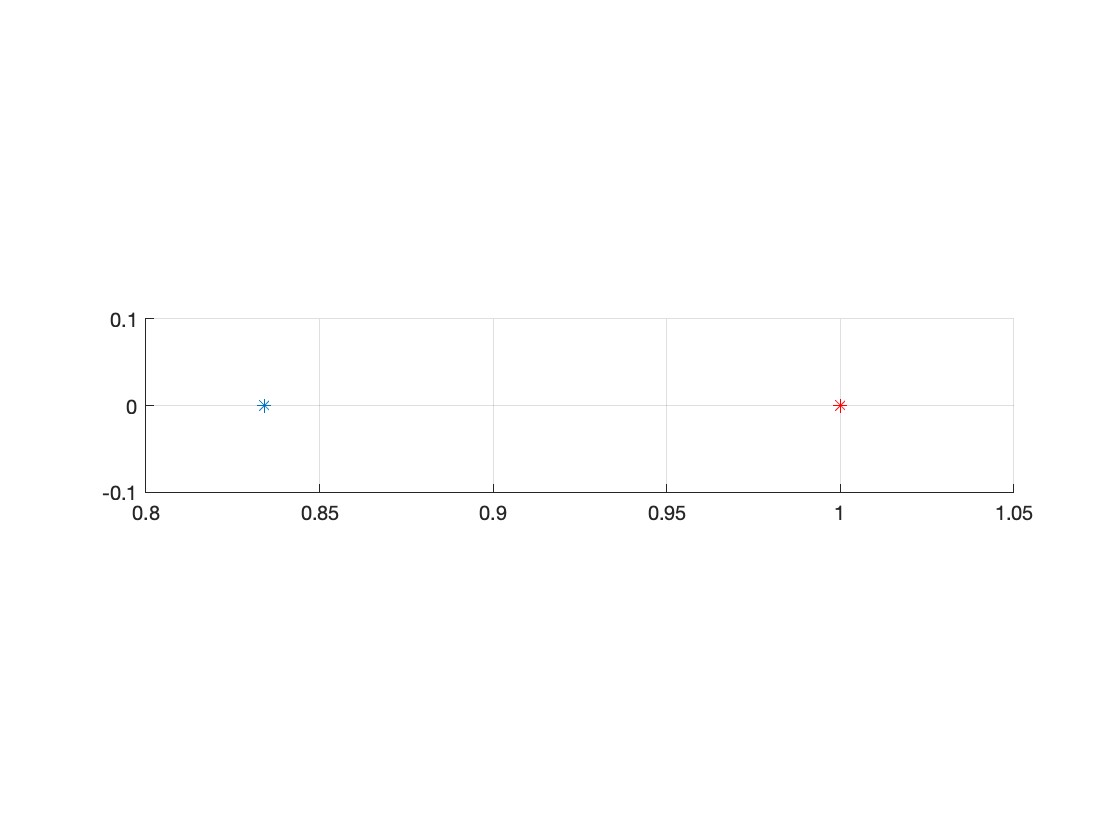}
				\end{minipage}
				\begin{minipage}[t]{0.4\hsize}
					\hspace{-2.2cm}\includegraphics[keepaspectratio, width=6cm]{beta240_1_1_40_20_5_5_40_20_eigenv.jpg}
					%	\vspace{5mm}
				\end{minipage}
			\begin{minipage}[t]{0.4\hsize}
			\hspace{-1.5cm}\includegraphics[keepaspectratio, width=6cm]{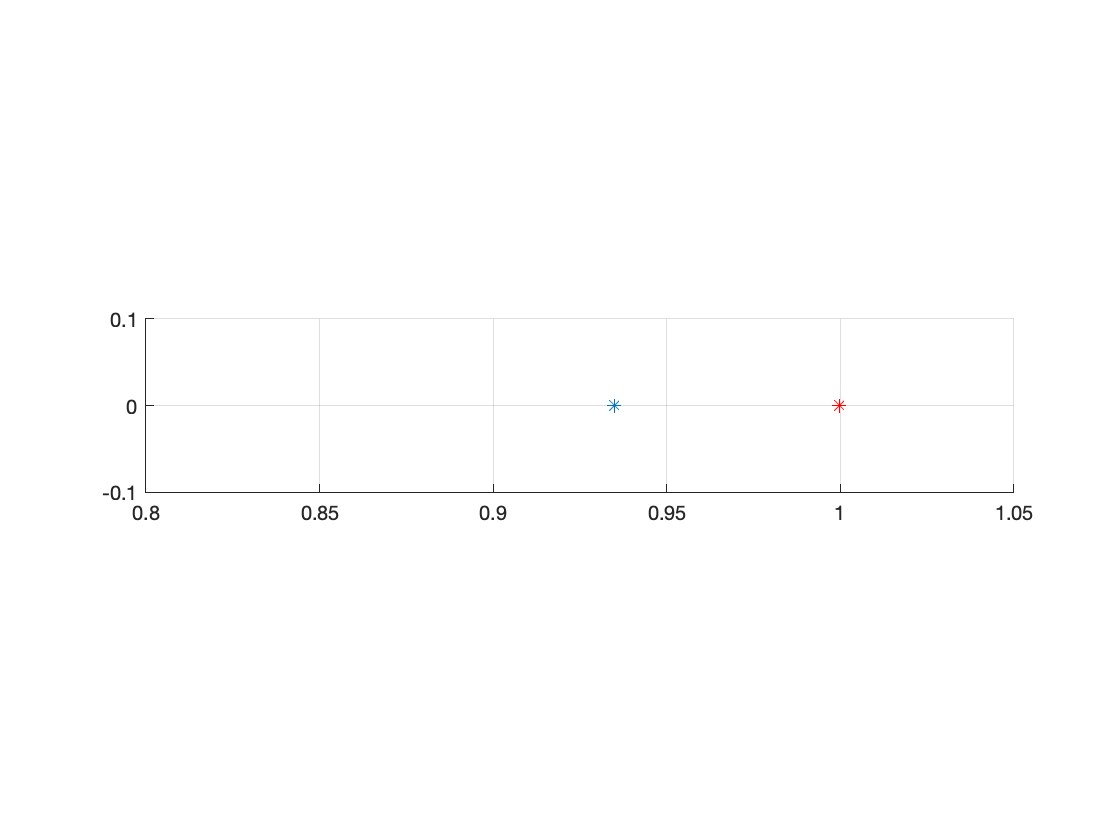}
			%	\vspace{5mm}
		\end{minipage}
			\end{tabular}
		\caption{Approximated eigenvalues near $1$ for $\alpha = 4.0, \beta = 2.4, L = 25$, $N=200(\text{left}), 800(\text{middle}), 1800(\text{right}),$ uniform weights.}
		\end{figure}

		\begin{figure}
			\centering
			\begin{tabular}{cc}
				\begin{minipage}[t]{0.5\hsize}
					\includegraphics[keepaspectratio, height=5cm]{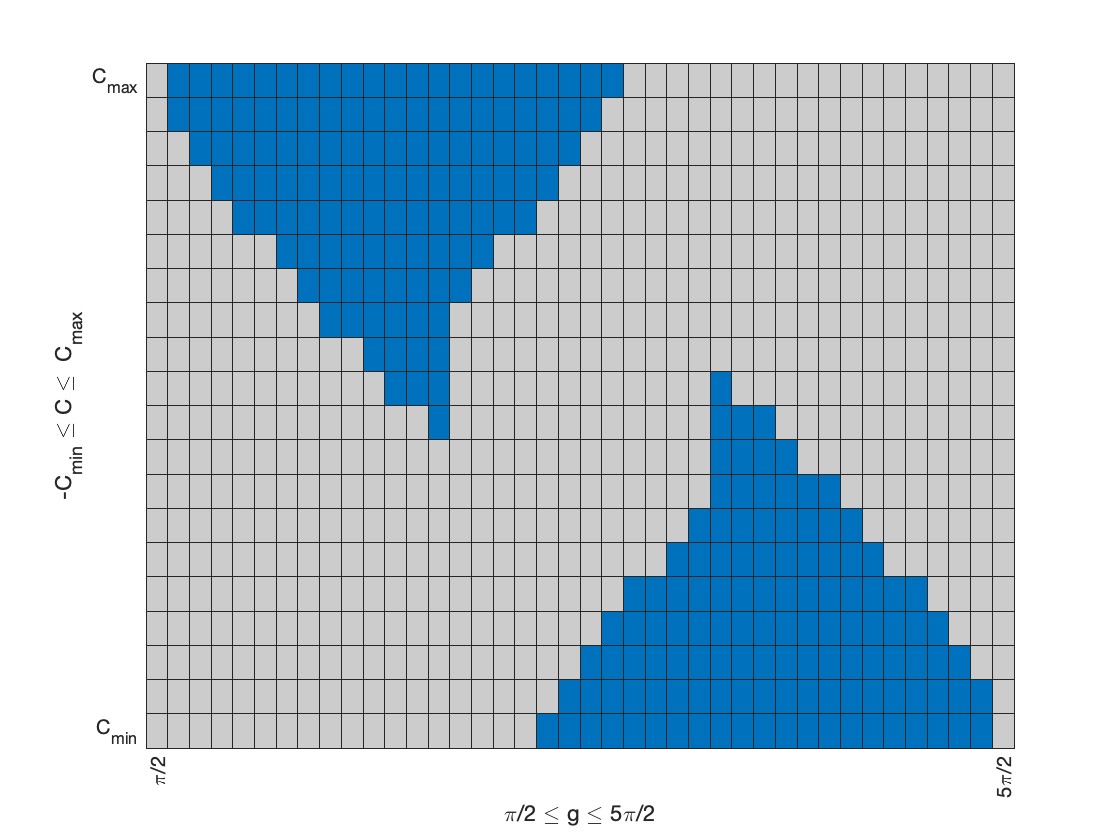}\\
					\centering
					\text{ a. Allowed regions (in blue)}
					\centering
					\includegraphics[keepaspectratio, height=5cm]{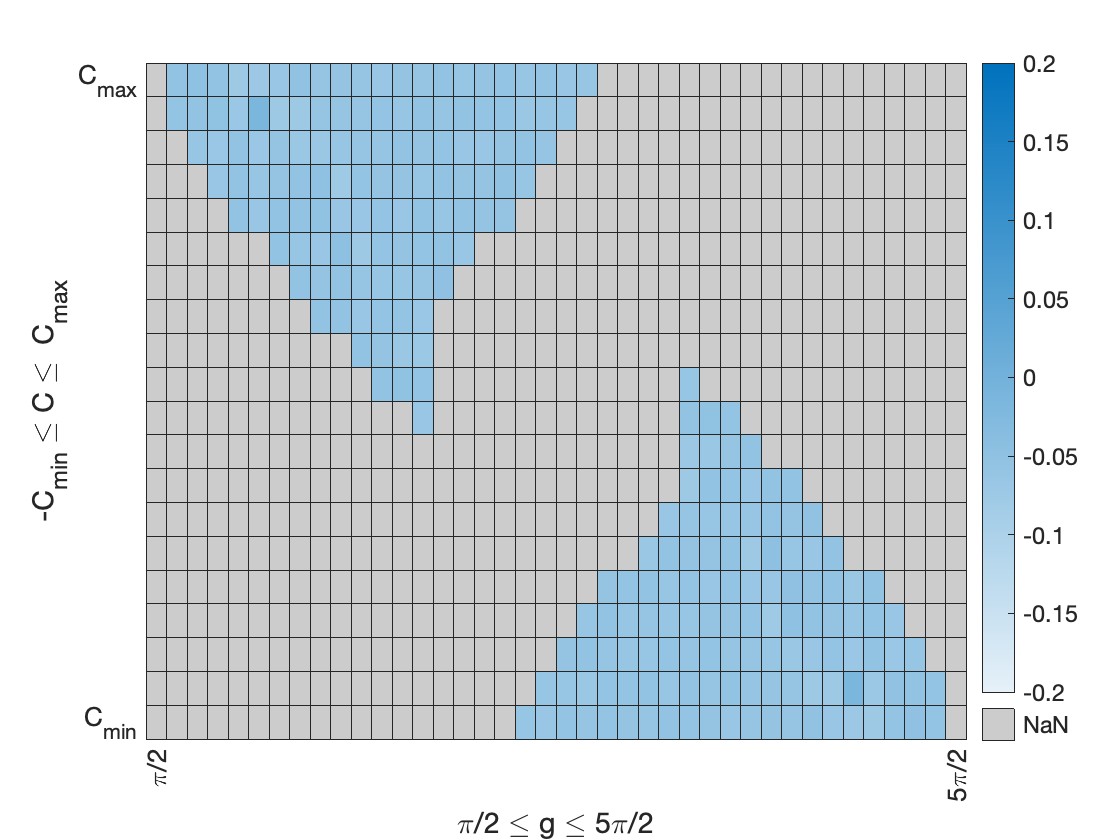}\\
					\centering
					\text{c. Eigenfunction for eigenvalue 1.00}
				\end{minipage}
				\begin{minipage}[t]{0.5\hsize}
					\includegraphics[keepaspectratio, height=5cm]{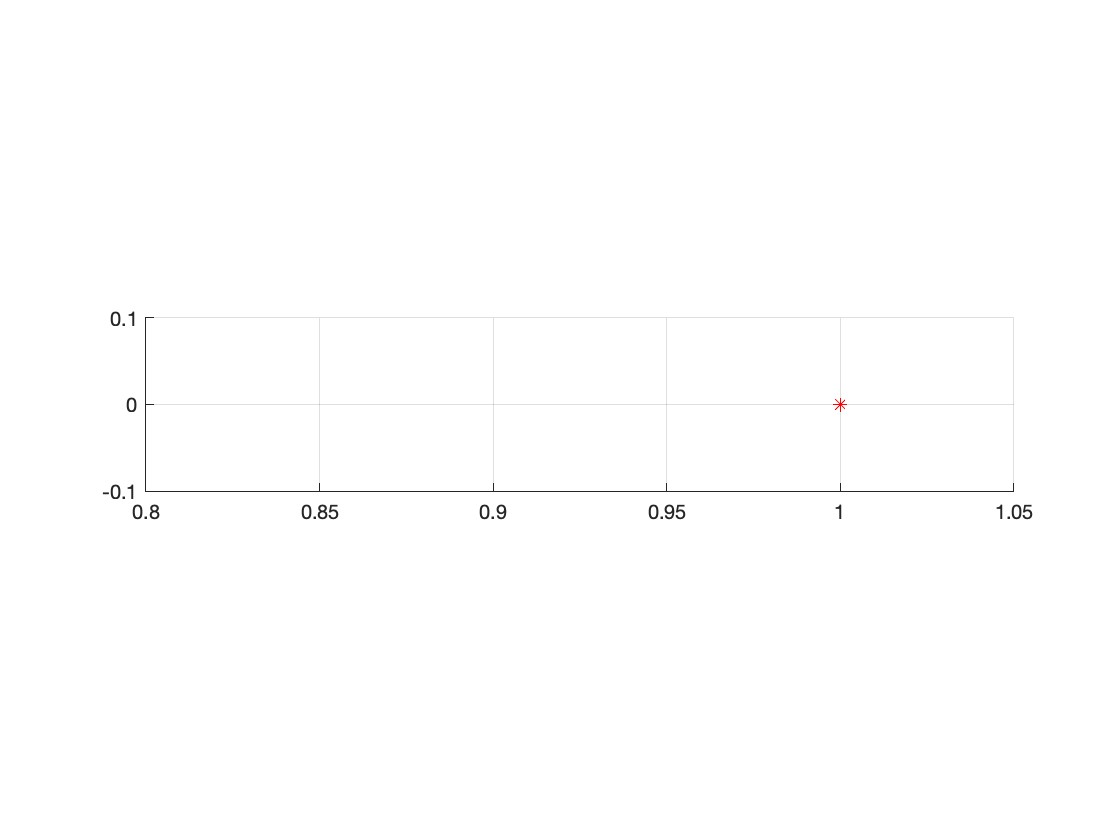}\\
					\vspace{-3.5mm}
					\centering
					\text{b. Approximated eigenvalues near  $1$}
					%	\vspace{5mm}
				\end{minipage}
			\end{tabular}
			\vspace{5mm}
			\centering
			%\text{c. level sets of independent eigenvectors of eigenvalue 1}
			\caption{Allowed regions in $(g,C)-$coordinates space $[0, 2 \pi ) \times [C_{min}, C_{max}]$.   Eigenvalues and eigenfunctions of approximated Koopman operator for $\alpha = 4.0, \beta = 2.6, N = 800, L = 25$, uniform weights.}
			\label{fig:eigen_beta260_gl}
		\end{figure}

		\paragraph{Galerkin method and Gauss-Legendre quadrature}
		The Gauss-Legendre quadrature approximates the integral of the function $f$ in the 
		domain $[-1, 1]$ with the sum of the values of the function at the Gauss points $\{x_k\}$, with the appropriate weights $\{w_k  \}$, as
		\begin{equation*}
			\int_{-1}^{1} f(x) dx \approx \sum_{k = 1}^{K } w_k f(x_k).
		\end{equation*}
		The Gauss node points can be defined as the roots of {the} Legendre polynomials 
		\begin{equation*}
			P_K(x) = \frac{1}{2^K K !} \frac{d^K}{d x^K}(x^2 -1)^K,
		\end{equation*}
		and the weights are {assigned as}:
		\begin{equation*}
			w_k = \frac{2}{(1- x_k^2)[ P_K' (x_k)]^2}.
		\end{equation*}
		The Gauss-Legendre quadrature can be extended to integration over {a surface as}:
		\begin{equation*}
			\int_{-1}^{1} \int_{-1}^{1} f(x,y) dx_1 dx_2 \approx \sum_{k_1 = 1}^{K } \sum_{k_2 = 1}^{K } w_{k_1 k_2} f(x_k, x_l),
		\end{equation*}
		where $w_{k_1 k_2}= w_{k_1} w_{k_2}$.
		In this way, we approximate each entry of the matrices in (\ref{eq: mat_eigen}) {and get}
		\begin{align}
			\label{eq: int_approx3}
			\begin{split}
				\langle f_n \circ S, f_m \rangle &= \int_{X} f_n(S(x)) \cdot f_m(x) dx\\
				& = \int_{\Omega_m} f_n(S(x)) dx \\
				& \approx \sum_{k_1 = 1}^{K} \sum_{k_2 = 1}^{K}  a_m w_{k_1 k_2}^{(m)}  f_n(S(x_{k_1 k_2})\\
				&= \sum_{\ell = 1}^{L= K \times K} w_{\ell}^{(m)} f_n(S(x_\ell)) \\
				& = \sum_{\ell  ~\text{s.t.} S(x_\ell) \in \Omega_n}  w_\ell^{(m)},
			\end{split}
		\end{align}
		where $a_m$ is the area of $\Omega_m$ and $w^{(m)}_{\ell}:=a_m  w^{(m)}_{k_1 k_2}$
		and 
		\begin{align}
			\label{eq: int_approx4}
			\begin{split}
				\langle f_n , f_m \rangle  & = \int_{X} f_n(x) \cdot f_m(x) dx \\
				& = \int_{\Omega_m} f_n (x) dx\\
				& \approx \sum_{\ell =1}^{L} w_{\ell}^{(m)} f _n(x_\ell)\\
				& = \sum_{\ell  ~\text{s.t.} x_\ell \in \Omega_n}  w_\ell^{(m)}.
			\end{split}
		\end{align}
		{Remind} that $f_n(x_\ell) = 1 $ if $x_\ell \in \Omega_n$ and $f_n(x_\ell) = 0$ otherwise.

		We again set $X= [0, 2 \pi ) \times [C_{min}, C_{max}]$, $x= (g,C)$ and consider the approximated eigenvalue problem of the corresponding Koopman operator with the Galerkin method using {the} Gauss-Legegendre quadrature. 
		
		In the following numerical computations, we divided 
		%the $(g,C)-$coordinate space $[0, 2 \pi ) \times [C_{min}, C_{max}]$
		{$X$} into $N=800$ partial sets. Recall that we need to restrict the space $[0, 2 \pi ) \times [C_{min}, C_{max}]$  into the subset where  the billiard mapping is well-defined.  In our computations, we set $\alpha = 4.0, E= -0.5, \gamma= 0.5$, and vary the parameter $\beta$. The number $L=25$ represents the number of Gauss nodes in each partition used to approximate integrals, which appear in the equations \eqref{eq: int_approx3} and \eqref{eq: int_approx4}.

		In the following figures, we illustrate the numerical results.
		% on the discretized eigenvalue problem \eqref{eq: mat_eigen} approximated by Galerkin method with Gauss-Legendre quadrature as is described above.

		In Figure \ref{fig:eigen_gl_beta000_gl}, Subfig. a shows the restricted region in a divided phase space $[0, 2 \pi ) \times [C_{min}, C_{max}]$ ($N =  800$)  where the billiard mapping is well-defined for $\alpha = 4.0, \beta = 0.0$, Subfig. b shows the all eigenvalues of the discretized Koopman eigenvalue problem, {and} Subfig. c, d, e, and f show the level sets of all independent eigenfunctions corresponding to the three closest eigenvalues from 1. 
		Figure \ref{fig:eigen_gl_beta050_gl}, Figure \ref{fig:eigen_gl_beta240_gl} and Figure \ref{fig:eigen_gl_beta260_gl} show the numerical results for $\beta = 0.5, \beta = 2.4$ and $\beta = 2.6$, respectively.
		% same information on the discretized Koopman eigenvalue problem as Figure \ref{fig:eigen_gl_beta000_gl} but for  different parameter values  

		%Figure \ref{fig:eigen_gl_beta000_gl} shows the divided restricted state space ($N =  662 $) for $\alpha = 4.0. \beta = 0.0$ and eigenfunctions corresponding to three closest real eigenvalues from 1.
		%Figure \ref{fig:eigen_gl_beta240_gl_1} shows the divided restricted state space ($N =  260 $) for $\alpha = 4.0. \beta = 0.0$ and eigenfunctions corresponding to three closest real eigenvalues from 1. 
		%Figure \ref{fig:eigen_gl_beta240_gl_2} shows the divided restricted state space ($N =  408 $) for $\alpha = 4.0. \beta = 2.4$ and eigenfunctions corresponding to three closest real eigenvalues from 1. 
		%Figure \ref{fig:eigen_gl_beta260_gl_1} shows the divided restricted state space ($N =  384 $) for $\alpha = 4.0. \beta = 2.6$ and eigenfunctions corresponding to two closest real eigenvalues from 1. 
		%Figure \ref{fig:eigen_gl_beta260_gl_2} shows the divided restricted state space ($N =  586 $) for $\alpha = 4.0. \beta = 2.6$ and eigenfunction corresponding to the closest real eigenvalues from 1. 

		\begin{figure}
			\centering
			\begin{tabular}{cc}
				\begin{minipage}[t]{0.5\hsize}
					\includegraphics[keepaspectratio, height=5cm]{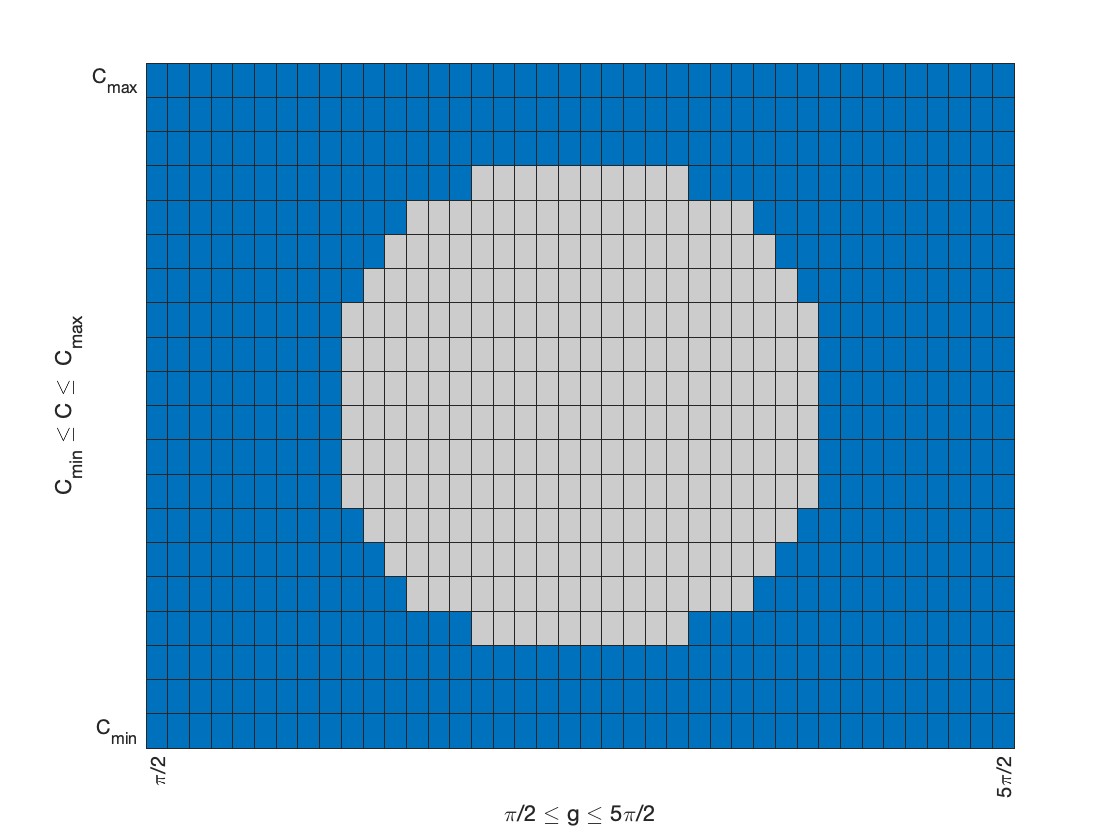}\\
					\centering
					\text{  a. Allowed regions (in blue)}
				\end{minipage}
				\begin{minipage}[t]{0.5\hsize}
					\includegraphics[keepaspectratio, height=5cm]{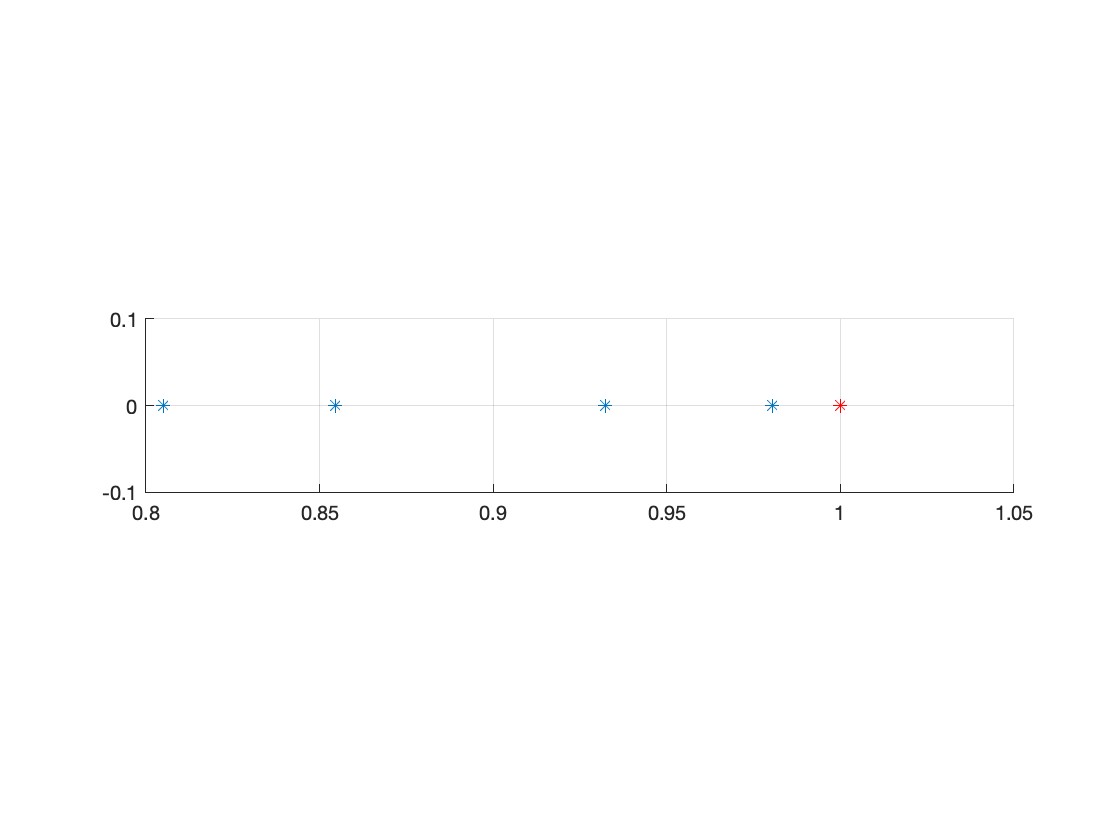}\\
					\vspace{-3.5mm}
					\centering
					\text{b. Approximated eigenvalues near $1$}
					%	\vspace{5mm}
				\end{minipage}
			\end{tabular}
			\begin{tabular}{cc}
				\begin{minipage}[t]{0.50\hsize}
					\centering
					\includegraphics[keepaspectratio, height=5cm]{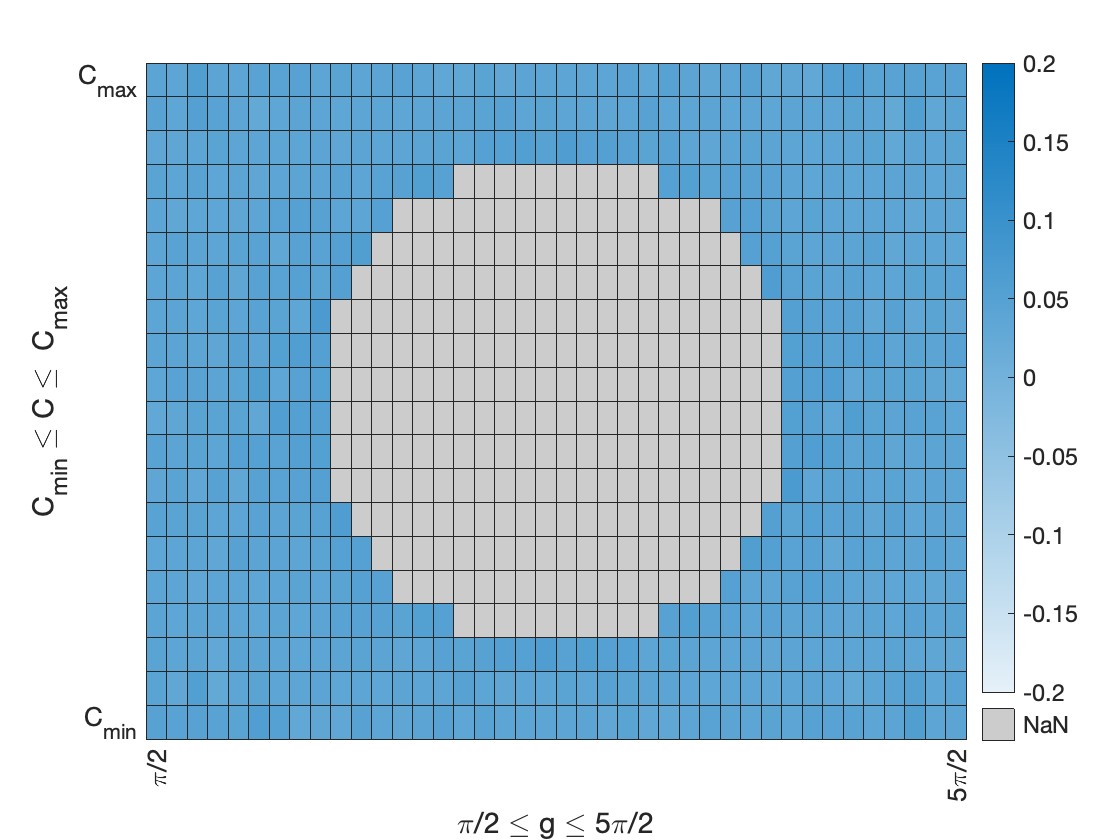}\\
					\centering
					\text{c. Eigenfunction for eigenvalue 1.00}
					\includegraphics[keepaspectratio, height=5cm]{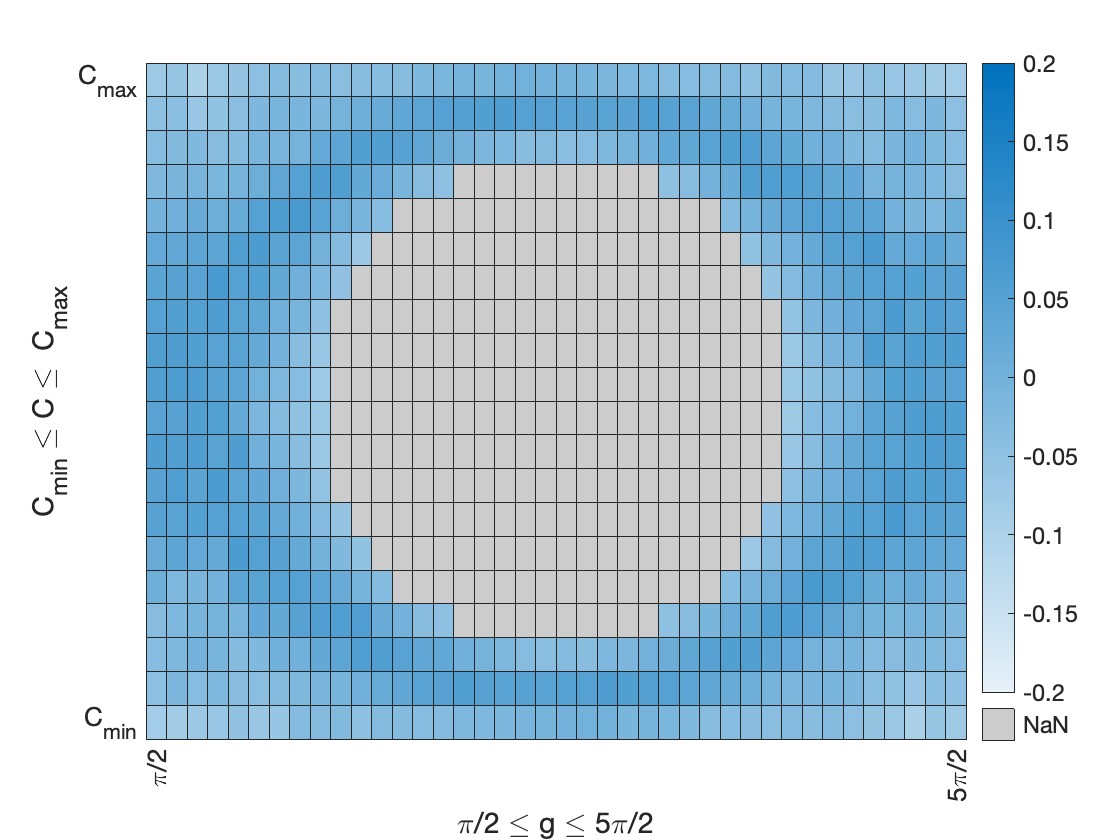}\\
					\centering
					\text{e. Eigenfunction for eigenvalue 0.93}
				\end{minipage} &
				\begin{minipage}[t]{0.50\hsize}
					\centering
					\includegraphics[keepaspectratio, height=5cm]{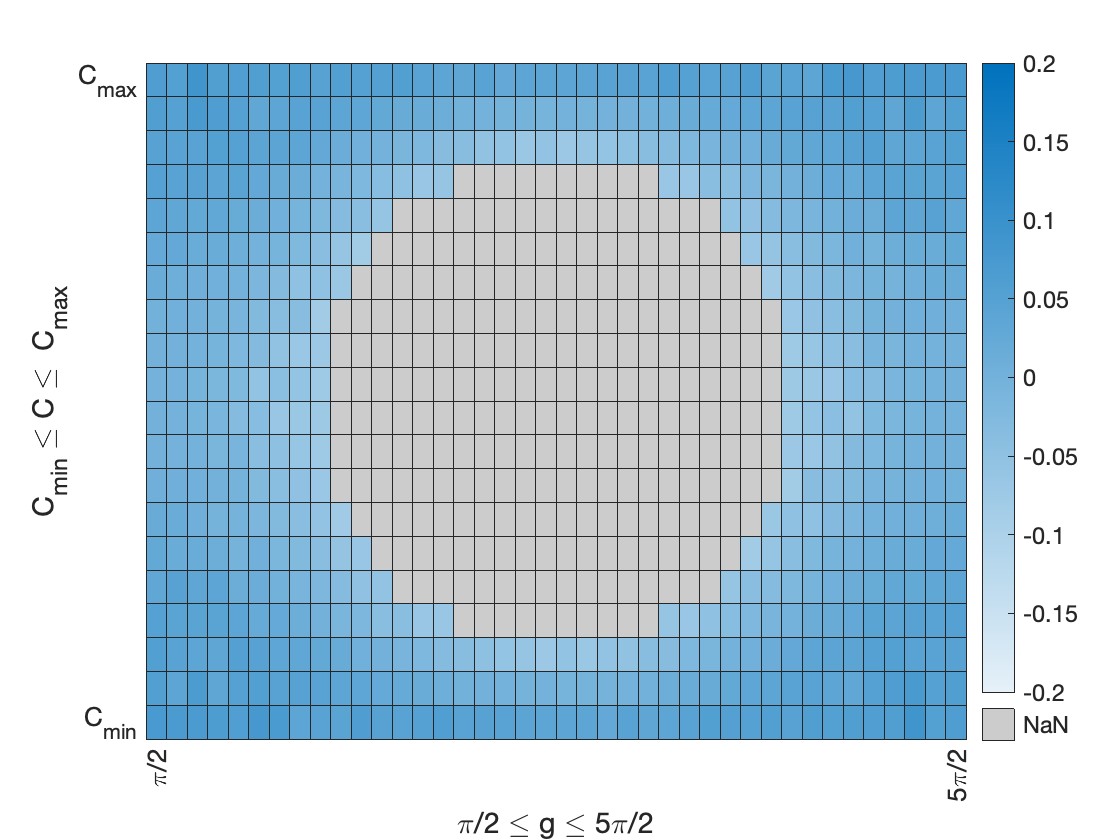}\\
					\centering
					\text{f. Eigenfunction for eigenvalue 0.98}
				\end{minipage}
			\end{tabular}
			\vspace{5mm}
			\centering
			%\text{c. level sets of independent eigenvectors of eigenvalue 1}
			\caption{Allowed regions in $(g,C)-$coordinates space $[0, 2 \pi ) \times [C_{min}, C_{max}]$.  Eigenvalues and eigenfunctions of approximated Koopman operator for $\alpha = 4.0, \beta = 0.0, N = 800, L = 25$, Gauss-Legendre quadrature.}
			\label{fig:eigen_gl_beta000_gl}
		\end{figure}

		\begin{figure}
			\centering
			\begin{tabular}{cc}
				\begin{minipage}[t]{0.5\hsize}
					\includegraphics[keepaspectratio, height=5cm]{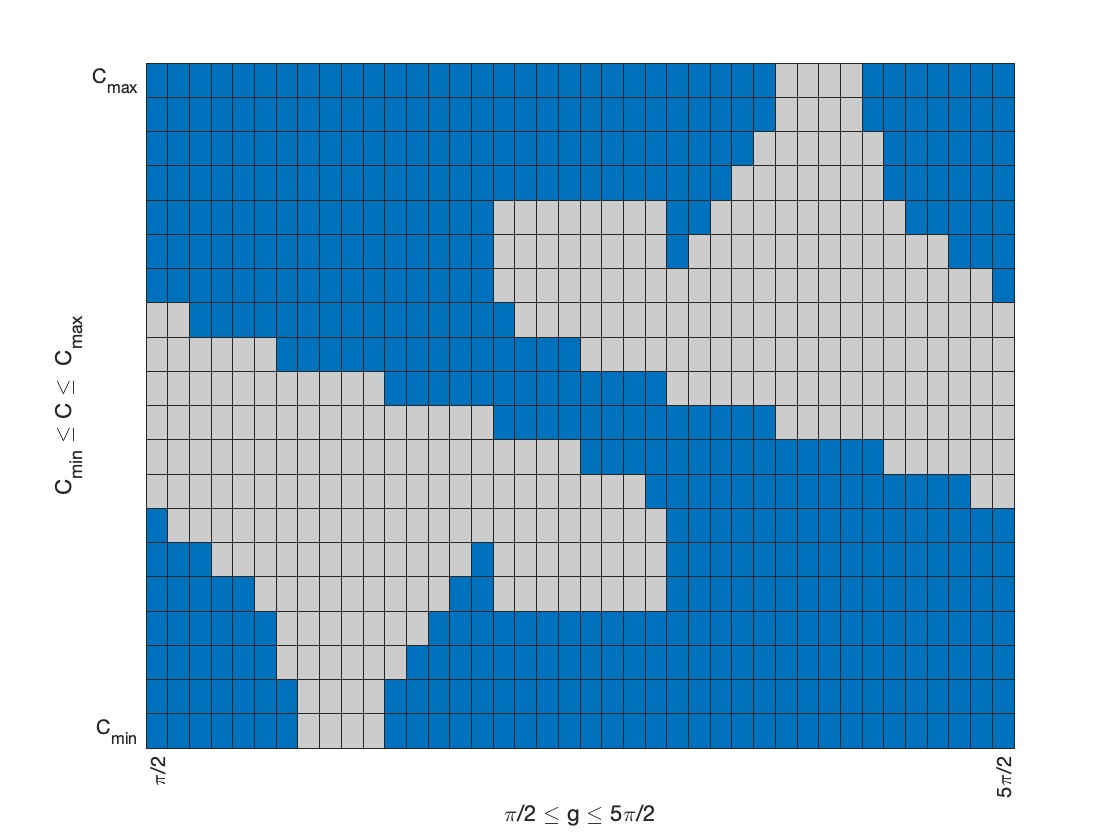}\\
					\centering
					\text{ \hspace{-10mm} a. Allowed regions (in blue)}
				\end{minipage}
				\begin{minipage}[t]{0.5\hsize}
					\includegraphics[keepaspectratio, height=5cm]{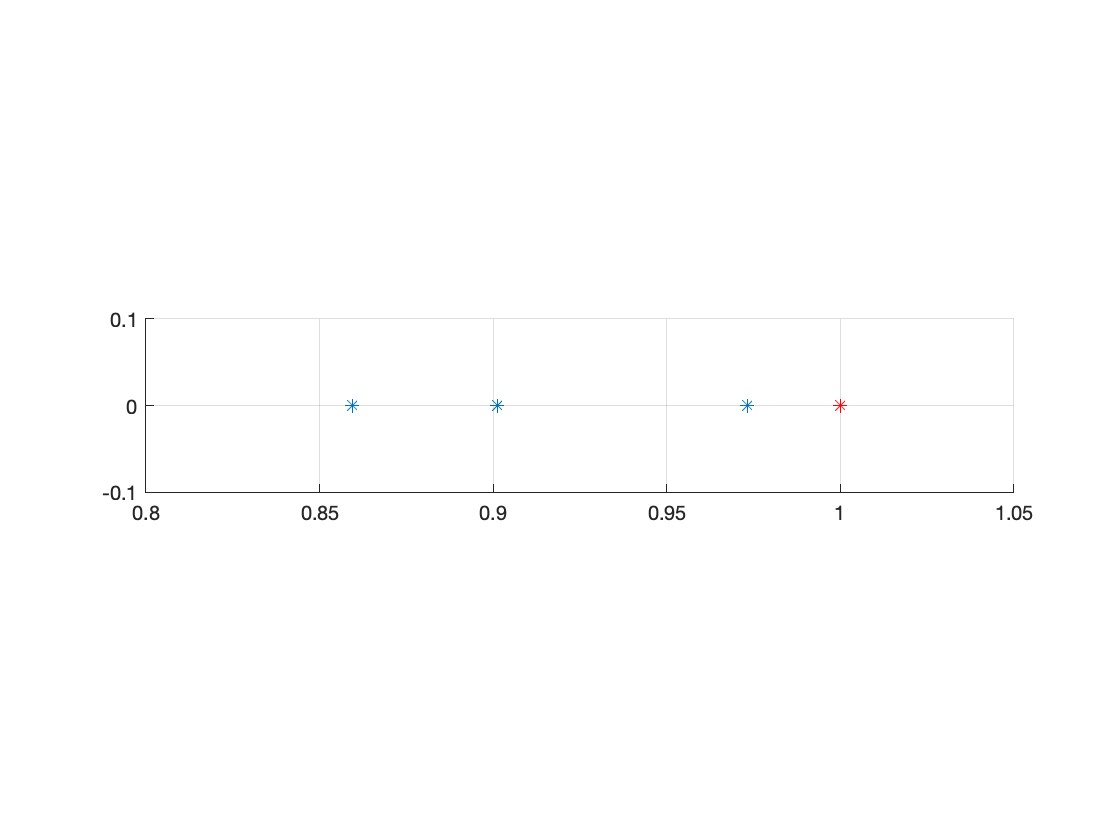}\\
					\vspace{-3.5mm}
					\centering
					\text{b. Approximated eigenvalues near $1$}
					%	\vspace{5mm}
				\end{minipage}
			\end{tabular}
			\begin{tabular}{cc}
				\begin{minipage}[t]{0.50\hsize}
					\centering
					\includegraphics[keepaspectratio, height=5cm]{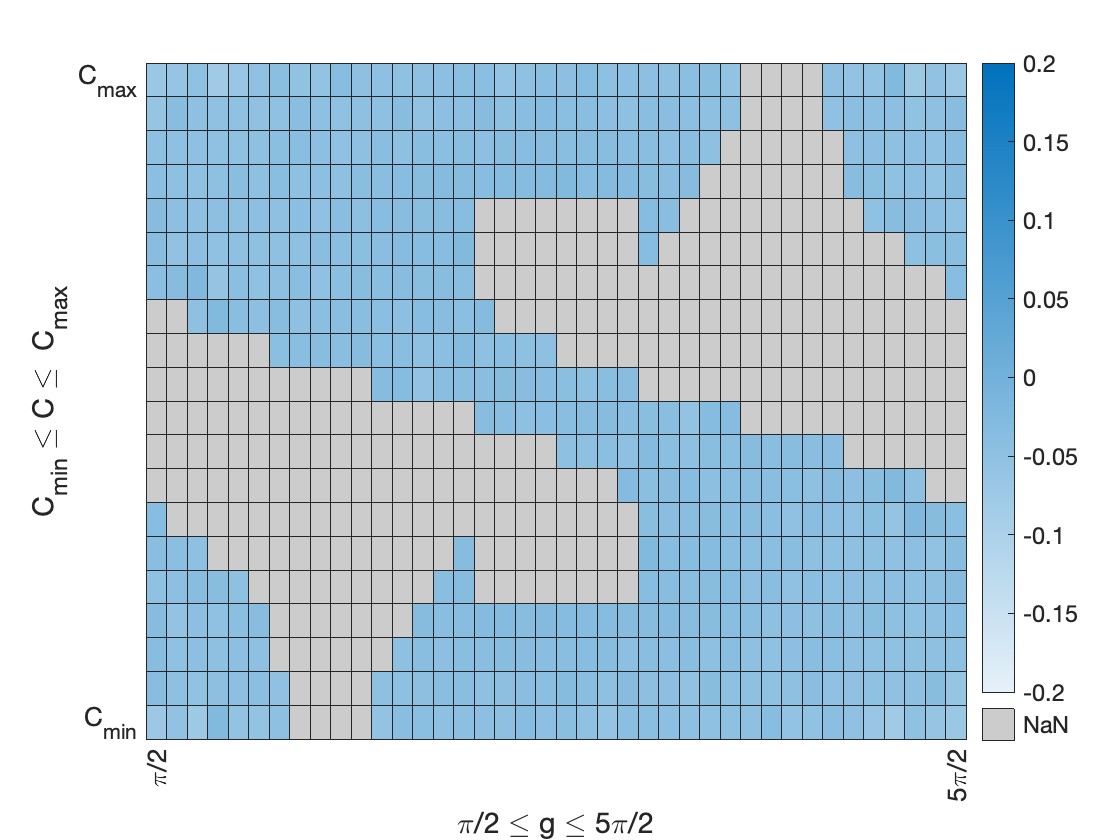}\\
					\centering
					\text{c. Eigenfunction for eigenvalue 1.00}
					\includegraphics[keepaspectratio, height=5cm]{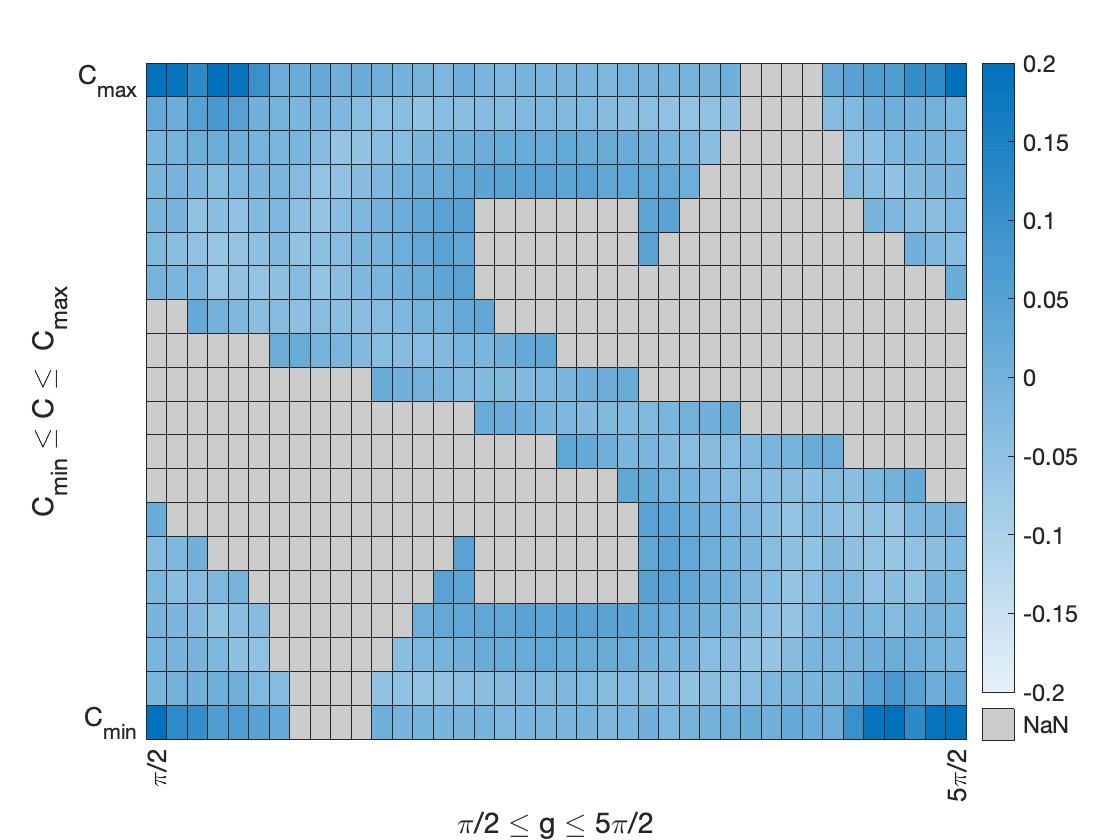}\\
					\centering
					\text{e. Eigenfunction for eigenvalue 0.90}
				\end{minipage} &
				\begin{minipage}[t]{0.50\hsize}
					\centering
					\includegraphics[keepaspectratio, height=5cm]{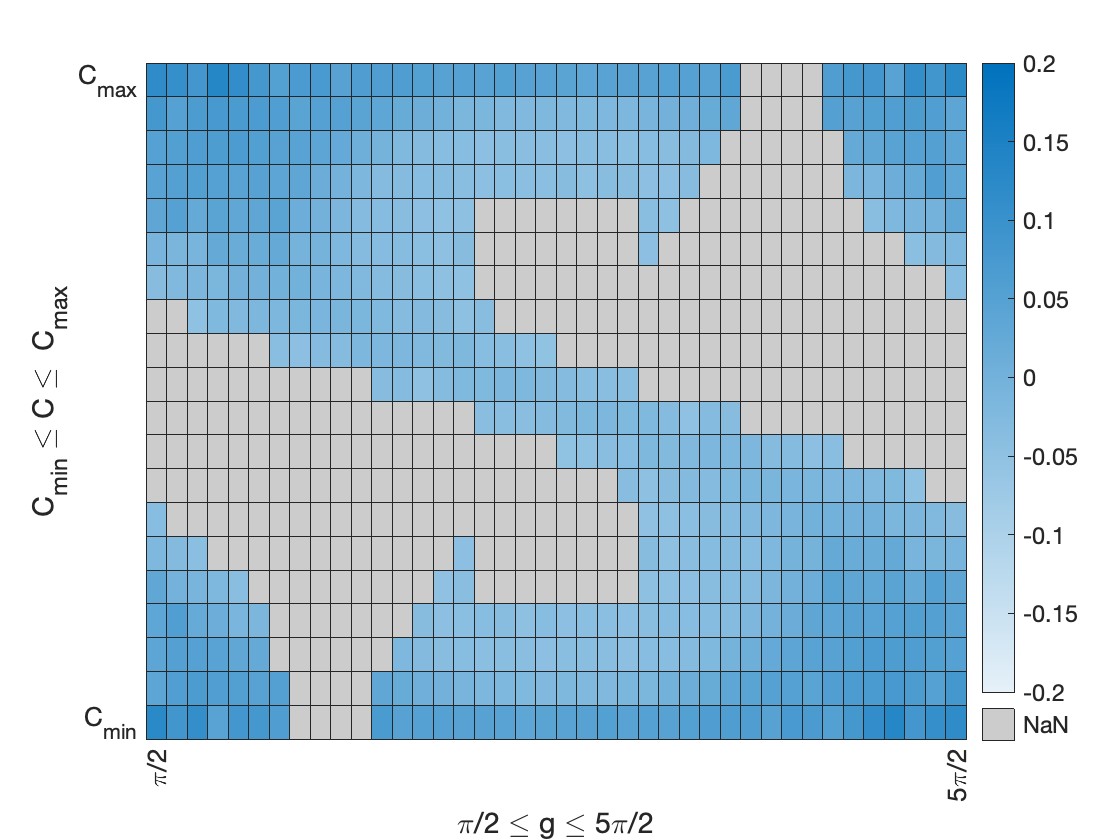}\\
					\centering
					\text{d. Eigenfunction for eigenvalue 0.97}
				\end{minipage}
			\end{tabular}
			\vspace{5mm}
			\centering
			%\text{c. level sets of independent eigenvectors of eigenvalue 1}
			\caption{Allowed regions in $(g,C)-$coordinates space $[0, 2 \pi ) \times [C_{min}, C_{max}]$.  Eigenvalues and eigenfunctions of approximated Koopman operator for $\alpha = 4.0, \beta = 0.5, N = 800, L = 25$, Gauss-Legendre quadrature.}
			\label{fig:eigen_gl_beta050_gl}
		\end{figure}

		\begin{figure}
			\centering
			\begin{tabular}{cc}
				\begin{minipage}[t]{0.5\hsize}
					\includegraphics[keepaspectratio, height=5cm]{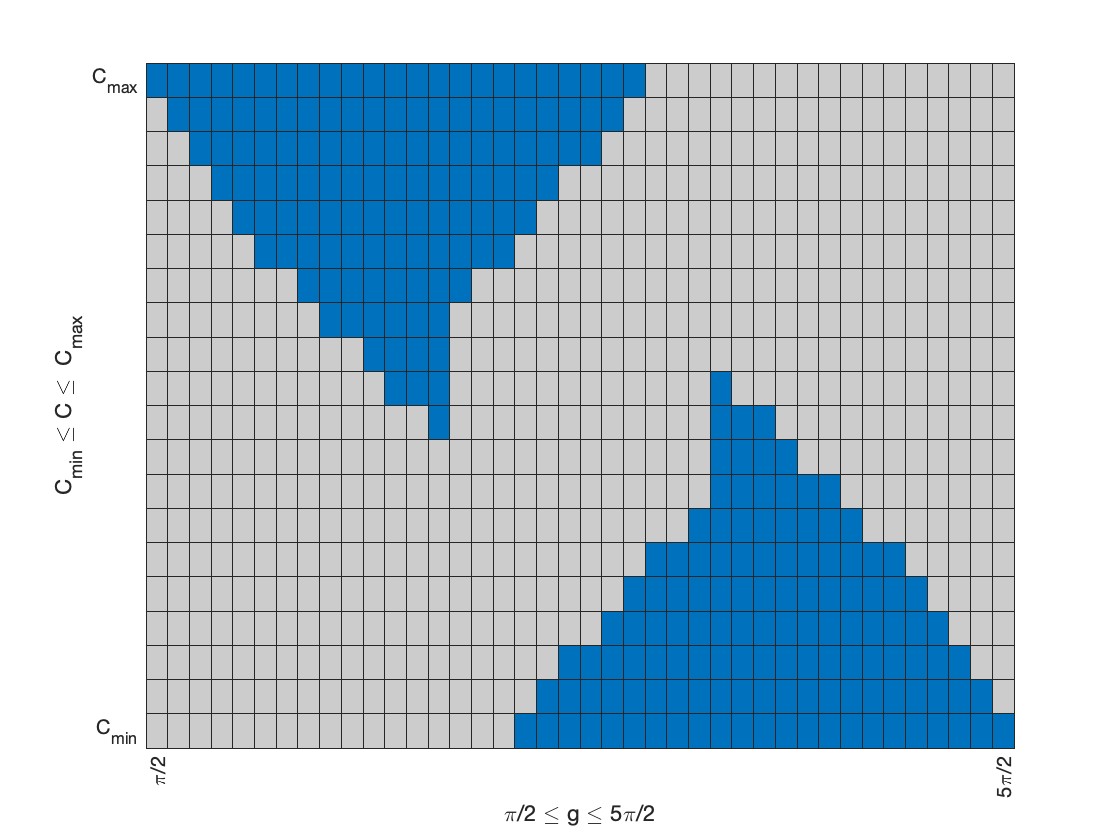}\\
					\centering
					\text{  a. Allowed regions (in blue)}
				\end{minipage}
				\begin{minipage}[t]{0.5\hsize}
					\includegraphics[keepaspectratio, height=5cm]{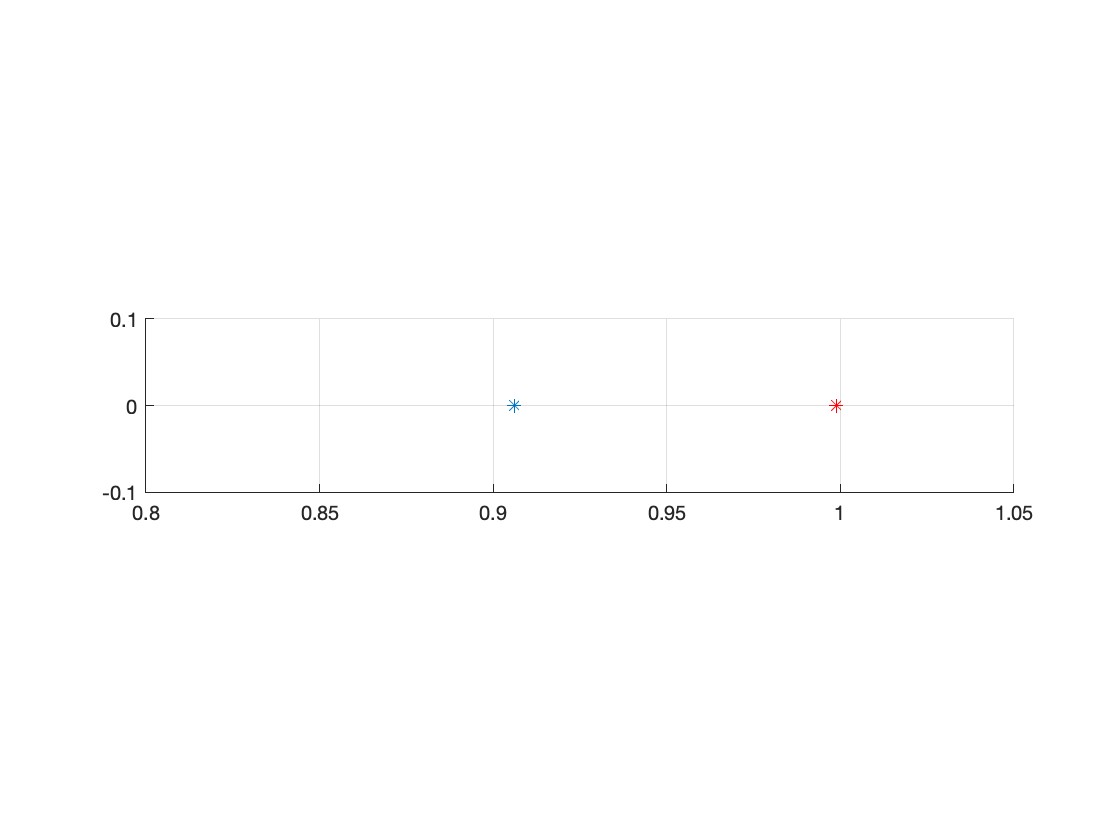}\\
					\vspace{-3.5mm}
					\centering
					\text{b. Approximated eigenvalues near $1$}
					%	\vspace{5mm}
				\end{minipage}
			\end{tabular}
			\begin{tabular}{cc}
				\begin{minipage}[t]{0.50\hsize}
					\centering
					\includegraphics[keepaspectratio, height=5cm]{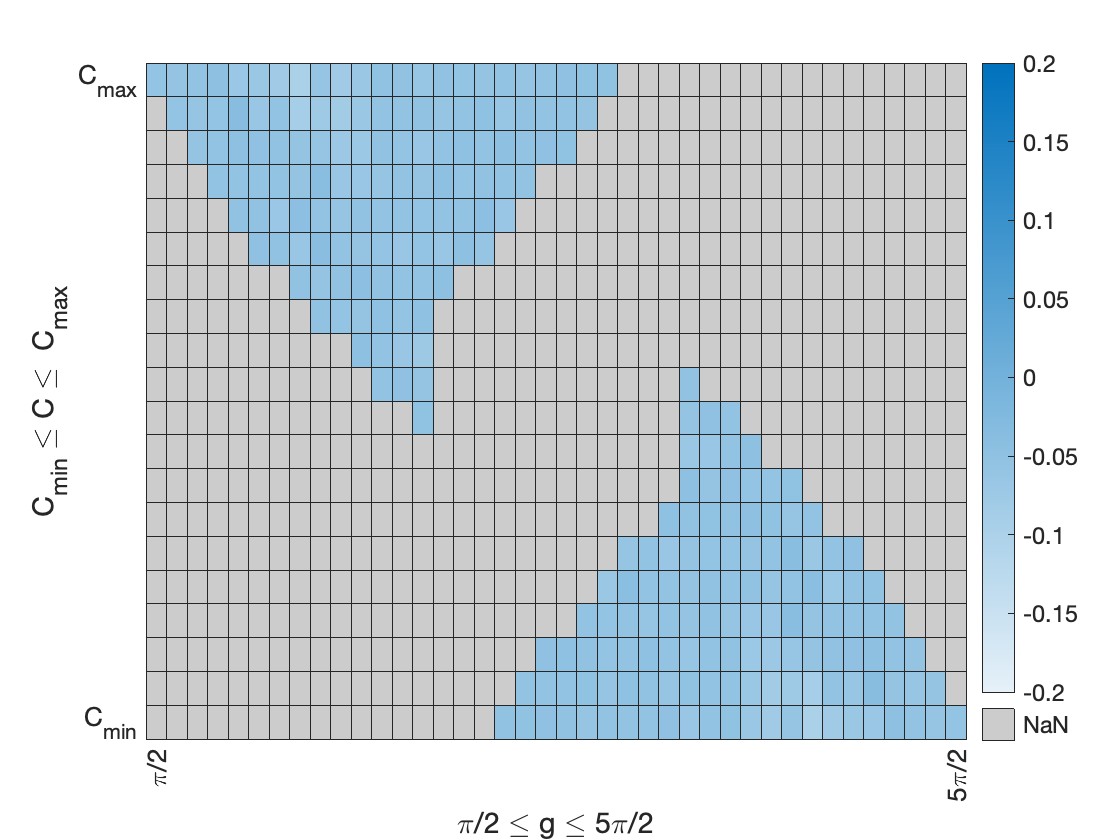}\\
					\centering
					\text{c. Eigenfunction for eigenvalue 1.00}
				\end{minipage} &
				\begin{minipage}[t]{0.50\hsize}
					\centering
					\includegraphics[keepaspectratio, height=5cm]{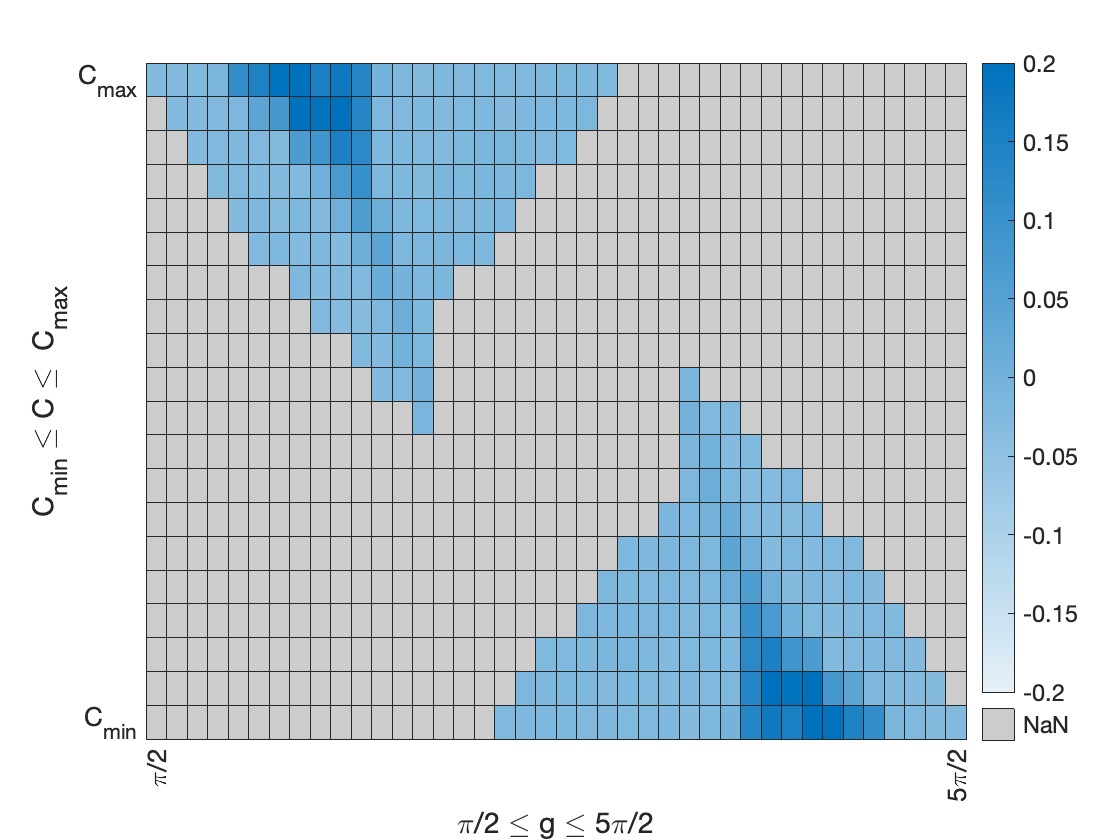}\\
					\centering
					\text{d. Eigenfunction for eigenvalue 0.91}
				\end{minipage}
			\end{tabular}
			\vspace{5mm}
			\centering
			%\text{c. level sets of independent eigenvectors of eigenvalue 1}
			\caption{Allowed regions in $(g,C)-$coordinates space $[0, 2 \pi ) \times [C_{min}, C_{max}]$.  Eigenvalues and eigenfunctions of approximated Koopman operator for $\alpha = 4.0, \beta = 2.4, N = 800, L = 25$, Gauss-Legendre quadrature.}
			\label{fig:eigen_gl_beta240_gl}
		\end{figure}

		\begin{figure}
			\centering
			\begin{tabular}{cc}
				\begin{minipage}[t]{0.5\hsize}
					\includegraphics[keepaspectratio, height=5cm]{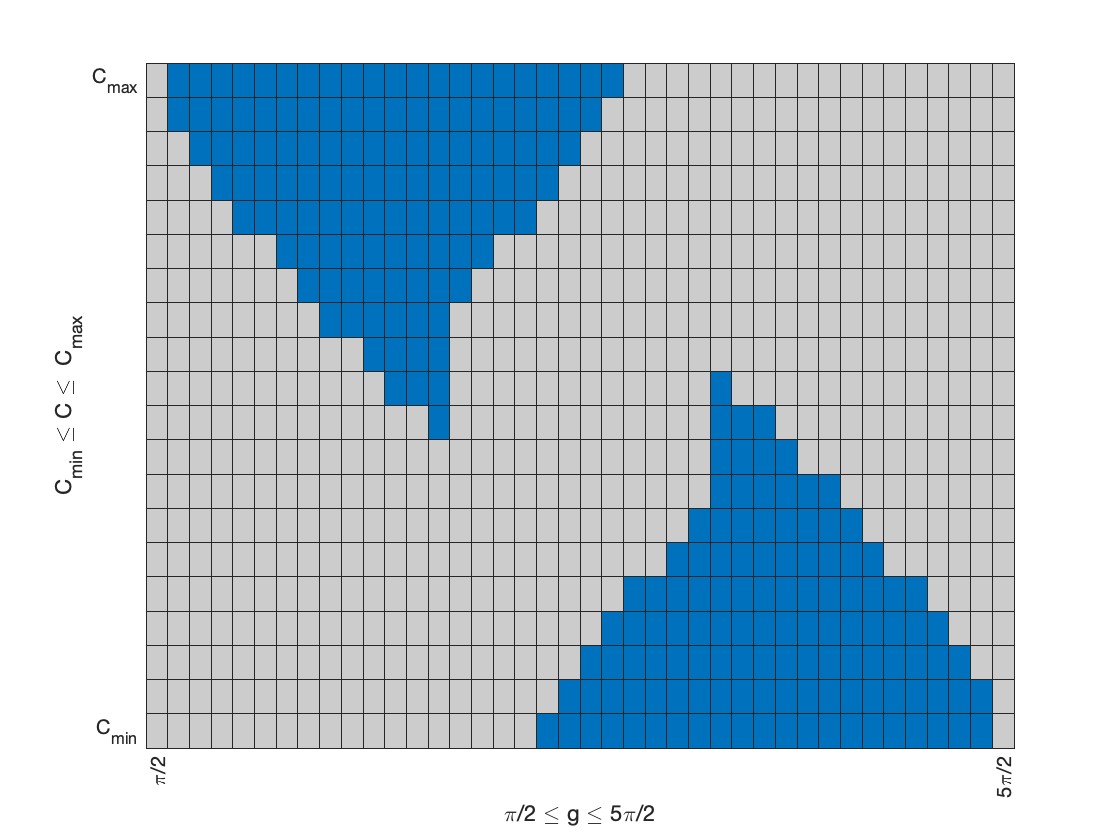}\\
					\centering
					\text{ \hspace{-10mm} a. Allowed regions (in blue) }
					\centering
					\includegraphics[keepaspectratio, height=5cm]{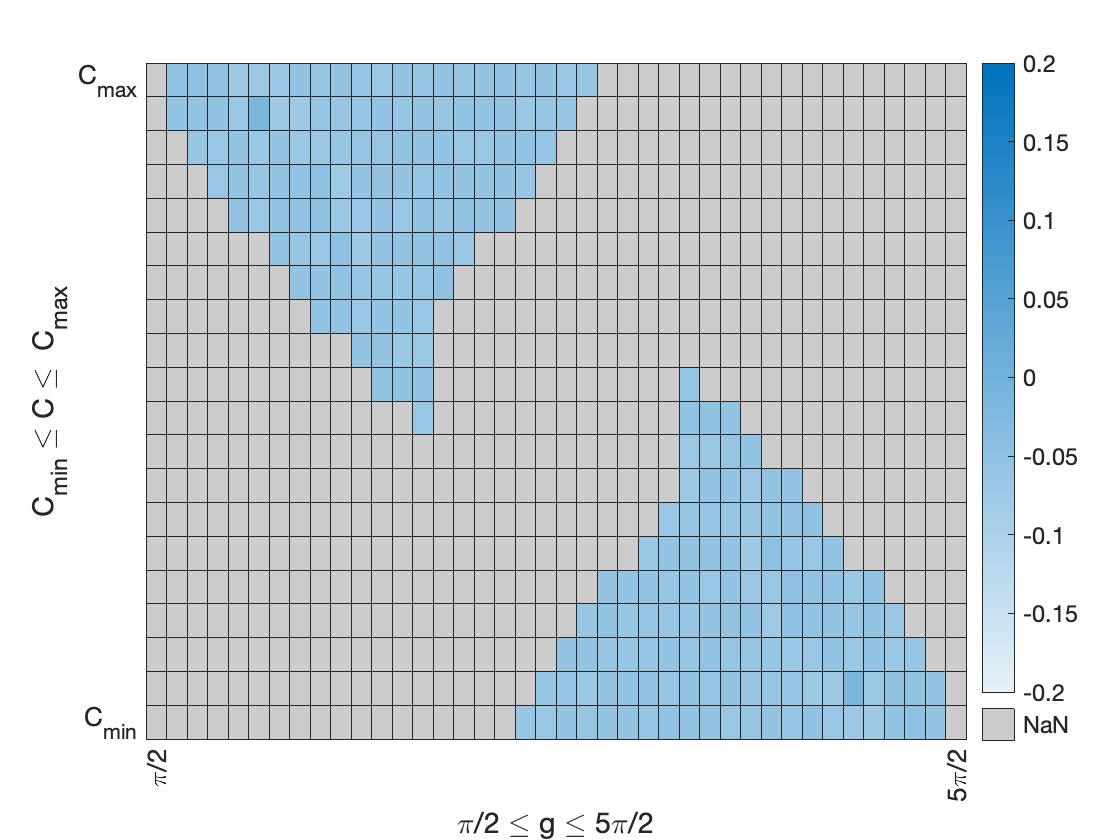}\\
					\centering
					\text{c. Eigenfunction for eigenvalue 1.00}
				\end{minipage}
				\begin{minipage}[t]{0.5\hsize}
					\includegraphics[keepaspectratio, height=5cm]{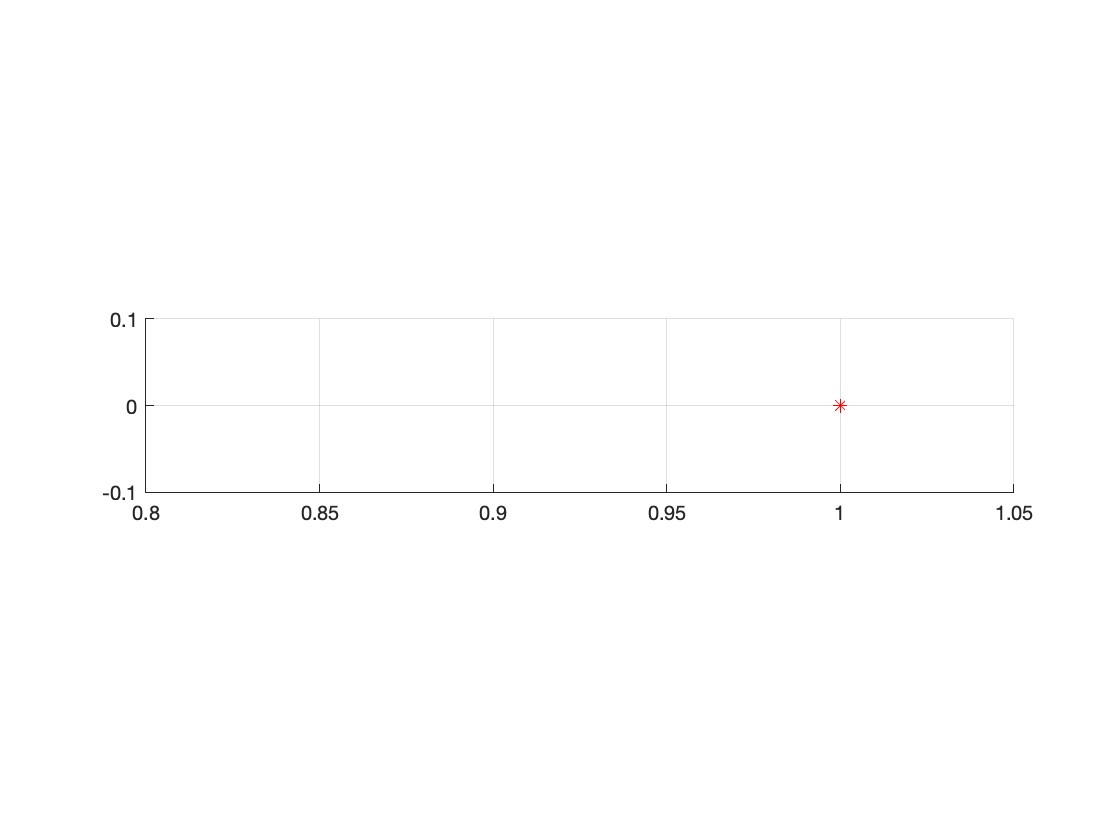}\\
					\vspace{-3.5mm}
					\centering
					\text{b. Approximated eigenvalues near $1$}
					%	\vspace{5mm}
				\end{minipage}
			\end{tabular}
			\vspace{5mm}
			\centering
			%\text{c. level sets of independent eigenvectors of eigenvalue 1}
			\caption{Allowed regions in $(g,C)-$coordinates space $[0, 2 \pi ) \times [C_{min}, C_{max}]$.   Eigenvalues and eigenfunctions of approximated Koopman operator for $\alpha = 4.0, \beta = 2.6, N = 800, L = 25$, Gauss-Legendre quadrature.}
			\label{fig:eigen_gl_beta260_gl}
		\end{figure}

			\subsection{{Discussions on the Numerical Results}}
			%Based on the numerical results we presented, in this section, we discuss ergodicity and other dynamical properties of Boltzmann's billiard systems.
			{The numerical results we have presented do not provide rigorous proofs, as the  Galerkin approximation might not be able to capture the true eigenfunctions corresponding to eigenvalue 1 with high oscillation terms. However they suggest what the true dynamics of the corresponding systems could be. }
			
			For the Kepler case ($\beta = 0$), {presented in} Figure \ref{fig:eigen_beta000_gl} and Figure \ref{fig:eigen_gl_beta000_gl}, our numerical study indicates that there is {a} large multiplicity for the eigenvalue 1. Also, these figures indicate that the level sets of {the eigenfunctions with eigenvalue (at least close to) 1 are invariant subsets consisting of periodic trajectories}. These results are compatible with the integrability of the billiard system for $\beta = 0$, as it has been shown in \cite{Gallavotti-Jauslin}. 
			
			For small values of $\beta$, our numerical results (Figure \ref{fig:eigen_beta050_gl} and Figure \ref{fig:eigen_gl_beta050_gl}) indicate that the there is still a large multiplicity for the eigenvalue 1 and the level sets of its eigenfunctions show many invariant subsets of the system. {The system is unlikely to be ergodic}. This {is in consistence }with the KAM stability of the integrable Boltzmann's billiard system ($\beta=0$) under the small perturbation by the additional centrifugal force $\beta/r^2$ {\cite{Felder}}. 
			% in a force function with $\beta \simeq 0$, which is again compatible with the KAM applicability shown in \cite{Felder}.
			
			For large values of $\beta$, {various types of dynamics may coexist. }As one can see from the level sets of eigenfunction depicted in Figure \ref{fig:eigen_beta240_gl} and Figure \ref{fig:eigen_gl_beta240_gl}, for $\beta = 2.4$, there exists small regions which are foliated by (quasi-)periodic trajectories and the left region is a large indecomposable invariant subset which is covered by a single chaotic trajectory. Our particular interest is the case $\beta = 2.6$, {presented in } Figure \ref{fig:eigen_beta260_gl} and Figure \ref{fig:eigen_gl_beta260_gl}, {in which} the discretized eigenvalue problem {seems to have}  only one simple eigenvalue in the neighborhood of 1, indicating the {potential} ergodicity of the system.

			\vspace{0.5cm}
			\textbf{Acknowledgement} A.T. and L.Z. are supported by DFG ZH 605/1-1, ZH
			605/1-2.
			
			%On the other hand, this Galerkin approximation might not be able to capture the true eigenfunctions corresponding to eigenvalue 1 with high oscillation terms. This means that the system might have invariant subsets which are finely distributed in the phase space, which indeed also implies the system is \emph{nearly} ergodic. 
			
			%Based on the above observations, it seems possible that for some parameter values of $\alpha, \beta, \gamma$, the corresponding Boltzmann's billiard system is ergodic.			

\hspace{-1cm}
\begin{tabular}{@{}l@{}}%
	Michael Plum\\
	\textsc{Karlsruhe Institute of Technology, Germany.}\\
	\textit{E-mail address}: \texttt{michael.plum@kit.edu}
\end{tabular}
\vspace{10pt}

\hspace{-1cm}
\begin{tabular}{@{}l@{}}%
	Airi Takeuchi\\
	\textsc{University of Augsburg, Augsburg, Germany.}\\
	\textit{E-mail address}: \texttt{airi1.takeuchi@uni-a.de}
\end{tabular}
\vspace{10pt}

\hspace{-1cm}
\begin{tabular}{@{}l@{}}%
	Lei Zhao\\
	\textsc{University of Augsburg, Augsburg, Germany.}\\
	\textit{E-mail address}: \texttt{lei.zhao@math.uni-augsburg.de}
\end{tabular}

\end{document}